\documentclass[10pt,a4paper]{article}

\usepackage[hidelinks]{hyperref}
\usepackage{lineno}

\usepackage{amssymb, bm}
\usepackage{amsmath}
\usepackage{amsthm}
\usepackage{mathrsfs}
\usepackage{siunitx}
\usepackage[margin=2cm]{geometry}
\usepackage{subfig}
\usepackage{float}
\usepackage{xcolor}
\usepackage{subfloat}

\usepackage{mathtools}
\usepackage{bbm}
\usepackage{multirow}
\usepackage{stmaryrd}

\usepackage{algorithm}     
\usepackage{algpseudocode} 

\usepackage[utf8]{inputenc} 
\usepackage[T1]{fontenc}

\usepackage{tikz}
\usetikzlibrary{decorations.markings}
\usetikzlibrary{shapes.geometric}
\usetikzlibrary{shapes.arrows}
\usetikzlibrary{fit,calc}
\usetikzlibrary{quotes,angles}
\usetikzlibrary{arrows.meta}
\usepackage{animate}  
\usepackage{graphicx}
\usepackage{pgf}

\usepackage{pgfplots}
\usepackage[normalem]{ulem}

\newcommand{\vertiii}[1]{{\left\vert\kern-0.25ex\left\vert\kern-0.25ex\left\vert #1 
   \right\vert\kern-0.25ex\right\vert\kern-0.25ex\right\vert}}

\usepackage{geometry}
\makeatletter
\def\blfootnote{\xdef\@thefnmark{$\star$}\@footnotetext}
\makeatother
\newenvironment{Authors}%
  {\begin{center}\begin{bfseries}}%
  {\end{bfseries}\end{center}}
\newenvironment{Addresses}%
  {\begin{flushleft}\begin{itshape}}%
  {\end{itshape}\end{flushleft}}
  \newcommand{\email}[1]{\hspace*{\stretch{1}}\emph{\texttt{#1}}}

 \usepackage{fancyhdr}  
  \fancypagestyle{plain}{
\fancyhead{}
\fancyhead[C]{\hfill Submitted to Elsevier, June 2026}
 }

\newtheorem*{remark*}{Remark}  

\begin{document}

 \title{A domain-decomposition reduced-order model for two-phase CO$_2$ sequestration with moving wells}

\date{}

 \maketitle
\vspace{-50pt} 

\begin{Authors}
Xin He$^{1}$,
Xianda Shen$^{2,3}$,
Lei Zhang$^1$.
\end{Authors}

\begin{Addresses}
$^1$
School of Mathematical Sciences, Key Laboratory of Intelligent Computing and Applications (Ministry of Education), Tongji University, Shanghai 200092, China, \\
$^2$
Department of Geotechnical Engineering, College of Civil Engineering, Tongji University, Shanghai 200092, China, \\
$^3$
Key Laboratory of Geotechnical and Underground Engineering of Ministry of Education, Tongji University, Shanghai 200092, China, 
\\ \email{2433969@tongji.edu.cn, xshen@tongji.edu.cn, 22210@tongji.edu.cn} \\ \email{} \\[3mm]

\end{Addresses}

\begin{abstract}
We present a domain-decomposition reduced-order model for two-phase geological CO$_2$ storage with fixed and moving wells. The approach builds on a nonoverlapping Schur formulation of a conservative RT0--P0 discretization solved by a sequential pressure--transport algorithm. The distinctive feature of the formulation is that local reduced models are coupled through physically meaningful interface variables: pressure traces for the pressure solve and inflow saturation traces for the transport solve. To address moving wells, we introduce a well-aware strategy that combines source-response enrichment, HAPOD-based far-field compression, and hybrid full-order/reduced-order treatment of well and near-well subdomains. We also formulate conservation-oriented hyper-reduction rules where interface fluxes, boundary terms and source terms are protected.  Numerical results are presented to demonstrate accuracy, computational efficiency, robustness of the overall methodology. 

\end{abstract}


\section{Introduction}
\label{sec:introduction}

Geological CO$_2$ storage is an important technology for reducing
anthropogenic carbon emissions and supporting long-term climate-change mitigation strategies \cite{bui2018ccs,bensonCole2008,bachu2008storage}. In a storage operation, CO$_2$ is
injected into deep geological formations, where it migrates as a
multiphase plume and interacts with resident brine, heterogeneous rock properties, capillary forces, pressure constraints, and injection controls. Numerical simulation is therefore central to storage-capacity assessment, pressure-management analysis, plume-migration prediction, leakage-risk evaluation, and monitoring design \cite{pruessGarcia2002,class2009benchmark}. 
However, high-fidelity simulations, such as finite element, finite volume methods, are computationally expensive,
especially when they needs to be repeatedly performed in parameter studies, optimization, uncertainty quantification, and well-control design \cite{jinDurlofsky2018,feng2024encoder,
pruessGarcia2002,class2009benchmark}.

Reduced-order models (ROMs) and surrogate models provide a natural route
to accelerate repeated CO$_2$ storage simulations. Projection-based methods, such as proper orthogonal decomposition (POD) and reduced basis (RB) methods, have demonstrated significant computational savings by approximating the high-dimensional discrete solutions within low-dimensional subspaces \cite{benner2015survey,quarteroni2016reduced}. For nonlinear time-discrete residuals,
least-squares Petrov--Galerkin (LSPG) and GNAT-type methods improve robustness by computing reduced states through residual minimization 
\cite{carlberg2013gnat,carlberg2017galerkin}. 
In reservoir and CO\(_2\) storage applications, projection-based ROMs have been used to accelerate repeated storage-operation simulations \cite{jinDurlofsky2018}. Data-driven surrogates further extend this direction to reservoir simulations with time-varying well controls \cite{jin2020deep}, including encoder--decoder ConvLSTM architectures for geological CO\(_2\) sequestration under dynamic well-control settings \cite{feng2024encoder}.

Projection-based MOR reduces the dimension of the state variables, but the computational complexity of the nonlinear residual and Jacobian evaluations often remains tied to the dimension of the underlying high-fidelity model. Hyper-reduction techniques address this bottleneck by approximating nonlinear operators using a reduced amount of high-fidelity information. For example, DEIM approximates nonlinear functions in a low-dimensional collateral basis and evaluates them only at selected interpolation points \cite{chaturantabut2010deim}. GNAT adopts a related residual-based viewpoint, combining residual sampling and reconstruction with least-squares Petrov--Galerkin projection \cite{carlberg2013gnat}. Sampling-and-weighting approaches such as ECSW \cite{farhat2015ecsw} and empirical cubature \cite{hernandez2017empiricalCubature} instead construct reduced quadrature rules that approximate nonlinear residual and Jacobian contributions through weighted subsets of elements or integration points. More recently, component-based hyper-reduction strategies have been developed for reduced-basis-element formulations, enabling localized nonlinear approximation within domain-decomposition ROM frameworks \cite{ebrahimi2024hyperreduced}. For conservative porous-media flow, hyper-reduction must be designed carefully. If interface fluxes, physical boundary fluxes, well/source terms, or subdomain-integrated mass balances are sampled inconsistently, the reduced model may introduce artificial sources or sinks. Such errors are particularly problematic in long-time CO$_2$ storage simulations, where even small conservation violations may accumulate and distort plume migration and storage accounting. 

Despite these advances, 
standard global ROMs are most reliable when the online solution remains close to the offline snapshot manifold used for training. This assumption becomes restrictive when injection wells or source locations change. A basis trained for one well
configuration may still represent the far-field pressure response, but the near-well pressure buildup, source-induced flux pattern, saturation front, and interface traces generated by a new well location may be poorly represented. 
Predict-and-correct strategies for changing well locations have been proposed in this context \cite{zalavadia2021two}, highlighting the need for reduced models that can handle local modifications without global retraining. In addition, a reduced model for CO$_2$ storage must
not only approximate the solution fields, but also preserve flux
continuity and component-mass balance, since small artificial sources or
sinks can accumulate over long-time plume simulations.

Domain decomposition (DD) provides a natural framework for addressing these limitations. 
Classical DD approaches, including transmission conditions \cite{alboin2000domain} and mortar-based methods \cite{ganis2014global}, have been extensively applied to porous-media flow problems;
multiscale and high-contrast problems have motivated reduced-dimensional coarse spaces and interface-aware preconditioners \cite{galvisEfendiev2010};
 nonlinear and adaptive DD methods have also been developed to localize challenging computations, exploit spatially localized nonlinearities, and enhance solver scalability  \cite{skogestad2013domain,cheng2024adaptive}.
Beyond its solver capabilities, DD provides a natural component-based architecture for ROM: local reduced spaces can approximate each subdomain independently, while interface unknowns maintain global coupling. This idea is closely related to localized reduced-basis methods for porous-media flow \cite{kaulmann2015lrbms,ohlbergerSchindler2015}, arbitrary local modification and randomized local ROM strategies \cite{buhrEngwerOhlbergerRave2017,buhrSmetana2018}, port-reduced static-condensation methods \cite{smetana2015certification,smetanaPatera2016}, and DD-LSPG formulations for nonlinear model reduction \cite{hoang2021domain}. 
These component-based constructions are particularly attractive for changing well locations: only the well and near-well subdomains may require enrichment or
high-fidelity treatment, whereas far-field subdomains can still be
represented by reduced spaces learned from previous simulations.
Even when the well location is fixed, combining parametric MOR with DD is useful for heterogeneous CO$_2$ storage problems, because permeability variations, localized sources, and parameter-dependent flow responses often have strong spatial structure. Reducing subdomain interiors separately while coupling them through interface variables also facilitates localized basis construction, reuse of far-field information, and parallel implementation.

This paper develops a domain-decomposition reduced-order model (DD-ROM) for two-phase CO$_2$ storage simulations, capable of handling both fixed and moving wells. The underlying high-fidelity discretization uses the RT0--P0 mixed finite-element method~\cite{raviartThomas1977,chenEwing1997mixed,
hoteitFiroozabadi2006,radu2018robust,caoMitra2019}, providing explicit normal flux representation and local mass conservation essential for accurate porous-media flow simulations~\cite{chenEwing1997mixed,hoteitFiroozabadi2006}. The resulting nonlinear algebraic problem is solved using a sequential fully implicit (SFI) approach~\cite{moncorge2018sfi,moncorge2020upwinding,jiangTchelepi2019}, which separates the coupled problem into pressure/flow and transport subproblems.  
This splitting is advantageous for model reduction, as it yields subproblems with distinct mathematical structures,
interface variables, and reduced-basis requirements. 

The main contributions of this work are as follows: 
we develop a DD-ROM framework for two-phase CO$_2$
storage based on pressure and transport Schur complement formulations. The method combines local POD spaces, interface-pressure lifting, Galerkin/LSPG projection, and either full-interface or port-reduced coupling while preserving the conservative structure of the underlying high-fidelity discretization;
we introduce a well-aware model reduction strategy for changing
well locations. The approach combines source-response enrichment,
HAPOD-based compression \cite{himpe2018hapod}, and hybrid FOM--ROM coupling to capture
localized near-well dynamics while maintaining reduced-order efficiency in the far field;
we propose a hyper-reduction framework in
which only interior residual contributions are sampled, while interface fluxes, boundary terms, well/source contributions, and selected mass-balance constraints are treated as protected quantities. This ensures compatibility between hyper-reduction and conservative porous media flow simulation.

The remainder of the paper is organized as follows.
Section~\ref{sec:model} presents the governing equations and the
RT0--P0 high-fidelity discretization. Section~\ref{sec:sfi-dd} introduces the
sequential pressure/transport splitting and its DD Schur formulation.
Section~\ref{sec:ddrom-compact} presents the pressure and transport
DD-ROMs. Section~\ref{sec:well-aware-compact} describes the well-aware
basis construction, hybrid FOM--ROM strategy, and
hyper-reduction strategy. Section~\ref{sec:numerical} reports numerical results,
and Section~\ref{sec:conclusion} concludes the paper.

\section{Governing equations and high-fidelity discretization}
\label{sec:model}

We consider an isothermal two-phase, two-component flow model for CO$_2$ sequestration \cite{pruessGarcia2002,class2009benchmark,
chenEwing1997formulations,hoteitFiroozabadi2006,
chen2006computational}.
The two phases are the aqueous phase and the gas phase, denoted by $\beta\in\{a,g\}$, and the two components are water and CO$_2$, denoted by $\kappa\in\{\mathrm{H}_2\mathrm{O},\mathrm{CO}_2\}$. 
The governing equations are expressed in terms of total Darcy velocity, gas pressure, and gas saturation $(\mathbf{u}_t,p,S_g)$. This choice naturally facilitates the residual splitting used in the sequential fully implicit (SFI)
method presented in Section~\ref{sec:sfi-dd}. Specifically, with saturation frozen, the Darcy equation and
the total-mass balance define a pressure/flow subsystem for $(\mathbf{u}_t,p)$, while the
CO$_2$ component balance gives a transport subsystem for $S_g$ once $(\mathbf{u}_t,p)$
has been updated \cite{moncorge2018sfi,moncorge2020upwinding,jiangTchelepi2019}.

\subsection{Governing Equations}
Let $\Omega\subset\mathbb{R}^d$, $d=2,3$, be a bounded porous-medium domain with boundary $\Gamma=\partial\Omega$ and outward unit normal $\mathbf{n}$. The boundary is decomposed as $\Gamma=\Gamma_D\cup\Gamma_N$, where pressure is prescribed on $\Gamma_D$ and normal total flux is prescribed on $\Gamma_N$. 
For each component $\kappa$, mass conservation in the porous medium is
\begin{equation}
\frac{\partial}{\partial t}\!\left(\phi \sum_{\beta} S_\beta \rho_\beta X_\beta^\kappa\right)
+\nabla\!\cdot\!\left(\sum_{\beta} \rho_\beta X_\beta^\kappa\,\mathbf{u}_\beta\right)
= Q^\kappa,
\label{eq:mass_balance_component}
\end{equation}
where $\phi$ is the porosity, $S_{\beta}$ is the phase saturation, $\rho_{\beta}$ is the phase density, $X_{\beta}^{\kappa}$ is the mass fraction of component $\kappa$ in phase $\beta$, 
$\mathbf{u}_{\beta}$ is the phase Darcy velocity, 
and $Q^\kappa$ is the component source term.
The pore space is assumed fully saturated,
\begin{equation*}
    S_a+S_g=1,
\end{equation*}
and the mass fractions in each phase sum to one:
\begin{equation*}
    \sum_{\kappa} X_\beta^{\kappa}=1.
\end{equation*}
The phase Darcy velocities follow Darcy's law
\begin{equation}
\mathbf{u}_\beta
= -\,\mathbf{K}\,\lambda_\beta\Big(\nabla p+\nabla p_{c\beta}-\rho_\beta \mathbf{g}\Big), 
\label{eq:darcy_flux_component}
\end{equation}
where $\mathbf{K}$ is the absolute permeability, 
$\mathbf{g}$ is the gravity vector, and
$$\lambda_{\beta}=k_{r\beta}/\mu_{\beta}$$ is the phase mobility, with relative permeability $k_{r\beta}$ and viscosity $\mu_\beta$. We take the gas pressure as the reference pressure, and the capillary pressures satisfy $p_{cg}\equiv0$ and $p_c:=p_{ca}\leq0$. 

To derive a total-velocity formulation advantageous for the sequential solution strategy presented in Section~\ref{sec:sfi-dd}, we first define total mobility and fractional-flow functions:
\begin{equation*}
    \lambda_t=\lambda_a+\lambda_g,\qquad
    f_\beta=\frac{\lambda_\beta}{\lambda_t},\quad \beta\in\{a,g\}.
    \label{eq:mobility-fractional-flow}
\end{equation*}
Summing the two phase Darcy equations (\ref{eq:darcy_flux_component}), and defining the total velocity as $\mathbf{u}_t=\mathbf{u}_a+\mathbf{u}_g$, we obtain the Darcy equation for total velocity
 \cite{chenEwing1997formulations,hoteit2008numerical}:
\begin{equation*}
\mathbf{u}_t
= -\,\mathbf{K}\lambda_t\,\nabla p
- \mathbf{K}\,\lambda_a\,\nabla p_c
+ \mathbf{K}\,\big(\lambda_a\rho_a+\lambda_g\rho_g\big)\,\mathbf{g}.
\label{eq:darcy_flux_total}
\end{equation*}
In this work, we neglect the gravitational effect $\mathbf{g}=\mathbf{0}$. 
We then rewrite the mass conservation equations \eqref{eq:mass_balance_component} by explicitly expressing the phase velocities in terms of the total velocity $\mathbf{u}_t$. Moreover, we sum these mass conservation equations over the two components to obtain the total-mass conservation equation, and retain the CO$_2$-component mass conservation separately to ensure completeness of the compositional description. We finally obtain the following governing system for the unknowns $(\mathbf{u}_t,p,S_g)$ considered in this work:
\begin{equation}
\left\{
\begin{aligned}
(\lambda_t(S_g)\mathbf{K})^{-1}\mathbf{u}_t
+\nabla p+f_a(S_g)\nabla p_c(S_g)
&=\mathbf{0}
&&\text{in }\Omega\times(0,T],\\
\partial_t M_t(p,S_g)+\nabla\!\cdot\mathbf{m}_t(\mathbf{u}_t,p,S_g)
&=Q_t
&&\text{in }\Omega\times(0,T],\\
\partial_t M_{\mathrm{CO}_2}(p,S_g)
+\nabla\!\cdot\mathbf{m}_{\mathrm{CO}_2}(\mathbf{u}_t,p,S_g)
&=Q_{\mathrm{CO}_2}
&&\text{in }\Omega\times(0,T],
\end{aligned}
\right.
\label{eq:continuous-system}
\end{equation}
where 
the mass accumulations are
\begin{equation}
\begin{aligned}
M_t(p,S_g)&=\phi\sum_{\beta\in\{a,g\}}S_\beta\rho_\beta,\\
M_{\mathrm{CO}_2}(p,S_g)&=\phi\sum_{\beta\in\{a,g\}}S_\beta\rho_\beta X_\beta^{\mathrm{CO}_2},
\end{aligned}
\label{eq:accumulation-definitions}
\end{equation}
the corresponding fluxes, written in terms of $(\mathbf{u}_t,p,S_g)$, are
\begin{equation*}
\mathbf{m}_t(\mathbf{u}_t,p,S_g)
= (\rho_a f_a+\rho_g f_g)\,\mathbf{u}_t
+\mathbf{K}\frac{\lambda_a\lambda_g}{\lambda_t}(\rho_g-\rho_a)\nabla p_c,
\label{eq:total-mass-flux}
\end{equation*}
\begin{equation*}
\mathbf{m}_{\mathrm{CO}_2}(\mathbf{u}_t,p,S_g)
= (\rho_a f_a X_a^{\mathrm{CO}_2}+\rho_g f_g X_g^{\mathrm{CO}_2})\,\mathbf{u}_t
+\mathbf{K}\frac{\lambda_a\lambda_g}{\lambda_t}
\big(\rho_g X_g^{\mathrm{CO}_2}-\rho_a X_a^{\mathrm{CO}_2}\big)\nabla p_c .
\label{eq:co2-mass-flux}
\end{equation*}
Here, the relative permeabilities and capillary pressure are assumed to depend only on saturation, following Corey-type or van Genuchten-type closures \cite{vanGenuchten1980,hoteit2008numerical}. 
The boundary and initial data are
\begin{equation*}
\left\{
\begin{aligned}
p&=p_D &&\text{on }\Gamma_D\times(0,T],\\
\mathbf{u}_t\cdot\mathbf{n}&=g_N &&\text{on }\Gamma_N\times(0,T],\\
S_g&=S_D &&\text{on }\Gamma_S\times(0,T],\\
(\mathbf{u}_t,p,S_g)(\cdot,0)&=(\mathbf{u}_t^0,p^0,S_g^0) &&\text{in }\Omega .
\end{aligned}
\right.
\label{eq:boundary-initial-data}
\end{equation*}
Here $\Gamma_S$ denotes the inflow part of the boundary where saturation is prescribed. On outflow boundaries, the saturation is determined by the interior/upwind state.

\subsection{Weak Formulation}
We introduce function spaces
\begin{align*}
H(\mathrm{div},\Omega)&:=\{\mathbf{v}\in L^{2}(\Omega)^{d}:\mathrm{div~}\mathbf{v}\in L^{2}(\Omega)\}, \\
H_{0,N}(\mathrm{div},\Omega)&:=\{\mathbf{v}\in H(\mathrm{div},\Omega):\mathbf{v}\cdot\mathbf{n}=0\mathrm{~on~}\Gamma_{N}\}, \\
H_{g_N,N}(\mathrm{div},\Omega)&:=\{\mathbf{v}\in H(\mathrm{div},\Omega):\mathbf{v}\cdot\mathbf{n}=g_N\mathrm{~on~}\Gamma_{N}\}.
\end{align*}
Then the weak formulation of \eqref{eq:continuous-system} reads \cite{raviartThomas1977,chenEwing1997mixed,radu2018robust,caoMitra2019}.
: for $t\in(0,T]$, seek $\mathbf{u}_t\in H_{g_N,N}(\mathrm{div},\Omega)$ and $p,S_g\in L^2(\Omega)$ such that, for all $\mathbf{v}\in H_{0,N}(\mathrm{div},\Omega)$ and $w\in L^2(\Omega)$,
\begin{equation}
\left\{
\begin{aligned}
\int_\Omega (\lambda_t(S_g)\mathbf{K})^{-1}\mathbf{u}_t\!\cdot\!\mathbf{v}
-\int_\Omega p\,\nabla\!\cdot\mathbf{v}
+ \int_\Omega \mathbf{f}_c(S_g)\!\cdot\!\mathbf{v}
&= - \int_{\Gamma_D} p_D(\mathbf{v}\!\cdot\!\mathbf{n}),\\
\int_\Omega w\,\partial_t M_t(p,S_g)
+\int_\Omega w\,\nabla\!\cdot \mathbf{m}_t(\mathbf{u}_t,p,S_g)
&= \int_\Omega w\,Q_t,\\
\int_\Omega w\,\partial_t M_{\mathrm{CO}_2}(p,S_g)
+\int_\Omega w\,\nabla\!\cdot \mathbf{m}_{\mathrm{CO}_2}(\mathbf{u}_t,p,S_g)
&= \int_\Omega w\,Q_{\mathrm{CO}_2},
\end{aligned}
\right.
\label{eq:weak-formulation}
\end{equation}
where $\mathbf{f}_c(S_g)=f_a(S_g)\nabla p_c(S_g)$.

\subsection{RT0--P0 Fully Discrete Formulation}
\label{sec:disc}
Let $\mathcal{T}_h$ be a triangular mesh of the domain $\Omega$, consisting of cells ${K\in\mathcal{T}_h}$, faces 
${e\in\mathcal{E}_h}$, with a total of $N_c$ cells and $N_f$ faces. 
The outward unit normal to a cell $K$ on a face $e\subset\partial K$ is denoted by $\mathbf{n}_{K,e}$.
We approximate the total velocity using the lowest-order Raviart--Thomas element (RT0) and represent both the pressure and gas saturation as cellwise constants (P0).
The RT0--P0 discretization is a natural mixed formulation for Darcy flows \cite{chenEwing1997mixed,hoteitFiroozabadi2006,
radu2018robust,caoMitra2019,
arbogastTao2019}. 
RT0 elements yield flux continuity across cell interfaces, while P0 test functions enforce cellwise mass conservation. 
With RT0 basis functions $\{\boldsymbol{\psi}_i\}_{i=1}^{N_f}$ associated with mesh faces and
piecewise constant basis functions 
 (cell indicators) $\{\chi_j\}_{j=1}^{N_c}$, we express the discrete unknowns as
\[
\mathbf{u}_{t,h}=\sum_{i=1}^{N_f}U_i\boldsymbol{\psi}_i,\qquad
p_h=\sum_{j=1}^{N_c}P_j\chi_j,\qquad
S_{g,h}=\sum_{j=1}^{N_c}S_j\chi_j .
\]
The coefficient vectors
\[
\mathbf{U}=[U_1,\ldots,U_{N_f}]^T\in\mathbb{R}^{N_f},\qquad
\mathbf{P}=[P_1,\ldots,P_{N_c}]^T\in\mathbb{R}^{N_c},\qquad
\mathbf{S}=[S_1,\ldots,S_{N_c}]^T\in\mathbb{R}^{N_c},
\]
collect the degrees of freedom associated with the total velocity, pressure, and gas saturation, respectively. 
Let $\mathbf{D}=\mathrm{diag}(|K_1|,\ldots,|K_{N_c}|)$ and define the discrete divergence operator and RT0 mass matrix by
\[
(\mathbf{C})_{ij}=\int_{K_j}\nabla\!\cdot\boldsymbol{\psi}_i\,\mathrm{d}x,\qquad
(\mathbf{M}_\lambda(\mathbf{S}))_{ij}=\int_{\Omega}(\lambda_t(S_{g,h})\mathbf{K})^{-1}\boldsymbol{\psi}_j\!\cdot\!\boldsymbol{\psi}_i\,\mathrm{d}x .
\]
The vectors $\mathbf{M}_t(\mathbf{P},\mathbf{S})$ and $\mathbf{M}_{\mathrm{CO}_2}(\mathbf{P},\mathbf{S})$ collect the cell values of the accumulations in \eqref{eq:accumulation-definitions}, and $\mathbf{Q}_t,\mathbf{Q}_{\mathrm{CO}_2}$ collect the corresponding cell-average sources.

We discretize the time interval $[0,T]$ into subintervals with endpoints $0=t^0<t^1<\cdots<t^{N_T}=T$ and define time step sizes $\Delta t^n=t^{n+1}-t^n$. Using backward Euler integration in time, the fully discrete high-fidelity system at time $t^{n+1}$ reads
\begin{equation}
\left\{
\begin{aligned}
\mathbf{R}_{\mathrm{mom}}^{n+1}
&:=\mathbf{M}_\lambda(\mathbf{S}^{n+1})\mathbf{U}^{n+1}
-\mathbf{C}\mathbf{P}^{n+1}
+\mathbf{b}_c(\mathbf{S}^{n+1})
+\mathbf{b}_D
=\mathbf{0},\\
\mathbf{R}_{t}^{n+1}
&:=\frac{1}{\Delta t^n}\mathbf{D}
\Big(\mathbf{M}_t(\mathbf{P}^{n+1},\mathbf{S}^{n+1})
-\mathbf{M}_t(\mathbf{P}^{n},\mathbf{S}^{n})\Big)
+\boldsymbol{\mathcal{F}}_t(\mathbf{U}^{n+1},\mathbf{P}^{n+1},\mathbf{S}^{n+1})
-\mathbf{D}\mathbf{Q}_t^{n+1}
=\mathbf{0},\\
\mathbf{R}_{\mathrm{CO}_2}^{n+1}
&:=\frac{1}{\Delta t^n}\mathbf{D}
\Big(\mathbf{M}_{\mathrm{CO}_2}(\mathbf{P}^{n+1},\mathbf{S}^{n+1})
-\mathbf{M}_{\mathrm{CO}_2}(\mathbf{P}^{n},\mathbf{S}^{n})\Big)
+\boldsymbol{\mathcal{F}}_{\mathrm{CO}_2}(\mathbf{U}^{n+1},\mathbf{P}^{n+1},\mathbf{S}^{n+1})
-\mathbf{D}\mathbf{Q}_{\mathrm{CO}_2}^{n+1}
=\mathbf{0},
\end{aligned}
\right.
\label{eq:global residual}
\end{equation}
where boundary term $(\mathbf{b}_D)_i=\int_{\Gamma_D}p_D(\boldsymbol{\psi}_i\!\cdot\!\mathbf{n})\,\mathrm{d}s$.
To explicitly define the discrete flux terms appearing in~\eqref{eq:global residual}, we first introduce face-based notation. 
For an interior face $e=K^-\cap K^+$, we fix an orientation by defining a unit normal vector $\mathbf{n}_e$ from cell $K^-$ toward $K^+$. The RT0 normal total flux across face $e$ is given by 
\[
F_e=\int_e\mathbf{u}_{t,h}\!\cdot\!\mathbf{n}_e\,\mathrm{d}s .
\]
For a generic cellwise quantity $a$, we define the upwind and centered face values by
\[
a^\uparrow_e =
\begin{cases}
a_{K^-},&F_e\ge 0,\\
a_{K^+},&F_e<0,
\end{cases}
\qquad
\overline{a}_e=\frac{a_{K^-}+a_{K^+}}{2}.
\]
With these notations, we define 
the capillary vector $\mathbf{b}_c$ in the RT0 load associated with $\mathbf{f}_c=f_a\nabla p_c$ is:
\[
(\mathbf{b}_c(\mathbf{S}))_i
\approx
\sum_{e\in\mathcal{E}_h}|e|\,
\overline{\left(\frac{\lambda_a}{\lambda_t}\right)}_e
(\partial_n p_c)_e
(\boldsymbol{\psi}_i\!\cdot\!\mathbf{n}_e)_e ,
\]
where $(\boldsymbol{\psi}_i\!\cdot\!\mathbf{n}_e)_e=|e|^{-1}\int_e\boldsymbol{\psi}_i\!\cdot\!\mathbf{n}_e\,\mathrm{d}s$, 
the normal derivative
\[
(\partial_n p_c)_e=\frac{p_c(S_{K^+})-p_c(S_{K^-})}{d_e},
\]
where $d_e$ is the projected distance between the two cell centers along $\mathbf{n}_e$. On boundary faces, $\mathbf{n}_e$ is the outward boundary normal and the boundary saturation is used in place of the missing neighboring cell. This face-based treatment is consistent with mixed/DG capillary discretizations for two-phase flow, where advective coefficients are upwinded and capillary forces are evaluated from pressure jumps across faces \cite{hoteitFiroozabadi2006,
hoteit2008numerical,
moncorge2020upwinding,arbogastTao2019}.
Next, the discrete flux vectors $\boldsymbol{\mathcal{F}}_t$ and $\boldsymbol{\mathcal{F}}_{\mathrm{CO}_2}$ assemble the oriented face fluxes in the cell balances:
\[
(\boldsymbol{\mathcal{F}}_t)_K=\sum_{e\subset\partial K}|e|(\mathbf{m}_t\!\cdot\!\mathbf{n}_{K,e})_e,\qquad
(\boldsymbol{\mathcal{F}}_{\mathrm{CO}_2})_K=\sum_{e\subset\partial K}|e|(\mathbf{m}_{\mathrm{CO}_2}\!\cdot\!\mathbf{n}_{K,e})_e .
\]
With $\overline{K}_{n,e}$ denoting the harmonically averaged normal permeability on $e$, the face fluxes used in the two balance equations are
\begin{align}
\label{dfn:m_t}
(\mathbf{m}_t\!\cdot\!\mathbf{n}_e)_e
&\approx
\Big(\rho_a f_a+\rho_g f_g\Big)^{\uparrow}_e\frac{F_e}{|e|}
+ \overline{K}_{n,e}\,
\overline{\left(\frac{\lambda_a\lambda_g}{\lambda_t}\right)}_e
\overline{(\rho_g-\rho_a)}_e\,
(\partial_n p_c)_e,\\
\label{dfn:m_co2}
(\mathbf{m}_{\mathrm{CO}_2}\!\cdot\!\mathbf{n}_e)_e
&\approx
\Big(\rho_a f_a X_a^{\mathrm{CO}_2}+\rho_g f_g X_g^{\mathrm{CO}_2}\Big)^{\uparrow}_e\frac{F_e}{|e|}
+ \overline{K}_{n,e}\,
\overline{\left(\frac{\lambda_a\lambda_g}{\lambda_t}\right)}_e
\overline{(\rho_g X_g^{\mathrm{CO}_2}-\rho_a X_a^{\mathrm{CO}_2})}_e\,
(\partial_n p_c)_e .
\end{align}
For boundary faces, the same notation is used with the prescribed boundary state or normal flux replacing the missing neighboring cell.

\section{Sequential Fully Implicit Domain-Decomposition Formulation}
\label{sec:sfi-dd}

\subsection{Sequential Fully Implicit Splitting}
\label{subsec:sfi-splitting}

The fully discrete high-fidelity residual system \eqref{eq:global residual} can, in principle, be solved monolithically by a fully implicit Newton method (FI), where all unknowns are updated simultaneously in a single coupled nonlinear solve. 
Such a formulation is robust and directly targets the fully coupled residual, but it also leads to large and strongly coupled Jacobian systems. For large-scale compositional and multiphase porous-media simulations, the solution of these monolithic systems may become computationally demanding, especially in the presence of strong nonlinearities, phase behavior effects, and heterogeneous coefficients \cite{chen2006computational,
moncorge2018sfi,jiangTchelepi2019,jiangTominZhou2021}. 
A widely used alternative in reservoir simulation is the Sequential Fully Implicit (SFI) method \cite{chen2006computational,moncorge2018sfi,moncorge2020upwinding,jiangTchelepi2019}.  In SFI, the coupled flow and transport system is split into two implicit subproblems — a pressure/flow solve and a transport solve — solved sequentially within an outer nonlinear loop. Each subproblem remains fully implicit in its own variables, while the global coupling is treated through outer iterations that update the frozen variables. 
Separating these two physics-based subproblems can leverage the distinct mathematical properties of each subsystem,  thus 
enhance computational efficiency compared to monolithic FI approaches. Additionally, this separation is natural for reduced-order modeling, allowing distinct reduced-order representations and tailored approximations for each subsystem as presented in Section \ref{sec:ddrom-compact}. 
We emphasize, however, that the efficiency of SFI depends on the convergence of the outer coupling iteration; for strongly coupled problems, basic SFI may require many outer iterations, and acceleration or inexact inner-solve strategies can be beneficial \cite{moncorge2018sfi,jiangTchelepi2019,jiangTominZhou2021,liTominTchelepi2021}.

We now formulate briefly the SFI method for time step $t^{n+1}$. Let $k$ denote the SFI outer iteration counter within this time step. We collect the total velocity and pressure unknowns as
\[
\mathbf{X}:=(\mathbf{U},\mathbf{P}), \qquad
\mathbf{R}_{p}(\mathbf{X};\mathbf{S})
:=
\begin{bmatrix}
\mathbf{R}_{\mathrm{mom}}(\mathbf{U},\mathbf{P},\mathbf{S})\\
\mathbf{R}_{t}(\mathbf{U},\mathbf{P},\mathbf{S})
\end{bmatrix},
\]
where the semicolon separates the active unknowns from the variables frozen in the current subproblem solve. Given an iterate $\mathbf{S}^{(n+1,k)}$ for gas saturation, the \textit{pressure step} consists of solving
\begin{equation}
    \mathbf{R}_{p}(\mathbf{X}^{(n+1,k+1)};\mathbf{S}^{(n+1,k)})=\mathbf{0}.
    \label{eq:sfi-pressure-compact}
\end{equation}
This nonlinear system is solved by Newton's method with respect to $\mathbf{X}$ only. Equivalently, the pressure-step Jacobian contains derivatives with respect to $\mathbf{U}$ and $\mathbf{P}$, while derivatives with respect to $\mathbf{S}$ are omitted because saturation is fixed during this subproblem.

After the pressure step, the updated pressure and total velocity are kept fixed, and the transport subproblem is solved for the new saturation. Since the transport equation is the \(\mathrm{CO}_2\)-component balance with \(\mathbf X\) frozen, we introduce the notation
\[
\mathbf R_s(\mathbf S;\mathbf X)
:=
\mathbf R_{\mathrm{CO}_2}(\mathbf U,\mathbf P,\mathbf S).
\]
Here, \(\mathbf R_{\mathrm{CO}_2}\) denotes the physical \(\mathrm{CO}_2\)-component residual in the fully coupled system~\eqref{eq:global residual}, whereas \(\mathbf R_s\) denotes the same residual viewed as a nonlinear equation for \(\mathbf S\) with \(\mathbf X\) fixed. The \textit{transport step} is therefore written as
\begin{equation}
    \mathbf R_s\left(
    \mathbf S^{(n+1,k+1)};
    \mathbf X^{(n+1,k+1)}
    \right)
    =
    \mathbf 0 .
    \label{eq:sfi-transport-compact}
\end{equation}
The corresponding Newton linearization is taken with respect to \(\mathbf S\) only; the total velocity and pressure fields are treated as fixed coefficients in the accumulation and flux terms.

The outer SFI iteration proceeds
between \eqref{eq:sfi-pressure-compact} and \eqref{eq:sfi-transport-compact} 
 until a convergence criterion is satisfied,
for example a tolerance on saturation increment or the coupled residuals.  
Once converged, the accepted solution is denoted by
\[
    \mathbf{X}^{n+1}=\mathbf{X}^{(n+1,k_\star)},
    \qquad
    \mathbf{S}^{n+1}=\mathbf{S}^{(n+1,k_\star)},
\]
where $k_\star$ is the final outer iteration. If the outer iteration is converged to the prescribed coupled tolerance, the SFI solution is consistent with the fully implicit residual \eqref{eq:global residual}, while retaining a block structure that is useful for both solver design and reduced-order approximation.

\subsection{Domain-Decomposition Schur Formulation}
\label{subsec:dd-schur}

We now recast the two SFI subproblems in a nonoverlapping domain-decomposition (DD) form \cite{alboin2000domain,ganis2014global,skogestad2013domain,cheng2024adaptive}. Let
$\Omega=\cup_{i=1}^{N_s}\Omega_i$ be a conforming nonoverlapping partition, and let $\Gamma$ be the union of all internal interfaces. For each subdomain, $\Gamma_{D,i}$ and $\Gamma_{N,i}$ denote the portions of the physical Dirichlet and Neumann boundaries on $\partial\Omega_i$. The local RT0--P0 unknowns are
$\mathbf{U}_i\in\mathbb{R}^{N_{f,i}}$, $\mathbf{P}_i\in\mathbb{R}^{N_{c,i}}$, and $\mathbf{S}_i\in\mathbb{R}^{N_{c,i}}$.
The global interface mesh contains $N_\Gamma$ degrees of freedom. For each subdomain, 
let $\mathcal{E}_{\Gamma,i}$ be the set of interface faces on $\partial\Omega_i\cap\Gamma$ and let $\mathbf{n}_{i,e}$ be the unit normal on $e\in\mathcal{E}_{\Gamma,i}$ pointing outward from $\Omega_i$. For a local interface face $e$, the adjacent cell in $\Omega_i$ is denoted by $K_i(e)$.
We define a global interface pressure
\[
\boldsymbol{\lambda}_p\in\mathbb{R}^{N_\Gamma},
\]
whose restriction to subdomain interface
$\partial\Omega_i\cap\Gamma$
 is given by
\[
\boldsymbol{\lambda}_{p,i}=\mathbf{H}_i\boldsymbol{\lambda}_p,
\qquad
\mathbf{H}_i\in\{0,1\}^{N_{\Gamma,i}\times N_\Gamma},
\]
where $\mathbf{H}_i$ is Boolean restriction matrix which maps the global interface degree of freedoms to the local ordering. This global interface pressure serves as the Dirichlet boundary condition for the local pressure problems, as detailed in Section \ref{subsubsec:dd-pressure}. Similarly, we introduce a global interface saturation 
\[
\boldsymbol{\lambda}_s\in\mathbb{R}^{N_\Gamma}, 
\]
with local restriction to each subdoamin is similarly defined using $\mathbf{H}_i$. 
This global saturation also serves as boundary data for the local transport problems, as described in Section \ref{subsubsec:dd-transport}.

\subsubsection{Pressure Step}
\label{subsubsec:dd-pressure}

Given a fixed saturation iterate $\mathbf{S}^{(n+1,k)}$, the pressure step solves for the local states
\[
\mathbf{X}_i:=\begin{bmatrix}\mathbf{U}_i\\ \mathbf{P}_i\end{bmatrix},
\qquad i=1,\ldots,N_s,
\]
and for the global interface pressure trace $\boldsymbol{\lambda}_p\in\mathbb{R}^{N_\Gamma}$. 
If $\boldsymbol{\lambda}_p$ is prescribed, the subdomains are independent: on $\Omega_i$ the local pressure problem is
\[
\mathbf{R}_{p,i}(\mathbf{X}_i,\boldsymbol{\lambda}_{p,i};\mathbf{S}_i)
:=
\begin{bmatrix}
\mathbf{R}_{\mathrm{mom},i}(\mathbf{U}_i,\mathbf{P}_i,\boldsymbol{\lambda}_{p,i};\mathbf{S}_i)\\
\mathbf{R}_{t,i}(\mathbf{U}_i,\mathbf{P}_i;\mathbf{S}_i)
\end{bmatrix} = 0,
\]
where the interface pressure acts as a weak Dirichlet contribution in the local momentum equation. Specifically, using the notation from Section \ref{sec:disc}, the local momentum residual can be written as
\[
\mathbf{R}_{\mathrm{mom},i}
=
\mathbf{M}_{\lambda,i}(\mathbf{S}_i)\mathbf{U}_i
-\mathbf{C}_i\mathbf{P}_i
+\mathbf{b}_{c,i}(\mathbf{S}_i)
+\mathbf{b}_{D,i}
+\mathbf{M}_{p,i}\boldsymbol{\lambda}_{p,i},
\]
where $\mathbf{b}_{D,i}$ contains the physical Dirichlet-boundary contribution and $\mathbf{M}_{p,i}\boldsymbol{\lambda}_{p,i}$ is the analogous contribution induced by
interface pressure. The matrix $\mathbf{M}_{p,i}$ is constant and is determined by the normal traces of the RT0 basis functions on the interface edges. 

However, the global interface pressure is not prescribed data; it is determined by enforcing total mass flux conservation across all interfaces. To this end, we define local interface total mass fluxes by
\[
\mathbf{F}_{t,i}(\mathbf{X}_i,\boldsymbol{\lambda}_{p,i};\mathbf{S}_i)
:=
\big(F_{t,i,e}(\mathbf{X}_i,\boldsymbol{\lambda}_{p,i};\mathbf{S}_i)\big)_{e\in\mathcal{E}_{\Gamma,i}},
\]
with each face flux given by
\begin{equation*}
F_{t,i,e}(\mathbf{X}_i,\boldsymbol{\lambda}_{p,i};\mathbf{S}_i)
:=
|e|\,\big(\mathbf{m}_{t}(\mathbf{u}_{t,h,i},p_{h,i},S_{g,h,i})
\cdot\mathbf{n}_{i,e}\big)_e.
\label{eq:local-total-interface-flux}
\end{equation*}
The numerical flux $(\mathbf{m}_{t}\cdot\mathbf{n}_{i,e})_e$ is evaluated by the formula \eqref{dfn:m_t}. 
The flux continuity condition across the interfaces is enforced by the global interface residual
\begin{equation*}
    \mathbf{R}_{p,\Gamma}(\{\mathbf{X}_i\},\boldsymbol{\lambda}_p)
    :=
    \sum_{i=1}^{N_s}
    \mathbf{H}_i^{\top}\mathbf{F}_{t,i}(\mathbf{X}_i,\mathbf{H}_i\boldsymbol{\lambda}_p;\mathbf{S}_i) = 0.
    \label{eq:pressure-interface-compact}
\end{equation*}

The high-fidelity DD pressure problem is therefore
\begin{equation*}
\left\{
\begin{aligned}
\mathbf{R}_{p,i}(\mathbf{X}_i,\mathbf{H}_i\boldsymbol{\lambda}_p;\mathbf{S}_i)&=\mathbf{0},
&& i=1,\ldots,N_s,\\
\mathbf{R}_{p,\Gamma}(\{\mathbf{X}_i\},\boldsymbol{\lambda}_p)&=\mathbf{0}.
\end{aligned}
\right.
\label{eq:DD pressure FOM}
\end{equation*}
We solve this coupled system using a Newton--Schur method. 
For notational brevity, we omit the global time step index $n$ and the Picard outer iteration index $k$. 
At Newton step $m$, the local residual linearization is
\begin{equation*}
    \mathbf{K}_{p,i}^{(m)}\delta\mathbf{X}_i
    +\widetilde{\mathbf{M}}_{p,i}\delta\boldsymbol{\lambda}_{p,i}
    =
    -\mathbf{R}_{p,i}^{(m)},
    \qquad i=1,\ldots,N_s,
    \label{eq:local-pressure-linear-compact}
\end{equation*}
where
\[
\mathbf{K}_{p,i}^{(m)}
:=
\frac{\partial\mathbf{R}_{p,i}}{\partial\mathbf{X}_i}
\Big|_{(\mathbf{X}_i^{(m)},\boldsymbol{\lambda}_{p,i}^{(m)})},
\qquad
\widetilde{\mathbf{M}}_{p,i}
:=
\frac{\partial\mathbf{R}_{p,i}}{\partial\boldsymbol{\lambda}_{p,i}}
=
\begin{bmatrix}
\mathbf{M}_{p,i}\\
\mathbf{0}
\end{bmatrix}.
\]
The matrix $\widetilde{\mathbf{M}}_{p,i}$ is independent of the Newton iterate because the interface pressure enters the local momentum residual linearly. Eliminating the interior increment gives
\begin{equation}
    \delta\mathbf{X}_i
    =
    -(\mathbf{K}_{p,i}^{(m)})^{-1}\mathbf{R}_{p,i}^{(m)}
    -(\mathbf{K}_{p,i}^{(m)})^{-1}\widetilde{\mathbf{M}}_{p,i}\mathbf{H}_i\delta\boldsymbol{\lambda}_p .
    \label{eq:pressure-interior-elimination}
\end{equation}
The linearization of the interface residual is
\begin{equation}
    \sum_{i=1}^{N_s}
    \mathbf{H}_i^{\top}
    \left(
    \mathbf{A}_{t,i}^{(m)}\delta\mathbf{X}_i
    +
    \mathbf{B}_{t,i}^{(m)}\mathbf{H}_i\delta\boldsymbol{\lambda}_p
    \right)
    =
    -\mathbf{R}_{p,\Gamma}^{(m)},
    \label{eq:pressure-interface-linear}
\end{equation}
where
\[
\mathbf{A}_{t,i}^{(m)}
:=
\frac{\partial\mathbf{F}_{t,i}}{\partial\mathbf{X}_i},
\qquad
\mathbf{B}_{t,i}^{(m)}
:=
\frac{\partial\mathbf{F}_{t,i}}{\partial\boldsymbol{\lambda}_{p,i}},
\]
with $\mathbf{S}_i$ fixed in both derivatives.
Substituting \eqref{eq:pressure-interior-elimination} into \eqref{eq:pressure-interface-linear} yields the global pressure Schur complement system
\begin{equation}
    \mathbf{S}_{p,\Gamma}^{(m)}\delta\boldsymbol{\lambda}_p
    =
    \mathbf{g}_{p,\Gamma}^{(m)},
    \label{eq:pressure-schur-compact}
\end{equation}
where
\begin{align*}
\mathbf{S}_{p,\Gamma}^{(m)}
&=
\sum_{i=1}^{N_s}
\mathbf{H}_i^{\top}
\left[
\mathbf{B}_{t,i}^{(m)}\mathbf{H}_i
-
\mathbf{A}_{t,i}^{(m)}
(\mathbf{K}_{p,i}^{(m)})^{-1}
\widetilde{\mathbf{M}}_{p,i}\mathbf{H}_i
\right],
\notag\\
\mathbf{g}_{p,\Gamma}^{(m)}
&=
-\mathbf{R}_{p,\Gamma}^{(m)}
+
\sum_{i=1}^{N_s}
\mathbf{H}_i^{\top}\mathbf{A}_{t,i}^{(m)}
(\mathbf{K}_{p,i}^{(m)})^{-1}
\mathbf{R}_{p,i}^{(m)}.
\label{eq:pressure-schur-blocks}
\end{align*}
After solving \eqref{eq:pressure-schur-compact} for $\delta\boldsymbol{\lambda}_p$, the local increments $\delta\mathbf{X}_i$ are recovered from \eqref{eq:pressure-interior-elimination}. In Section \ref{subsec:pressure-rom}, we derived a similar Schur complement system tailored for the reduced-order problem.

\subsubsection{Transport Step}
\label{subsubsec:dd-transport}
We now derive the DD formulation of the transport step. 
Although the resulting algebraic structure is again a Schur complement, the role of the interface trace is different from that in the pressure step. In the pressure problem, the interface pressure $\boldsymbol\lambda_p$ is required for every local pressure solve. In contrast, in the transport problem, after the pressure step has provided the frozen velocity and pressure fields, the numerical $\mathrm{CO}_2$ flux is evaluated with an upwind saturation state. Hence an interface saturation trace is required only on inflow faces. On outflow faces, the flux is fully determined by the interior saturation of the upstream subdomain.

To determine the interface saturation, similar to the pressure step, we impose CO$_2$ mass flux continuity across the interface.  
Define the local interface CO$_2$ flux vector by
\[
\mathbf{F}_{c,i}(\mathbf{S}_i,\boldsymbol{\lambda}_{s,i};\mathbf{X}_i)
:=
\big(F_{c,i,e}(\mathbf{S}_i,\boldsymbol{\lambda}_{s,i};\mathbf{X}_i)\big)_{e\in\mathcal{E}_{\Gamma,i}},
\qquad
F_{c,i,e}
:=
|e|\,\big(\mathbf{m}_{\mathrm{CO}_2}\cdot\mathbf{n}_{i,e}\big)_e,
\]
where the numerical flux $(\mathbf{m}_{\mathrm{CO}_2}\cdot\mathbf{n}_{i,e})_e$ follows the discrete flux formula given in \eqref{dfn:m_co2}. 
The global interface CO$_2$  flux continuity can then be expressed as
\begin{equation}
\mathbf{R}_{s,\Gamma}(\{\mathbf{S}_i\},\boldsymbol{\lambda}_s)
:=
\sum_{i=1}^{N_s}
\mathbf{H}_i^\top
\mathbf{F}_{c,i}(\mathbf{S}_i,\mathbf{H}_i\boldsymbol{\lambda}_s;\mathbf{X}_i) = 0.
\label{eq:transport-dd-compact}
\end{equation}
Combining this interface flux continuity condition with the local CO$_2$ residual equations, we obtain the coupled DD transport system
\begin{equation*}
\left\{
\begin{aligned}
\mathbf{R}_{s,i}(\mathbf{S}_i,\mathbf{H}_i\boldsymbol{\lambda}_s;\mathbf{X}_i)&=\mathbf{0},
&& i=1,\ldots,N_s,\\
\mathbf{R}_{s,\Gamma}(\{\mathbf{S}_i\},\boldsymbol{\lambda}_s)&=\mathbf{0}.
\end{aligned}
\right.
\label{eq:DD transport FOM}
\end{equation*}

This coupled nonlinear system is solved by a Newton--Schur method analogous to that in the pressure step. At Newton iteration $m$, the local transport linearization is
\begin{equation*}
    \mathbf{J}_{s,i}^{(m)}\delta\mathbf{S}_i
    +
    \mathbf{G}_{s,i}^{(m)}\mathbf{H}_i\delta\boldsymbol{\lambda}_s
    =
    -\mathbf{R}_{s,i}^{(m)},
    \qquad i=1,\ldots,N_s,
    \label{eq:transport-local-linear}
\end{equation*}
where
\[
\mathbf{J}_{s,i}^{(m)}
:=
\frac{\partial\mathbf{R}_{s,i}}{\partial\mathbf{S}_i},
\qquad
\mathbf{G}_{s,i}^{(m)}
:=
\frac{\partial\mathbf{R}_{s,i}}{\partial\boldsymbol{\lambda}_{s,i}}.
\]
Thus
\begin{equation}
\delta\mathbf{S}_i
=
-(\mathbf{J}_{s,i}^{(m)})^{-1}\mathbf{R}_{s,i}^{(m)}
-
(\mathbf{J}_{s,i}^{(m)})^{-1}\mathbf{G}_{s,i}^{(m)}
\mathbf{H}_i\delta\boldsymbol{\lambda}_s .
\label{eq:transport-interior-elimination}
\end{equation}
Substituting \eqref{eq:transport-interior-elimination} into the linearized interface flux equation,
\begin{equation*}
\sum_{i=1}^{N_s}
\mathbf{H}_i^\top
\left(
\mathbf{A}_{c,i}^{(m)}\delta\mathbf{S}_i
+
\mathbf{B}_{c,i}^{(m)}\mathbf{H}_i\delta\boldsymbol{\lambda}_s
\right)
=
-\mathbf{R}_{s,\Gamma}^{(m)},
\label{eq:transport-interface-linear}
\end{equation*}
with
\[
\mathbf{A}_{c,i}^{(m)}
:=
\frac{\partial\mathbf{F}_{c,i}}{\partial\mathbf{S}_i},
\qquad
\mathbf{B}_{c,i}^{(m)}
:=
\frac{\partial\mathbf{F}_{c,i}}{\partial\boldsymbol{\lambda}_{s,i}},
\]
we obtain the transport Schur complement system
\begin{equation}
\mathbf{S}_{s,\Gamma}^{(m)}
\delta\boldsymbol{\lambda}_s
=
\mathbf{g}_{s,\Gamma}^{(m)},
\label{eq:transport-schur-compact}
\end{equation}
where 
\begin{align*}
\mathbf{S}_{s,\Gamma}^{(m)}
&=
\sum_{i=1}^{N_s}
\mathbf{H}_i^\top
\left[
\mathbf{B}_{c,i}^{(m)}\mathbf{H}_i
-
\mathbf{A}_{c,i}^{(m)}
(\mathbf{J}_{s,i}^{(m)})^{-1}
\mathbf{G}_{s,i}^{(m)}\mathbf{H}_i
\right],
\notag\\
\mathbf{g}_{s,\Gamma}^{(m)}
&=
-\mathbf{R}_{s,\Gamma}^{(m)}
+
\sum_{i=1}^{N_s}
\mathbf{H}_i^\top\mathbf{A}_{c,i}^{(m)}
(\mathbf{J}_{s,i}^{(m)})^{-1}
\mathbf{R}_{s,i}^{(m)}.
\label{eq:transport-schur-blocks}
\end{align*}

Equations \eqref{eq:pressure-schur-compact} and \eqref{eq:transport-schur-compact} are the high-fidelity interface systems that the DD-ROM formulation introduced in Section \ref{sec:ddrom-compact} preserves at the reduced-order level.

\section{DD-ROM Formulation}
\label{sec:ddrom-compact}

In this section, we derive the domain-decomposition reduced-order model (DD-ROM) based on the high-fidelity DD-SFI formulation introduced in Section~\ref{sec:sfi-dd}. The main idea is to replace the local subdomain unknowns by reduced coordinates while preserving the same interface coupling structure as in the high-fidelity Schur systems. Thus the pressure and transport subproblems are reduced locally, but the interface residuals are either retained in the full interface space or projected onto a low-dimensional port space.

The construction has three components. First, local reduced spaces are built from high-fidelity DD-SFI snapshots, with an interface lifting used for the pressure variable. Second, the pressure step is reduced by either Galerkin projection or least-squares Petrov--Galerkin (LSPG) projection, leading to full-interface and port-reduced Schur systems. Third, the transport step is reduced in the saturation variable while retaining the conservative interface-flux coupling.

\subsection{Local Reduced Spaces and Interface Lifting}
\label{subsec:local-spaces-lifting}
For each subdomain $\Omega_i$, the reduced spaces are constructed from high-fidelity DD-SFI snapshots. We denote by
\[
\{\mathbf{U}_i^k,\mathbf{P}_i^k,\mathbf{S}_i^k,\boldsymbol{\lambda}_{p,i}^k\}_{k=1}^{n_{\mathrm{train}}}
\]
the local total-velocity, pressure, saturation, and interface-pressure. 

In the pressure step, the interface pressure acts as nonhomogeneous Dirichlet data in the mixed pressure solve. Therefore before applying POD, we decompose the velocity and pressure snapshots into an interface lifting and a homogeneous fluctuation. This local-interior/interface split is related to localized reduced-basis, local-ROM, and static-condensation model reduction constructions \cite{kaulmann2015lrbms,ohlbergerSchindler2015,buhrSmetana2018,smetana2015certification,smetanaPatera2016}. 

At each training state we compute an auxiliary RT0--P0 Darcy extension on $\Omega_i$: 
\begin{equation*} 
\left\{ \begin{aligned} \mathbf u_{\mathrm{int},i}+\nabla p_{\mathrm{int},i} &=\mathbf 0 &&\text{in }\Omega_i,\\ \nabla\cdot\mathbf u_{\mathrm{int},i} &=0 &&\text{in }\Omega_i,\\ p_{\mathrm{int},i} &=p_D &&\text{on }\Gamma_{D,i},\\ p_{\mathrm{int},i} &=\boldsymbol\lambda_{p,i} &&\text{on }\partial\Omega_i\cap\Gamma,\\ \mathbf u_{\mathrm{int},i}\cdot\mathbf n &=g_N &&\text{on }\Gamma_{N,i}. \end{aligned} \right. 
\label{eq:darcy-interface-lifting-local} \end{equation*}
By linearity, its discrete solution
\[ \mathbf X_{\mathrm{int},i} := \begin{bmatrix} \mathbf U_{\mathrm{int},i}\\ \mathbf P_{\mathrm{int},i} \end{bmatrix} \]
can be written in affine form: 
\begin{equation*} 
\mathbf X_{\mathrm{int},i}(\boldsymbol\lambda_{p,i}) = \mathbf L_i\boldsymbol\lambda_{p,i}+\mathbf b_i. 
\end{equation*}
The matrix $\mathbf{L}_i$ is assembled once offline from local  solves associated with the interface pressure degrees of freedom, while $\mathbf{b}_i$ contains the response to the physical boundary data.

POD is then applied to the fluctuation snapshots
\[ \mathbf U_i^k-\mathbf U_{\mathrm{int},i}^k, \qquad \mathbf P_i^k-\mathbf P_{\mathrm{int},i}^k, \qquad \mathbf S_i^k . \]
This gives local bases
\[
\mathbf{Z}_{u,i}\in\mathbb{R}^{N_{f,i}\times n_{u,i}},
\qquad
\mathbf{Z}_{p,i}\in\mathbb{R}^{N_{c,i}\times n_{p,i}},
\qquad
\mathbf{Z}_{s,i}\in\mathbb{R}^{N_{c,i}\times n_{s,i}}.
\]
The dimensions are selected either by POD energy tolerances or by the HAPOD construction described in Section~\ref{sec:well-aware-compact}. Defining
\[ \mathbf a_i := \begin{bmatrix} \boldsymbol\alpha_{u,i}\\ \boldsymbol\alpha_{p,i} \end{bmatrix}, \qquad \boldsymbol\Phi_i := \operatorname{blkdiag}(\mathbf Z_{u,i},\mathbf Z_{p,i}), \] 
the reduced ansatz used throughout the pressure and transport ROMs is
\begin{equation}
\left\{
\begin{aligned}
\mathbf X_i &\approx \mathbf X_i^r := \boldsymbol\Phi_i\mathbf a_i + \mathbf L_i\mathbf H_i\boldsymbol\lambda_p + \mathbf b_i ,\\
\mathbf S_i 
&
\approx \mathbf S_i^r := \mathbf Z_{s,i}\boldsymbol\alpha_{s,i}. 
\end{aligned}
\right.
\label{eq:reduced-ansatz-local}
\end{equation}
We note that no analogous elliptic lifting is introduced for saturation, since the interface saturation is used only as inflow data in the transport residual, not as a two-sided Dirichlet trace.

\subsection{Pressure-Step Reduced Schur Systems}
\label{subsec:pressure-rom}

In the ROM pressure step, we retain the nonoverlapping interface formulation of Section~\ref{subsubsec:dd-pressure}, but replace the local pressure-step state \(\mathbf X_i\) by the reduced representation \eqref{eq:reduced-ansatz-local}. The interface pressure can either be kept in the full trace space, leading to full-interface coupling, or approximated in a port space, leading to port-reduced (PR) coupling. We first derive the full-interface formulation.

Inserting the reduced ansatz \eqref{eq:reduced-ansatz-local} into the local high-fidelity residual yields the reduced residual system
\begin{equation}
\left\{
\begin{aligned}
\mathbf{R}_{p,i}
\big(\mathbf{X}_{i}^{r},\boldsymbol{\lambda}_{p,i};\mathbf{S}_{i}^{r}\big)
&=\mathbf{0},
&& i=1,\ldots,N_s,\\
\mathbf{R}_{p,\Gamma}(\{\mathbf{a}_i\},\boldsymbol{\lambda}_p)&=\mathbf{0}.
\end{aligned}
\right.
\label{eq:full-interface-rom-compact}
\end{equation}
The first equation is overdetermined, because \(\mathbf R_{p,i}\) has the dimension of the full local RT0--P0 pressure residual, whereas \(\mathbf a_i\) contains only reduced coordinates.
The interface residual in the second equation is not projected in the full-interface formulation, so total flux continuity is still enforced on each interface face.

We consider both Galerkin and LSPG reductions of the local residual in \eqref{eq:full-interface-rom-compact}. The Galerkin variant enforces the projected condition
\begin{equation*}
\widehat{\mathbf{R}}_{p,i}^{G}
:=
\mathbf{\Phi}_i^T
\mathbf{R}_{p,i}
\big(\mathbf{X}_{i}^{r},\boldsymbol{\lambda}_{p,i};\mathbf{S}_{i}^{r}\big)
=\mathbf{0}.
\label{eq:pressure-galerkin-residual} 
\end{equation*}
The LSPG variant determines the reduced coordinates by minimizing the same local residual: \begin{equation} 
\mathbf a_i = \arg\min_{\mathbf y\in\mathbb R^{n_{u,i}+n_{p,i}}} \frac{1}{2} \left\| \mathbf R_{p,i} \left( \boldsymbol\Phi_i\mathbf y + \mathbf L_i\mathbf H_i\boldsymbol\lambda_p + \mathbf b_i, \boldsymbol\lambda_{p,i}; \mathbf S_i^r \right) \right\|_2^2 . \label{eq:pressure-lspg-nonlinear-ls} 
\end{equation}
The first-order optimality condition of \eqref{eq:pressure-lspg-nonlinear-ls} gives the nonlinear normal equation 
\begin{equation*} 
\widehat{\mathbf R}_{p,i}^{L} := \left( \frac{\partial \mathbf R_{p,i}}{\partial \mathbf a_i} \right)^\top \mathbf R_{p,i} \left( \mathbf X_i^r,\boldsymbol\lambda_{p,i};\mathbf S_i^r \right) = \mathbf 0 . 
\label{eq:pressure-lspg-normal-residual} 
\end{equation*}
Thus both local ROMs can be written as: 
\begin{equation*} 
\left\{ \begin{aligned} \widehat{\mathbf R}_{p,i}^{\star} (\mathbf a_i,\boldsymbol\lambda_p;\boldsymbol\alpha_{s,i}) &= \mathbf 0, && i=1,\ldots,N_s, \\ \mathbf R_{p,\Gamma} (\{\mathbf X_i^r\},\boldsymbol\lambda_p) &= \mathbf 0, \end{aligned} \right. \qquad \star\in\{G,L\}, \label{eq:pressure-rom-full-interface-system} 
\end{equation*}
where \(\star=G\) denotes Galerkin projection and \(\star=L\) denotes the LSPG projection. In both cases, the interface equation is kept in the full interface space.

As in the full-order pressure step in Section \ref{subsubsec:dd-pressure}, we now derive the Newton--Schur system. At Newton iteration \(m\), the linearization of the  local pressure residual is
 \begin{equation*} 
 \widehat{\mathbf K}_{p,i}^{\star,(m)}\delta\mathbf a_i + \widehat{\mathbf M}_{p,i}^{\star,(m)} \mathbf H_i\delta\boldsymbol\lambda_p = -\widehat{\mathbf R}_{p,i}^{\star,(m)}, \qquad \star\in\{G,L\}. \label{eq:pressure-local-reduced-linear} 
 \end{equation*}
where 
\begin{equation*} 
\widehat{\mathbf K}_{p,i}^{G,(m)} = \boldsymbol\Phi_i^\top\mathbf K_{p,i}^{(m)}\boldsymbol\Phi_i, \qquad \widehat{\mathbf M}_{p,i}^{G,(m)} = \boldsymbol\Phi_i^\top\mathbf D_{p,i}^{(m)}, \qquad \widehat{\mathbf R}_{p,i}^{G,(m)} = \boldsymbol\Phi_i^\top\mathbf R_{p,i}^{(m)},
 \label{eq:pressure-galerkin-blocks} \end{equation*} 
 \begin{equation*}  
 \widehat{\mathbf K}_{p,i}^{L,(m)} = \left(\mathbf K_{p,i}^{(m)}\boldsymbol\Phi_i\right)^\top \left(\mathbf K_{p,i}^{(m)}\boldsymbol\Phi_i\right), \quad \widehat{\mathbf M}_{p,i}^{L,(m)} = \left(\mathbf K_{p,i}^{(m)}\boldsymbol\Phi_i\right)^\top \mathbf D_{p,i}^{(m)}, \quad \widehat{\mathbf R}_{p,i}^{L,(m)} = \left(\mathbf K_{p,i}^{(m)}\boldsymbol\Phi_i\right)^\top \mathbf R_{p,i}^{(m)},  
 \label{eq:pressure-lspg-blocks} 
 \end{equation*}
 and 
\begin{equation*} 
\mathbf D_{p,i}^{(m)} := \mathbf K_{p,i}^{(m)}\mathbf L_i + \widetilde{\mathbf M}_{p,i}.
 \label{eq:pressure-Dpi-definition} \end{equation*} 
 which denotes the derivative of the local residual with respect to the global interface-pressure increment, after including the interface lifting.

 Assuming that \(\widehat{\mathbf K}_{p,i}^{\star,(m)}\) is nonsingular, the local reduced increment can be eliminated: \begin{equation} \delta\mathbf a_i = \delta\mathbf a_{i,\mathrm{load}}^{\star,(m)} + \mathbf Y_{p,i}^{\star,(m)} \mathbf H_i\delta\boldsymbol\lambda_p, \label{eq:pressure-reduced-local-elimination} \end{equation} where \[ \delta\mathbf a_{i,\mathrm{load}}^{\star,(m)} = - \left( \widehat{\mathbf K}_{p,i}^{\star,(m)} \right)^{-1} \widehat{\mathbf R}_{p,i}^{\star,(m)}, \qquad \mathbf Y_{p,i}^{\star,(m)} = - \left( \widehat{\mathbf K}_{p,i}^{\star,(m)} \right)^{-1} \widehat{\mathbf M}_{p,i}^{\star,(m)} . \]

We then linearize the interface flux equation and obtain 
\begin{equation*} 
\sum_{i=1}^{N_s} \mathbf H_i^\top \left[ \mathbf A_{t,i}^{(m)} \left( \boldsymbol\Phi_i\delta\mathbf a_i + \mathbf L_i\mathbf H_i\delta\boldsymbol\lambda_p \right) + \mathbf B_{t,i}^{(m)} \mathbf H_i\delta\boldsymbol\lambda_p \right] = -\mathbf R_{p,\Gamma}^{(m)}, 
\label{eq:pressure-interface-linear-rom}
 \end{equation*}
where $ \mathbf A_{t,i}^{(m)}$ and 
$ B_{t,i}^{(m)}$ are the high-fidelity interface derivative blocks defined in  Section~\ref{subsubsec:dd-pressure}, evaluated at the reconstructed reduced state.

  Substituting \eqref{eq:pressure-reduced-local-elimination} into the linearized interface flux equation yields the reduced full-interface pressure Schur complement 
  \begin{equation} \mathbf S_{p,\Gamma}^{\star,(m)} \delta\boldsymbol\lambda_p = \mathbf g_{p,\Gamma}^{\star,(m)} , 
  \label{eq:pressure-reduced-schur} 
  \end{equation} 
  with 
  \begin{align*} 
  \mathbf S_{p,\Gamma}^{\star,(m)} &= \sum_{i=1}^{N_s} \mathbf H_i^\top \left[ \mathbf A_{t,i}^{(m)} \left( \boldsymbol\Phi_i\mathbf Y_{p,i}^{\star,(m)} + \mathbf L_i \right) \mathbf H_i + \mathbf B_{t,i}^{(m)}\mathbf H_i \right], \notag \\ \mathbf g_{p,\Gamma}^{\star,(m)} &= -\mathbf R_{p,\Gamma}^{(m)} - \sum_{i=1}^{N_s} \mathbf H_i^\top \mathbf A_{t,i}^{(m)} \boldsymbol\Phi_i \delta\mathbf a_{i,\mathrm{load}}^{\star,(m)} . 
  \label{eq:pressure-reduced-schur-blocks} \end{align*} 
  The Schur system \eqref{eq:pressure-reduced-schur} has the same interface dimension as the high-fidelity pressure Schur system, but the expensive local inverses are replaced by reduced local solves.
 
In the port-reduced case, the interface pressure trace is further approximated in a low-dimensional port space. 
We introduce a port basis \(\mathbf W\in\mathbb R^{N_\Gamma\times n_\Gamma}\) and approximate \[ \boldsymbol\lambda_p\approx\mathbf W\boldsymbol\mu, \qquad \mathbf W_i:=\mathbf H_i\mathbf W . \] This is the analogue of port-reduced static-condensation reduced-basis methods \cite{smetana2015certification,smetanaPatera2016}. The port-reduced pressure system is \begin{equation*} 
\left\{ \begin{aligned} \widehat{\mathbf R}_{p,i}^{\star} (\mathbf a_i,\mathbf W\boldsymbol\mu;\boldsymbol\alpha_{s,i}) &= \mathbf 0, && i=1,\ldots,N_s, \\ \mathbf W^\top\mathbf M_\Gamma \mathbf R_{p,\Gamma} (\{\mathbf X_i^r\},\mathbf W\boldsymbol\mu) &= \mathbf 0, \end{aligned} \right. 
\label{eq:pressure-port-rom-system} 
\end{equation*} 
where \(\mathbf M_\Gamma\) is an interface mass matrix. At the Newton level, this corresponds to projecting the full-interface Schur equation: \begin{equation*} 
\mathbf W^\top\mathbf M_\Gamma \mathbf S_{p,\Gamma}^{\star,(m)} \mathbf W\,\delta\boldsymbol\mu = \mathbf W^\top\mathbf M_\Gamma \mathbf g_{p,\Gamma}^{\star,(m)} . 
\label{eq:pressure-port-schur} 
\end{equation*} 
The port-reduced formulation decreases the interface dimension from \(N_\Gamma\) to \(n_\Gamma\). However, unlike the full-interface formulation, it enforces only the projected interface flux balance. Hence edgewise interface flux continuity is generally not guaranteed, unless the port space is sufficiently rich to represent the relevant interface-pressure traces. For this reason, the port-reduced variant is appropriate when the offline data are expected to cover the online interface manifold. In moving-well prediction, where new interface pressure traces may appear outside the training manifold, we primarily use full-interface coupling.

\subsection{Transport-Step Reduced Schur System}
\label{subsec:transport-rom}

We now derive the reduced transport step. As discussed in Section~\ref{subsubsec:dd-transport}, the interface saturation trace has a different role from the interface pressure trace. The pressure trace is an elliptic interface variable required by all adjacent local pressure solves. By contrast, the saturation trace is an upwind inflow quantity: it enters a local transport residual only on interface faces where the total velocity points into the subdomain. This makes a port-reduced saturation trace less natural than a port-reduced pressure trace, because the active inflow set may depend on the velocity field, time, and parameter. We therefore retain the saturation trace in the full interface space and reduce only the local saturation unknowns.

The reduced transport system is
\begin{equation*}
\left\{
\begin{aligned}
\mathbf{Z}_{s,i}^T
\mathbf{R}_{s,i}
(\mathbf{Z}_{s,i}\boldsymbol{\alpha}_{s,i},
\mathbf{H}_i\boldsymbol{\lambda}_s;\mathbf{X}_i)
&=\mathbf{0},
&&i=1,\ldots,N_s,\\
\mathbf{R}_{s,\Gamma}
(\{\boldsymbol{\alpha}_{s,i}\},\boldsymbol{\lambda}_s)&=\mathbf{0}.
\end{aligned}
\right.
\label{eq:transport-rom-compact}
\end{equation*}
The first equation is the Galerkin projection of the local transport residual onto the saturation basis. The second equation is the unprojected interface \(\mathrm{CO}_2\)-flux continuity condition. It is evaluated using the reconstructed saturations \(\mathbf S_i^r\) and the same interface flux definition as in \eqref{eq:transport-dd-compact}. Thus, as in the full-interface pressure ROM, the interface residual is kept in the original trace space.

In our numerical experiments, this Galerkin transport ROM provides the same order of accuracy as the corresponding LSPG transport update, while avoiding the additional normal-equation solves required by LSPG. We therefore use Galerkin projection for the transport step.

At Newton iteration \(m\), the reduced local transport linearization is 
\begin{equation*} 
\widehat{\mathbf J}_{s,i}^{(m)} \delta\boldsymbol\alpha_{s,i} + \widehat{\mathbf G}_{s,i}^{(m)} \mathbf H_i\delta\boldsymbol\lambda_s = -\widehat{\mathbf R}_{s,i}^{(m)}, \qquad i=1,\ldots,N_s, 
\label{eq:transport-rom-local-linear}
 \end{equation*} 
where \begin{equation*} 
\widehat{\mathbf J}_{s,i}^{(m)} = \mathbf Z_{s,i}^\top \mathbf J_{s,i}^{(m)} \mathbf Z_{s,i}, \qquad \widehat{\mathbf G}_{s,i}^{(m)} = \mathbf Z_{s,i}^\top \mathbf G_{s,i}^{(m)}, \qquad \widehat{\mathbf R}_{s,i}^{(m)} = \mathbf Z_{s,i}^\top \mathbf R_{s,i}^{(m)} .
 \label{eq:transport-rom-local-blocks} 
 \end{equation*} 
 Here \(\mathbf J_{s,i}^{(m)}\) and \(\mathbf G_{s,i}^{(m)}\) are the local transport derivative blocks defined in Section~\ref{subsubsec:dd-transport}, evaluated at the reconstructed reduced saturation and the frozen pressure-step state. Because the interface saturation is used only on inflow faces, \(\mathbf G_{s,i}^{(m)}\) has nonzero contributions only from cells adjacent to inflow interface faces. Assuming that \(\widehat{\mathbf J}_{s,i}^{(m)}\) is nonsingular, the local reduced saturation increment can be eliminated: \begin{equation} \delta\boldsymbol\alpha_{s,i} = \delta\boldsymbol\alpha_{s,i,\mathrm{load}}^{(m)} + \mathbf Y_{s,i}^{(m)} \mathbf H_i\delta\boldsymbol\lambda_s, \label{eq:transport-rom-local-elimination} \end{equation} with \[ \delta\boldsymbol\alpha_{s,i,\mathrm{load}}^{(m)} = - \left( \widehat{\mathbf J}_{s,i}^{(m)} \right)^{-1} \widehat{\mathbf R}_{s,i}^{(m)}, \qquad \mathbf Y_{s,i}^{(m)} = - \left( \widehat{\mathbf J}_{s,i}^{(m)} \right)^{-1} \widehat{\mathbf G}_{s,i}^{(m)} . \] The interface residual is then linearized using the full-interface \(\mathrm{CO}_2\)-flux derivatives from Section~\ref{subsubsec:dd-transport}, evaluated at the reconstructed reduced state: 
\[ \mathbf A_{c,i}^{(m)} := \left. \frac{\partial \mathbf F_{c,i}}{\partial \mathbf S_i} \right|_{(m)}, \qquad \mathbf B_{c,i}^{(m)} := \left. \frac{\partial \mathbf F_{c,i}}{\partial\boldsymbol\lambda_{s,i}} \right|_{(m)} . 
\] 
Since 
\[ \delta\mathbf S_i = \mathbf Z_{s,i}\delta\boldsymbol\alpha_{s,i}, \qquad \delta\boldsymbol\lambda_{s,i} = \mathbf H_i\delta\boldsymbol\lambda_s, 
\] 
the linearized interface equation is 
\begin{equation} 
\sum_{i=1}^{N_s} \mathbf H_i^\top \left( \mathbf A_{c,i}^{(m)} \mathbf Z_{s,i}\delta\boldsymbol\alpha_{s,i} + \mathbf B_{c,i}^{(m)} \mathbf H_i\delta\boldsymbol\lambda_s \right) = -\mathbf R_{s,\Gamma}^{(m)} . 
\label{eq:transport-rom-interface-linear} 
\end{equation} 
Substituting \eqref{eq:transport-rom-local-elimination} into \eqref{eq:transport-rom-interface-linear} gives the reduced transport Schur complement \begin{equation} 
\mathbf S_{s,\Gamma}^{\mathrm{G},(m)} \delta\boldsymbol\lambda_s = \mathbf g_{s,\Gamma}^{\mathrm{G},(m)}, 
\label{eq:transport-rom-schur}
 \end{equation} 
 where 
 \begin{align*} 
 \mathbf S_{s,\Gamma}^{\mathrm{G},(m)} &= \sum_{i=1}^{N_s} \mathbf H_i^\top \left[ \mathbf A_{c,i}^{(m)} \mathbf Z_{s,i} \mathbf Y_{s,i}^{(m)} \mathbf H_i + \mathbf B_{c,i}^{(m)} \mathbf H_i \right], \notag \\ 
 \mathbf g_{s,\Gamma}^{\mathrm{G},(m)} &= -\mathbf R_{s,\Gamma}^{(m)} - \sum_{i=1}^{N_s} \mathbf H_i^\top \mathbf A_{c,i}^{(m)} \mathbf Z_{s,i} \delta\boldsymbol\alpha_{s,i,\mathrm{load}}^{(m)} .
  \label{eq:transport-rom-schur-blocks} 
  \end{align*} 
  The Schur system \eqref{eq:transport-rom-schur} has the same interface dimension as the high-fidelity transport Schur system, but the local cellwise saturation updates are performed in reduced coordinates. Since the interface residual is not projected, the full-interface \(\mathrm{CO}_2\)-flux balance is preserved at the reduced level.

For later reference, we label the DD-ROM variants according to the two choices made in the pressure step: the local pressure residual reduction and the interface-pressure treatment. The local pressure equation is reduced either by Galerkin projection or by LSPG projection, and the pressure interface trace is either retained in the full-interface space or approximated in a port-reduced space. This gives four variants:
\[
\mathrm{Galerkin\text{-}FI},\qquad
\mathrm{LSPG\text{-}FI},\qquad
\mathrm{Galerkin\text{-}PR},\qquad
\mathrm{LSPG\text{-}PR}.
\]
In all four variants, the transport step is treated by the Galerkin saturation ROM with full-interface saturation coupling described above.

\section{Well-Aware Basis Construction and Hyper-Reduction} 
\label{sec:well-aware-compact}

\begin{figure}[H]
\centering
\begin{tikzpicture}[
    font=\small,
    >=Latex,
    gridline/.style={draw=gray!55,line width=0.35pt},
    source/.style={fill=black,draw=black},
    near/.style={fill=red!13,draw=red!65,dashed,line width=0.8pt},
    far/.style={fill=blue!7,draw=blue!55,dashed,line width=0.8pt},
    arrow/.style={-Latex,thick}
]
\begin{scope}[shift={(0,0)}]
    \node[align=center] at (1.5,3.48) {offline training\\fixed well};
    \draw[thick] (0,0) rectangle (3.0,3.0);
    \foreach \x in {0.6,1.2,1.8,2.4} {\draw[gridline] (\x,0) -- (\x,3.0);}
    \foreach \y in {0.6,1.2,1.8,2.4} {\draw[gridline] (0,\y) -- (3.0,\y);}
    \fill[red!13] (1.2,1.2) rectangle (1.8,1.8);
    \draw[red!65,dashed,line width=0.8pt] (1.2,1.2) rectangle (1.8,1.8);
    \fill[source] (1.40,1.40) rectangle (1.60,1.60);
    \node[white,font=\scriptsize] at (1.50,1.50) {$q$};
    \node[align=center,font=\scriptsize] at (1.5,-0.32) {well-local snapshots};
\end{scope}

\draw[arrow] (3.55,1.55) -- node[above,align=center,font=\scriptsize] {local POD\\HAPOD} (4.95,1.55);

\begin{scope}[shift={(5.5,0)}]
    \node[align=center] at (1.5,3.48) {online prediction\\well moved};
    \fill[blue!7] (0,0) rectangle (3.0,3.0);
    \draw[thick] (0,0) rectangle (3.0,3.0);
    \foreach \x in {0.6,1.2,1.8,2.4} {\draw[gridline] (\x,0) -- (\x,3.0);}
    \foreach \y in {0.6,1.2,1.8,2.4} {\draw[gridline] (0,\y) -- (3.0,\y);}
    \draw[blue!55,dashed,line width=0.8pt] (0,0) rectangle (3.0,3.0);
    \fill[red!13] (0.0,1.2) rectangle (0.6,1.8);
    \fill[source] (0.20,1.40) rectangle (0.40,1.60);
    \node[white,font=\scriptsize] at (0.30,1.50) {$q$};
    \node[align=center,font=\scriptsize] at (1.5,-0.32) {FOM near new well, ROM elsewhere};
\end{scope}

\node[draw=red!65,fill=red!13,minimum width=0.35cm,minimum height=0.18cm] at (3.15,-0.82) {};
\node[anchor=west,font=\scriptsize] at (3.42,-0.82) {well subdomain};
\node[draw=blue!55,dashed,fill=blue!7,minimum width=0.35cm,minimum height=0.18cm] at (3.15,-1.15) {};
\node[anchor=west,font=\scriptsize] at (3.42,-1.15) {far-field DD-ROM region};
\end{tikzpicture}
\caption{Schematic of the moving-well setting on a $5\times5$ subdomain partition. Offline snapshots for one well location can be compressed by local POD and HAPOD to form far-field bases. When the source moves, the target-well source response is appended as an enrichment mode and the new well region may be kept at high fidelity or hybrid resolution, while DD-ROM bases are used in the far field.}
\label{fig:moving-well-hapod-schematic}
\end{figure}
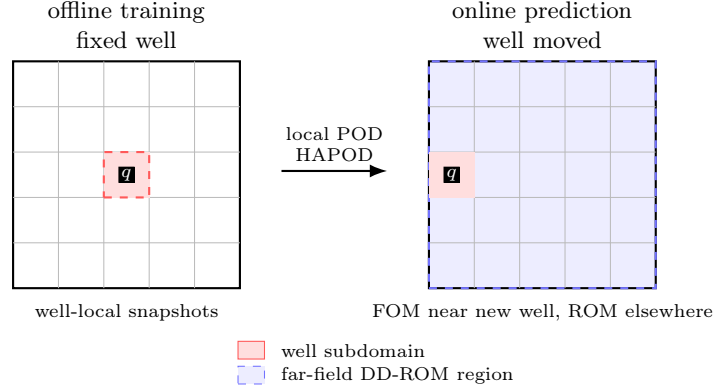

The DD-ROM formulation developed in Section~\ref{sec:ddrom-compact} is effective when the online pressure, saturation, and interface traces remain well represented by the offline snapshot spaces. 
 This assumption is reasonable for parameter variations that perturb the solution within the same spatial regions sampled during training, but it becomes more restrictive when the well location changes.
 In that case, the source term is displaced to a different subdomain, and the associated near-well pressure response, local flux pattern, and saturation front may appear in regions that were only weakly activated in the offline data. 
This observation motivates basis enrichment, but the enrichment should be localized. A uniform enrichment of all subdomain bases would be unnecessarily expensive, since the dominant loss of accuracy is typically localized around the new source. 
We therefore introduce the well-aware strategy illustrated in Figure~\ref{fig:moving-well-hapod-schematic}. The offline stage uses local POD and HAPOD to compress well-free or far-field information into reusable subdomain bases. When the well moves online, source-response enrichment adds target well pressure and velocity 
directions to the affected local spaces. If the new near-well response is still not reliably captured by the enriched reduced basis, the well and neighboring subdomains are treated by a hybrid FOM--ROM coupling, while the far-field region remains reduced. 

This section first describes the source-response enriched HAPOD basis construction and the hybrid coupling strategy, and then presents the conservation-aware hyper-reduction procedure used to reduce nonlinear residual assembly costs without sampling interface fluxes, boundary terms, or well/source contributions. 

\subsection{Source-response enrichment and HAPOD compression}
\label{subsec:well-aware-basis}

Let \(\xi\) denote an online well configuration. We denote by
\(\mathcal I_{\mathrm{well}}(\xi)\) the set of subdomains that contain 
the active well and by \(\mathcal I_{\mathrm{near}}(\xi)\) a prescribed
near-well neighborhood, for example the first layer of subdomains adjacent to \(\mathcal I_{\mathrm{well}}(\xi)\). The remaining subdomains form the far-field set
\[
\mathcal I_{\mathrm{far}}(\xi)
:=
\{1,\ldots,N_s\}
\setminus
\left(
\mathcal I_{\mathrm{well}}(\xi)
\cup
\mathcal I_{\mathrm{near}}(\xi)
\right).
\]
The far-field solution is typically less sensitive to the precise well
position and can often be represented by bases trained from previous
well configurations. By contrast, the well and near-well subdomains may
require additional enrichment or high-fidelity treatment. 
This localization principle is consistent with localized reduced basis and arbitrary local modification strategies, where the effect of a local change is isolated to a subset of components so that global retraining can be avoided 
 \cite{kaulmann2015lrbms,ohlbergerSchindler2015,buhrEngwerOhlbergerRave2017,buhrSmetana2018}.

Following the local snapshot preprocessing described in Section~\ref{subsec:local-spaces-lifting}, we first construct compact far-field bases from high-fidelity DD-SFI snapshots. For each variable \(v\in\{\mathbf{u}_t,p,S\}\) and each subdomain \(\Omega_i\), the velocity and pressure snapshots are taken after subtracting the interface lifting, whereas the saturation snapshots are used directly. A local POD compression with tolerance \(\varepsilon_{\mathrm{loc}}\) gives preliminary bases
\[
\left[
\mathbf Z_{v,i}^{\mathrm{loc}},
\boldsymbol\sigma_{v,i}
\right]
=
\mathrm{POD}
\left(
\{\mathbf V_{v,i}^k\}_{k=1}^{n_{\mathrm{train}}},
(\cdot,\cdot)_{\mathbf X_v},
\varepsilon_{\mathrm{loc}}
\right),
\]
where \(\mathbf V_{v,i}^k\) denotes the corresponding preprocessed local snapshot and \((\cdot,\cdot)_{\mathbf X_v}\) is the inner product used for that variable.
To aggregate information from multiple
parameters or multiple far-field subdomains, we apply a second HAPOD \cite{himpe2018hapod} 
compression to the singular-value-weighted local modes:
\[
\left[
\mathbf Z_{v}^{\mathrm{far}},
\boldsymbol\sigma_{v}^{\mathrm{far}}
\right]
=
\mathrm{POD}
\left(
\left[
\mathbf Z_{v,i_1}^{\mathrm{loc}}\boldsymbol\Sigma_{v,i_1},
\ldots,
\mathbf Z_{v,i_m}^{\mathrm{loc}}\boldsymbol\Sigma_{v,i_m}
\right],
(\cdot,\cdot)_{\mathbf X_v},
\varepsilon_{\mathrm{root}}
\right).
\]

For subdomains affected by the target well, we enrich the reduced spaces
using a local source response. Let \(\mathbf q_{t,i}(\xi)\) denote the local total mass source vector on \(\Omega_i\) for the target well
configuration. If \(\mathbf q_{t,i}(\xi)\neq\mathbf 0\), we compute a
local mixed Darcy source response
\[
\begin{bmatrix}
\mathbf M_{u,i} & -\mathbf C_{d,i}\\
\mathbf C_{d,i}^{\top} & \mathbf 0
\end{bmatrix}
\begin{bmatrix}
\mathbf U_{q,i}\\
\mathbf P_{q,i}
\end{bmatrix}
=
\begin{bmatrix}
\mathbf 0\\
\mathbf D_{c,i}\mathbf q_{t,i}(\xi)
\end{bmatrix}.
\]
The algebraic form is the same local RT0--P0 saddle-point structure as in the interface-lifting problem of Section~\ref{subsec:local-spaces-lifting}, but with a source term in the mass equation instead of nonhomogeneous interface-pressure data. 
Here \(\mathbf M_{u,i}\) is the local RT0 mass matrix, \(\mathbf C_{d,i}\) is the local divergence matrix, and \(\mathbf D_{c,i}\) is the cell-measure matrix. The resulting pair \((\mathbf U_{q,i},\mathbf P_{q,i})\) is then orthogonalized against the existing local bases and appended as an enrichment direction:
\[
\mathbf Z_{u,i}
\leftarrow
\mathrm{orth}_{\mathbf X_u}
\left[
\mathbf Z_{u,i},\mathbf U_{q,i}
\right],
\qquad
\mathbf Z_{p,i}
\leftarrow
\mathrm{orth}_{\mathbf X_p}
\left[
\mathbf Z_{p,i},\mathbf P_{q,i}
\right].
\]
If \(\mathbf q_{t,i}(\xi)=\mathbf 0\), no source-response enrichment is added on \(\Omega_i\).

For the moving-well tests in Section~\ref{sec:num_case3}, we use full-interface pressure coupling. A port basis trained from the original well configuration may not represent the interface-pressure traces induced by a relocated source, especially near the new well. Full-interface coupling avoids this additional approximation by retaining all interface pressure degrees of freedom and enforcing flux balance in the original trace space. Developing adaptive or well-aware port bases for moving wells is left for future work.

\subsection{Hybrid FOM--ROM Coupling}
\label{subsec:hybrid-coupling}

Source-response enrichment improves the local approximation space, but near-well extrapolation can still be severe when the source is moved far from the training configuration. We therefore allow selected subdomains to remain at high fidelity during the online solve. Let \(\mathcal I_{\mathrm{FOM}}(\xi)\) denote the well and near-well subdomains treated by the high-fidelity model, and define the reduced
far-field region by
\[
\mathcal I_{\mathrm{ROM}}(\xi)
:=
\{1,\ldots,N_s\}
\setminus
\mathcal I_{\mathrm{FOM}}(\xi).
\]
In the numerical experiments, \(\mathcal I_{\mathrm{FOM}}(\xi)\) can be
chosen as the active well subdomain, and, when needed, a prescribed near-well neighborhood, such as one layer of adjacent subdomains.

The hybrid construction is natural in the present DD formulation, because all subdomains are coupled only through the interface residuals. 
Thus a subdomain may be treated either at high fidelity or at reduced order without changing the global interface unknowns. We describe this hybrid coupling only for the pressure step, and the transport step follows similarly. 
Using the notations in Sections \ref{sec:sfi-dd} and \ref{sec:ddrom-compact}, The hybrid pressure problem is therefore
\begin{equation}
\left\{
\begin{aligned}
\mathbf R_{p,i}
\left(
\mathbf X_i,\mathbf H_i\boldsymbol\lambda_p;\mathbf S_i
\right)
&=
\mathbf 0,
&& i\in\mathcal I_{\mathrm{FOM}},
\\
\widehat{\mathbf R}_{p,i}^{\star}
\left(
\mathbf a_i,\boldsymbol\lambda_p;\boldsymbol\alpha_{s,i}
\right)
&=
\mathbf 0,
&& i\in\mathcal I_{\mathrm{ROM}},
\\
\mathbf R_{p,\Gamma}
\left(
\{\mathbf X_i^{\mathrm{hyb}}\},\boldsymbol\lambda_p
\right)
&=
\mathbf 0,
\end{aligned}
\right.
\label{eq:hybrid-pressure-system}
\end{equation}
where
\[
\mathbf X_i^{\mathrm{hyb}}
=
\begin{cases}
\mathbf X_i, & i\in\mathcal I_{\mathrm{FOM}},\\
\mathbf X_i^r, & i\in\mathcal I_{\mathrm{ROM}}.
\end{cases}
\]
The interface equation in \eqref{eq:hybrid-pressure-system} is the same DD interface residual used in the high-fidelity formulation and in the DD-ROM. This is what makes the hybrid construction straightforward: the model fidelity is chosen locally, while the global coupling is still imposed through the same interface variables. 
The resulting nonlinear system \eqref{eq:hybrid-pressure-system} is solved by the same Newton--Schur strategy as in Sections~\ref{subsubsec:dd-pressure} and~\ref{subsec:pressure-rom}; only the local elimination differs depending on whether a subdomain is treated at high fidelity or reduced order. We therefore omit the algebraic details.

\subsection{Conservation-aware hyper-reduction}
\label{subsec:hyperreduction}

In projection-based ROMs, the number of unknowns is reduced, but the online evaluation of nonlinear residuals and Jacobians can still scale with the high-fidelity dimension. To address this bottleneck, we introduce a local empirical-quadrature hyper-reduction procedure for the interior residual contributions. The construction follows the residual-minimization viewpoint of GNAT/LSPG \cite{carlberg2013gnat,carlberg2017galerkin}, while using nonnegative sampling weights in the spirit of ECSW and empirical cubature \cite{farhat2015ecsw,hernandez2017empiricalCubature}. The resulting sample-and-weight approximation is applied locally within each subdomain, whereas the interface residuals are kept in the original interface space. This component-wise organization is consistent with hyper-reduced reduced-basis-element methods, in which each component is equipped with its own reduced space, sampled operator, and quadrature rule before being coupled through interface variables \cite{ebrahimi2024hyperreduced}. 

We describe the algebraic form for a generic local subproblem
\(q\in\{p,s\}\). Let \(\mathbf R_{q,i}\) be the local residual,
\(\mathbf J_{q,i}\) its Jacobian with respect to the local state, and
\(\mathbf Z_{q,i}\) the corresponding trial basis, with
\[
\mathbf Z_{p,i}:=\boldsymbol\Phi_i,
\qquad
\mathbf Z_{s,i}:=\mathbf Z_{s,i}.
\]
Let \(\mathbf B_{q,i}\) denote the derivative of the local residual with
respect to the local interface trace after all lifting terms are
included. Thus,
\[
\mathbf B_{p,i}:=\mathbf D_{p,i},
\qquad
\mathbf B_{s,i}:=\mathbf G_{s,i}.
\]
We also introduce a test operator \(\mathbf V_{q,i}\). For pressure
Galerkin, \(\mathbf V_{p,i}=\boldsymbol\Phi_i\); for pressure LSPG,
\(\mathbf V_{p,i}=\mathbf K_{p,i}\boldsymbol\Phi_i\); for transport
Galerkin, \(\mathbf V_{s,i}=\mathbf Z_{s,i}\).

Let \(\mathcal I_{q,i}\) be a set of sampled local residual rows or row
blocks, let \(\boldsymbol\Pi_{\mathcal I_{q,i}}\) be the associated row
restriction, and let
\[
\mathbf W_{q,i}
=
\mathrm{diag}(\rho_\ell)_{\ell\in\mathcal I_{q,i}},
\qquad
\rho_\ell>0,
\]
be the diagonal matrix of empirical-quadrature weights. The
hyper-reduced local blocks are
\begin{align}
\widehat{\mathbf K}_{q,i}^{\mathrm{hr}}
&=
\left(
\boldsymbol\Pi_{\mathcal I_{q,i}}\mathbf V_{q,i}
\right)^{\top}
\mathbf W_{q,i}
\left(
\boldsymbol\Pi_{\mathcal I_{q,i}}\mathbf J_{q,i}\mathbf Z_{q,i}
\right),
\notag\\
\widehat{\mathbf M}_{q,i}^{\mathrm{hr}}
&=
\left(
\boldsymbol\Pi_{\mathcal I_{q,i}}\mathbf V_{q,i}
\right)^{\top}
\mathbf W_{q,i}
\left(
\boldsymbol\Pi_{\mathcal I_{q,i}}\mathbf B_{q,i}
\right),
\notag\\
\widehat{\mathbf R}_{q,i}^{\mathrm{hr}}
&=
\left(
\boldsymbol\Pi_{\mathcal I_{q,i}}\mathbf V_{q,i}
\right)^{\top}
\mathbf W_{q,i}
\left(
\boldsymbol\Pi_{\mathcal I_{q,i}}\mathbf R_{q,i}
\right).
\label{eq:hr-local-blocks}
\end{align}
These blocks replace the corresponding full reduced products in the
local elimination step. The interface residuals
\(\mathbf R_{p,\Gamma}\) and \(\mathbf R_{s,\Gamma}\), however, are not
sampled and are assembled in the full interface space.

The sampling sets are chosen according to the RT0--P0 structure. In the pressure step, the local residual contains RT0 momentum equations on faces and P0 mass-balance
equations on cells. We therefore write
\[
\mathcal I_{p,i}
=
\mathcal I_{u,i}^{p}
\cup
\left(
N_{f,i}+\mathcal I_{m,i}^{p}
\right),
\]
where \(N_{f,i}\) is the number of local face equations,
\(\mathcal I_{u,i}^{p}\) samples face/momentum rows, and
\(\mathcal I_{m,i}^{p}\) samples cell mass-balance rows. In the
transport step, the primary sample is a set of cells
\(\mathcal I_{s,i}\). Evaluating the residuals on these cells requires
the fluxes on all adjacent faces, so we use the induced edge closure
\[
\mathcal E(\mathcal I_{s,i})
=
\{e:\ e\subset\partial K,\ K\in\mathcal I_{s,i}\}.
\]

The empirical weights are fitted offline from training residual and
Jacobian data. For a candidate block \(\ell\) and a training state
indexed by \(m\), let
\[
\mathbf K_{\ell}^{(m)}
=
\left(\mathbf V_{\ell}^{(m)}\right)^{\top}
\left(\mathbf J_{\ell}^{(m)}\mathbf Z_{q,i}\right),
\qquad
\mathbf M_{\ell}^{(m)}
=
\left(\mathbf V_{\ell}^{(m)}\right)^{\top}
\mathbf B_{\ell}^{(m)},
\qquad
\mathbf r_{\ell}^{(m)}
=
\left(\mathbf V_{\ell}^{(m)}\right)^{\top}
\mathbf R_{\ell}^{(m)} .
\]
The candidate moment vector collects the contributions needed to
reconstruct the reduced blocks,
\[
\mathbf y_{\ell}^{(m)}
=
\begin{bmatrix}
\mathrm{vec}\left(\mathbf K_{\ell}^{(m)}\right)\\
\mathrm{vec}\left(\mathbf M_{\ell}^{(m)}\right)\\
\mathbf r_{\ell}^{(m)}\\
1
\end{bmatrix}.
\]
After aggregating over training states,
\[
\overline{\mathbf y}_{\ell}
=
\sum_m \omega_m\mathbf y_{\ell}^{(m)},
\]
the full moment target is
\[
\mathbf b_{q,i}
=
\sum_{\ell\in\mathcal F_{q,i}}
\overline{\mathbf y}_{\ell},
\]
where \(\mathcal F_{q,i}\) is the full candidate set. For a selected
sample \(\mathcal I_{q,i}=\{\ell_1,\ldots,\ell_r\}\), define
\[
\mathbf Y_{q,i,\mathcal I}
=
\left[
\overline{\mathbf y}_{\ell_1},
\ldots,
\overline{\mathbf y}_{\ell_r}
\right].
\]
The nonnegative empirical-quadrature weights are obtained from
\[
\boldsymbol\rho_{q,i}
=
\arg\min_{\boldsymbol\rho\ge \mathbf 0}
\left\|
\begin{bmatrix}
\mathbf Y_{q,i,\mathcal I}\\
\sqrt{\eta}\mathbf I
\end{bmatrix}
\boldsymbol\rho
-
\begin{bmatrix}
\mathbf b_{q,i}\\
\sqrt{\eta}\mathbf 1
\end{bmatrix}
\right\|_2^2 ,
\]
where \(\eta>0\) is a small regularization parameter. 

We emphasize that hyper-reduction is applied only to the projected local interior residual contributions used in the reduced local elimination. Interface flux residuals, physical boundary terms, and well/source contributions are excluded from the empirical-quadrature approximation and are assembled with their original weights. In particular, in the full-interface variants, the interface residual is neither sampled nor projected, so edgewise interface flux balance is retained at the algebraic level. Preserving the well/source terms exactly is also important for moving-well simulations, since the localized forcing near the active well should not be distorted by the sampling procedure.

\section{Numerical Results}
\label{sec:numerical}

We assess the proposed DD-ROM framework on three two-phase flow test cases. Section \ref{sec:num_case1} uses a manufactured solution to verify the high-fidelity RT0--P0/SFI discretization and to compare the DD-ROM variants. Section \ref{sec:num_case2} studies out-of-sample parametric prediction, first for injection-rate variation and then for permeability variation. Section \ref{sec:num_case3} considers fixed- and moving-well configurations and evaluates the robustness of the HAPOD, and hybrid FOM--ROM strategies.

\noindent\textbf{Error metric.}
All ROM errors are measured relative to the high-fidelity DD solution computed on the same mesh. 
For a variable \(q\in\{p,S\}\), let \(\mathbf q_i^{\mathrm{hf},(k)}\) and \(\mathbf q_i^{\mathrm{rom},(k)}\) denote the high-fidelity and reduced coefficient vectors on subdomain \(\Omega_i\) at time \(t^k\). 
We define the relative error by
\begin{equation}
E_q^{\mathrm{ROM},(k)}
=
\left(
\frac{
\sum_{i=1}^{N_s}
\left\|
\mathbf q_i^{\mathrm{rom},(k)}
-
\mathbf q_i^{\mathrm{hf},(k)}
\right\|_{0,\Omega_i}^2
}{
\sum_{i=1}^{N_s}
\left\|
\mathbf q_i^{\mathrm{hf},(k)}
\right\|_{0,\Omega_i}^2
}
\right)^{1/2},
\qquad
q\in\{p,S\}.
\label{eq:rom-relative-errors}
\end{equation}
Here \(q=p\) corresponds to pressure and \(q=S\) to gas saturation. 
For the manufactured-solution verification in Section~\ref{sec:num_case1}, analogous relative \(L^2\) errors with respect to the exact solution are also reported to document the accuracy of the reference high-fidelity solver.

\subsection{Test case 1: manufactured solution}
\label{sec:num_case1}

We consider a manufactured solution problem on the dimensionless square $\Omega=[0,1]\times[0,1]$. The domain is divided into a $3\times3$ non-overlapping partition, and each subdomain is discretized by a uniform Cartesian mesh with 800 triangular cells. We use time step $\Delta t=10^{-2}$ and integrate to $t=10$. We use the following constitutive relations for the capillary pressure \cite[Chapter~8]{riviere2008discontinuous} and the relative permeability \cite{hoteit2008numerical}
\[
  P_c=-0.1(1-S_g)^{-1/2},\qquad
  k_{ra}=(1-S_g)^2,\qquad
  k_{rg}=S_g^2 .
\]
The permeability is isotropic, $\mathbf{K}=\mathbf{I}$. The remaining dimensionless parameters are $\phi=0.2$, $\mu_a=1$, $\mu_g=0.8$, $\rho_a=1.1$, $\rho_g=1$, $X_a^{\mathrm{CO}_2}=0.1$, and $X_g^{\mathrm{CO}_2}=0.8$.

The analytical expressions of pressure $p$ and saturation $S_g$ are prescribed as 
\begin{align*}
    p&=(1 + 0.2t)\sin{\pi x} \sin{\pi y},\\
    S_g&= 0.5 + 0.025t \sin{\pi x} \sin{\pi y}.
\end{align*}
The source terms \(Q_t\) and \(Q_{\mathrm{CO}_2}\), together with the total velocity \(\mathbf u_t\), are obtained by substituting these fields into the governing equations.  Boundary and initial data are imposed consistently from these analytical fields. 

\begin{figure}[H]
    \centering
    \subfloat[$t=2.5$]{
        \includegraphics[width=0.225\textwidth]{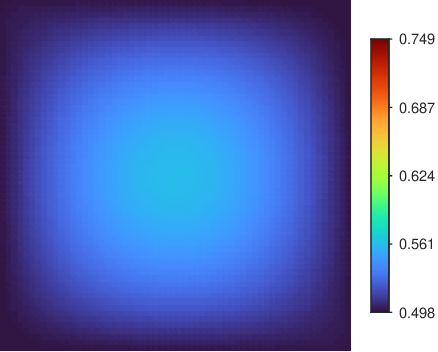}
    }\hspace{0.005\textwidth}
    \subfloat[$t=5.0$]{
        \includegraphics[width=0.225\textwidth]{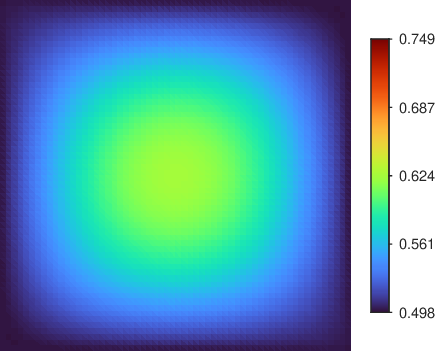}
    }\hspace{0.005\textwidth}
    \subfloat[$t=7.5$]{
        \includegraphics[width=0.225\textwidth]{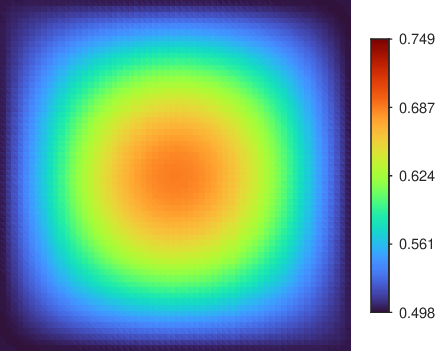}
    }\hspace{0.005\textwidth}
    \subfloat[$t=10.0$]{
        \includegraphics[width=0.225\textwidth]{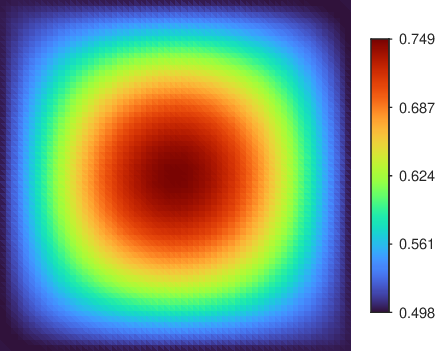}
    }\\
    \vspace{0.05cm}
    \subfloat[$t=2.5$]{
        \includegraphics[width=0.225\textwidth]{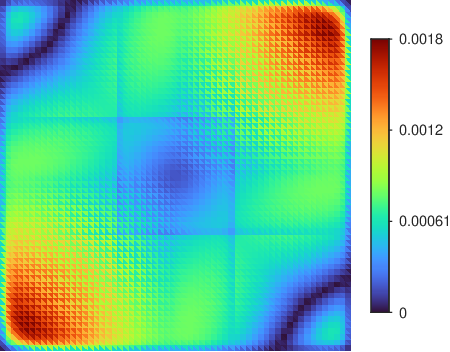}
    }\hspace{0.005\textwidth}
    \subfloat[$t=5.0$]{
        \includegraphics[width=0.225\textwidth]{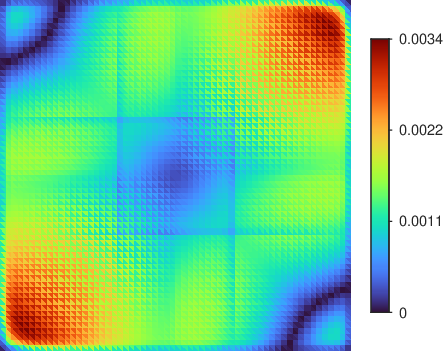}
    }\hspace{0.005\textwidth}
    \subfloat[$t=7.5$]{
        \includegraphics[width=0.225\textwidth]{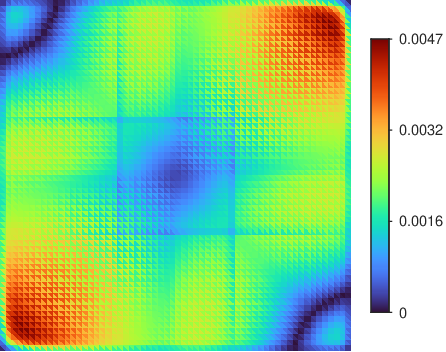}
    }\hspace{0.005\textwidth}
    \subfloat[$t=10.0$]{
        \includegraphics[width=0.225\textwidth]{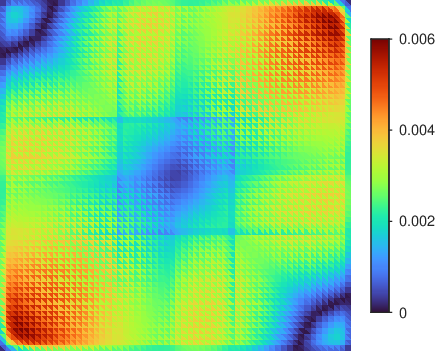}
    }
    \caption{Manufactured solution case: high-fidelity saturation. Top row: $S_g^{\mathrm{hf}}$ . Bottom row: absolute analytical error $|S_g^{\mathrm{hf}}-S_g^{\mathrm{ex}}|$. }
    \label{fig:case1-fom-saturation-time-error}
\end{figure}

We first verify the RT0--P0/SFI high-fidelity solver against the manufactured solution. 
Figure~\ref{fig:case1-fom-saturation-time-error} shows the computed saturation field at four representative times, together with the corresponding absolute error with respect to the exact solution. 
The saturation remains smooth and preserves the prescribed sinusoidal spatial structure, while its amplitude increases gradually in time according to the manufactured expression. 
The error fields remain small over the full simulation interval, confirming that the high-fidelity solver provides an accurate reference for the subsequent ROM training and evaluation.

\begin{figure}[H]
    \centering
    \subfloat[Pressure\label{fig:FOM-ex errline P}]{
        \includegraphics[width=0.48\textwidth]{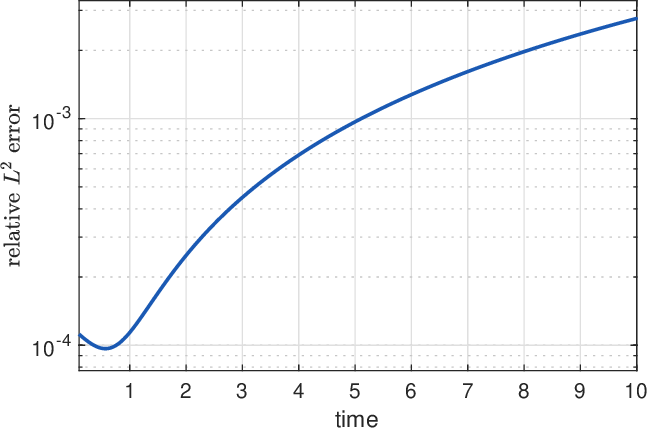}
    }
    \hfill
    \subfloat[Saturation\label{fig:FOM-ex errline S}]{
        \includegraphics[width=0.48\textwidth]{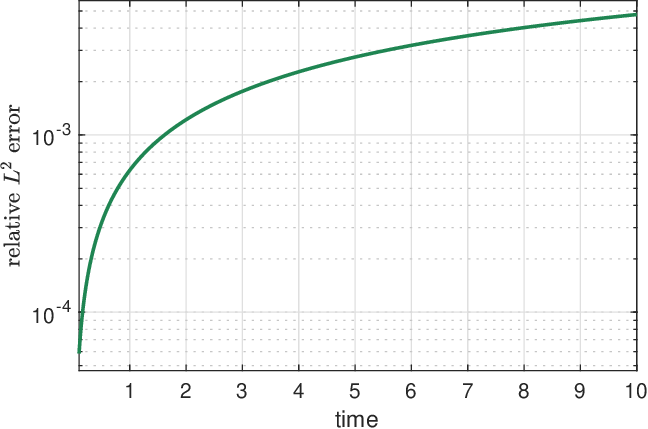}
    }
    \caption{Relative $L^2$ high-fidelity errors with respect to the analytical manufactured solution: (a) pressure and (b) gas saturation.}
    \label{fig:FOM-ex_errline}
\end{figure}

Figure~\ref{fig:FOM-ex_errline} reports the high-fidelity pressure and saturation errors with respect to the exact solution over time. Both quantities are measured in the discrete $L^2$ norm, as in \eqref{eq:rom-relative-errors} with the high-fidelity solution replaced by the exact field.
The error histories remain small and stable, which verifies the accuracy of the reference discretization for this smooth test. 

\begin{figure}[H]
    \centering
    \makebox[\textwidth][c]{%
    \subfloat[$t=2.5$]{
        \includegraphics[width=0.235\textwidth]{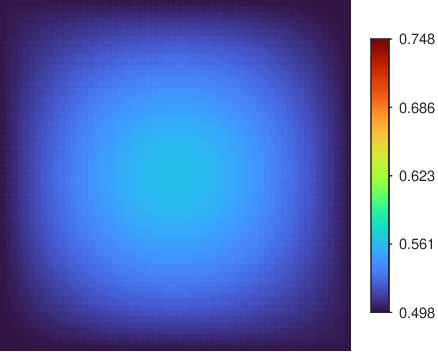}
    }\hspace{0.005\textwidth}
    \subfloat[$t=5.0$]{
        \includegraphics[width=0.235\textwidth]{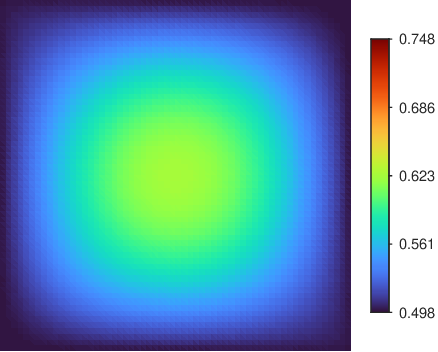}
    }\hspace{0.005\textwidth}
    \subfloat[$t=7.5$]{
        \includegraphics[width=0.235\textwidth]{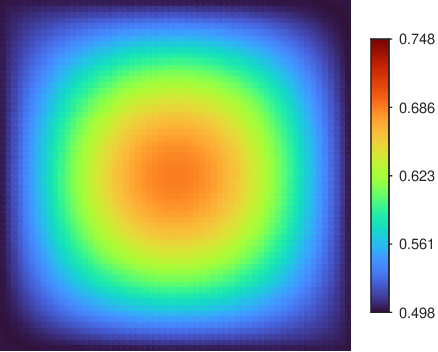}
    }\hspace{0.005\textwidth}
    \subfloat[$t=10.0$]{
        \includegraphics[width=0.235\textwidth]{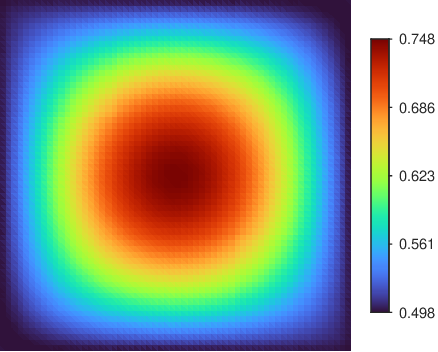}
    }}
    \vspace{0.05cm}
    \makebox[\textwidth][c]{%
    \subfloat[$t=2.5$]{
        \includegraphics[width=0.235\textwidth]{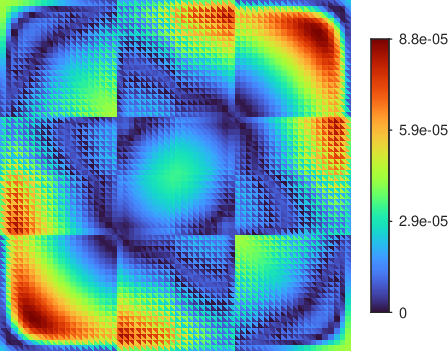}
    }\hspace{0.005\textwidth}
    \subfloat[$t=5.0$]{
        \includegraphics[width=0.235\textwidth]{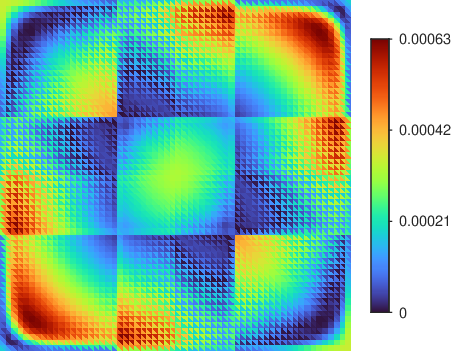}
    }\hspace{0.005\textwidth}
    \subfloat[$t=7.5$]{
        \includegraphics[width=0.235\textwidth]{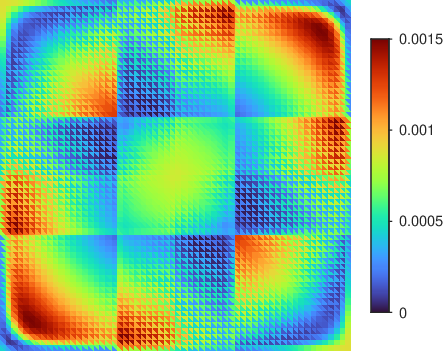}
    }\hspace{0.005\textwidth}
    \subfloat[$t=10.0$]{
        \includegraphics[width=0.235\textwidth]{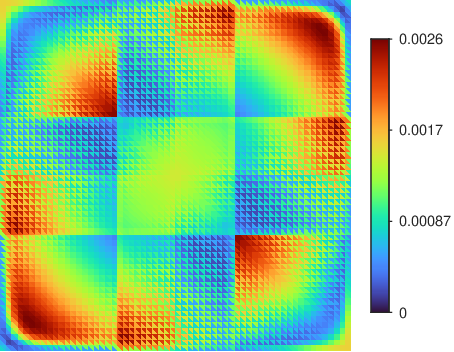}
    }}
    \caption{LSPG-FI saturation prediction. Top row: $S_g^{\mathrm{rom}}$. Bottom row: absolute ROM error $|S_g^{\mathrm{rom}}-S_g^{\mathrm{hf}}|$.}
    \label{fig:case1-lspgfi-saturation-time-error}
\end{figure}

We next compare the four DD-ROM variants introduced in Section~\ref{subsec:transport-rom}. 
The reduced bases are constructed from 200 high-fidelity time snapshots, after which the ROM is restarted from \(t=0\) and advanced over the full time interval. 
Figure~\ref{fig:case1-lspgfi-saturation-time-error} shows the LSPG-FI saturation prediction at the same representative times as in Figure~\ref{fig:case1-fom-saturation-time-error}. 
The bottom row reports the absolute difference with respect to the high-fidelity solution. 
The ROM accurately reproduces the smooth saturation evolution, and the error remains small throughout the domain.

\begin{figure}[H]
    \centering
    \subfloat[Pressure\label{fig:ex1-rom4-pressure}]{
        \includegraphics[width=0.48\textwidth]{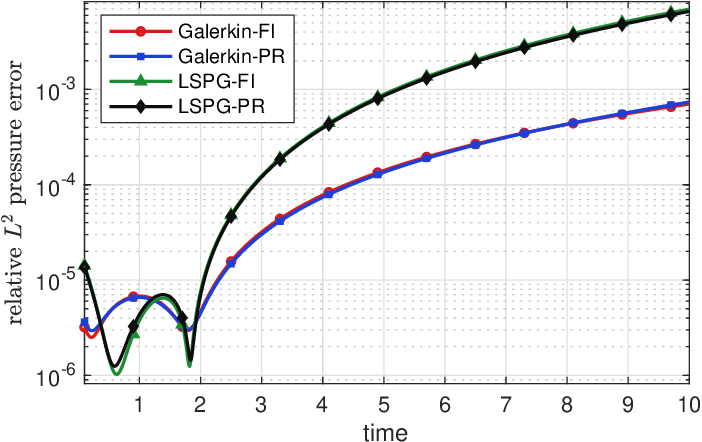}
    }
    \hfill
    \subfloat[Saturation\label{fig:ex1-rom4-saturation}]{
        \includegraphics[width=0.48\textwidth]{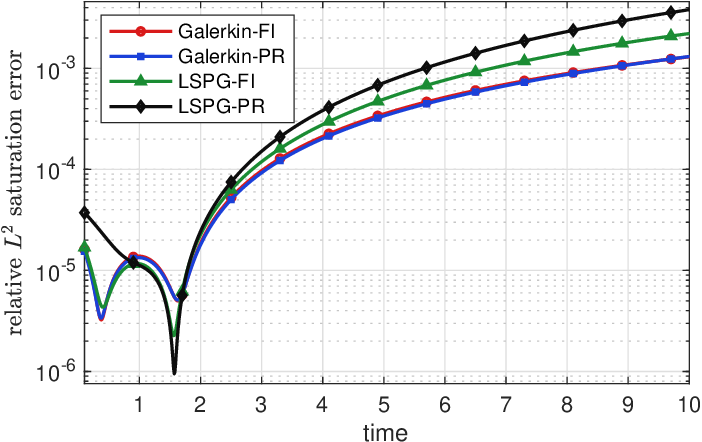}
    }
    \caption{ROM relative-error histories for the four DD-ROM variants: (a) relative $L^2$ pressure error and (b) relative $L^2$ saturation error.}
    \label{fig:ex1_ROM_errline_4methods}
\end{figure}

Figure~\ref{fig:ex1_ROM_errline_4methods} compares the temporal pressure and saturation error histories for the four DD-ROM variants. 
All methods accurately reconstruct the high-fidelity solution in this smooth setting. 
The full-interface and port-reduced variants have comparable errors, indicating that the interface-pressure manifold is well represented by the port basis for this training configuration. 

\subsection{Test case 2: parameterized permeability}
\label{sec:num_case2}

The test evaluates out-of-sample prediction for a parameterized heterogeneous permeability field. This setting is designed to test whether the local HAPOD bases can represent changes in the pressure response, flux pattern, and saturation plume induced by variations in the rock properties. The domain is \(\Omega=[0,1]\times[0,1]\) with  \(3\times3\) nonoverlapping partition, each subdomain contains \(800\) triangular cells. Dirichlet data are imposed on the outer boundary with \(p_D=0\) and \(S_D=0\). The final time is \(T=1\), the time step is \(\Delta t=10^{-2}\), and the well support is
\[
  \omega_w=[0.44,0.56]\times[0.44,0.56],
\]
located in the central subdomain. We prescribe a constant gas source
$
Q_{\mathrm{CO}_2}(t)=|\omega_w|.
$
The capillary-pressure and relative-permeability laws are the same as test case 1 in Section \ref{sec:num_case1}. 
The remaining dimensionless parameters are
$
\phi=0.3,
\mu_a=1,
\mu_g=0.8,
\rho_a=1.1,
\rho_g=1,
X_a^{\mathrm{CO}_2}=0.1,
X_g^{\mathrm{CO}_2}=0.8.
$
The permeability is isotropic but spatially heterogeneous. We parameterize the logarithm of the scalar permeability by
\begin{equation}
\log_{10} K(\mathbf{x};\theta_K)
=
Z_{\rm perm}(\mathbf{x})
+0.28\,\theta_K\,\psi_1(\mathbf{x})
+0.16\,\theta_K^2\,\psi_2(\mathbf{x})
+0.10\,\theta_K,
\label{eq:case2-perm-family}
\end{equation}
where \(Z_{\rm perm}\) is a fixed smooth surrogate for a Gaussian log-permeability realization with standard deviation about \(0.3\), and \(\psi_1\) and \(\psi_2\) are normalized smooth spatial perturbations. 
We randomly generate 10 training points and 5 test points for $ \theta_K\in [-1,1]$. 
Local reduced bases are first constructed from the high-fidelity snapshots at the training parameters and are then merged by a second HAPOD level to obtain one parameter-independent basis per subdomain and variable.

\begin{figure}[H]
    \centering
    \subfloat[\(\theta_K=-1\)]{
        \includegraphics[width=0.225\textwidth]{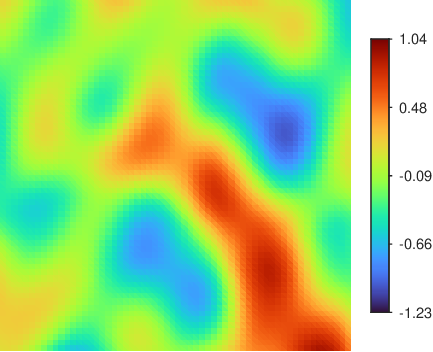}
    }\hspace{0.005\textwidth}
    \subfloat[\(\theta_K=-0.75\)]{
        \includegraphics[width=0.225\textwidth]{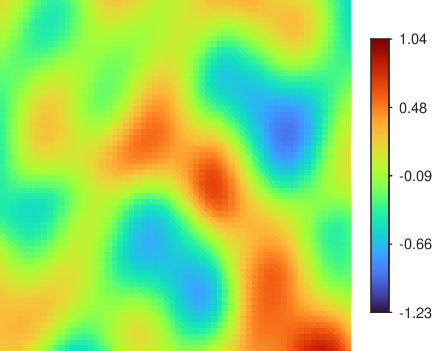}
    }\hspace{0.005\textwidth}
    \subfloat[\(\theta_K=0\)]{
        \includegraphics[width=0.225\textwidth]{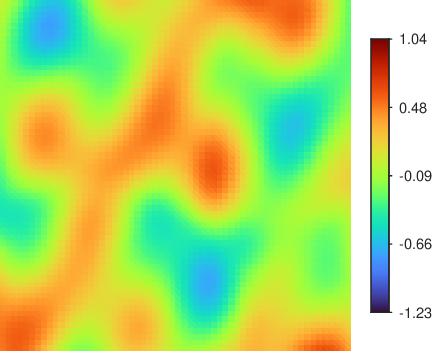}
    }\hspace{0.005\textwidth}
    \subfloat[\(\theta_K=1\)]{
        \includegraphics[width=0.225\textwidth]{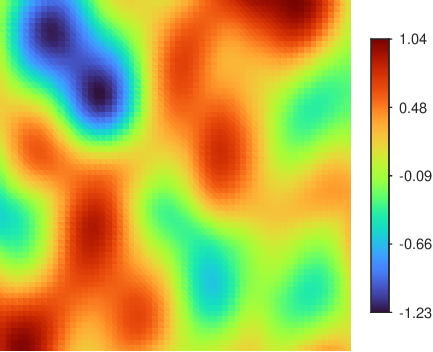}
    }
    \caption{Permeability-parametric case: representative fields of \(\log_{10}K(\mathbf{x};\theta_K)\).}
    \label{fig:case2-perm-logk-fields}
\end{figure}

\begin{figure}[H]
    \centering
    \subfloat[\(\theta_K=-1\)]{
        \includegraphics[width=0.225\textwidth]{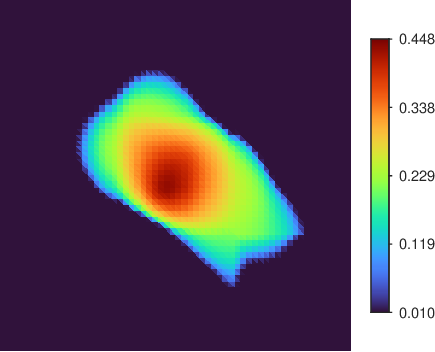}
    }\hspace{0.005\textwidth}
    \subfloat[\(\theta_K=-0.75\)]{
        \includegraphics[width=0.225\textwidth]{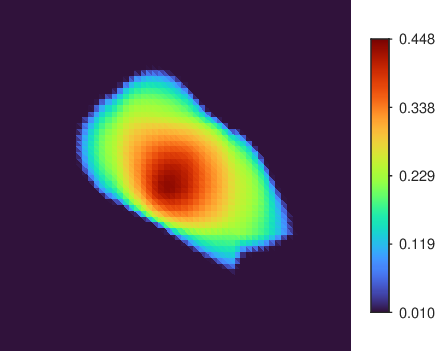}
    }\hspace{0.005\textwidth}
    \subfloat[\(\theta_K=0\)]{
        \includegraphics[width=0.225\textwidth]{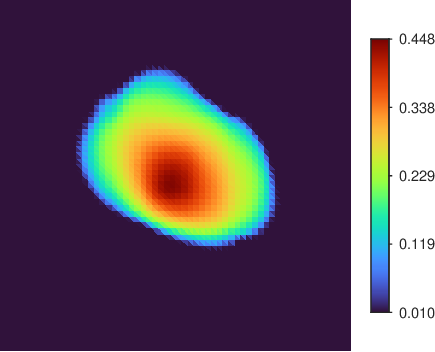}
    }\hspace{0.005\textwidth}
    \subfloat[\(\theta_K=1\)]{
        \includegraphics[width=0.225\textwidth]{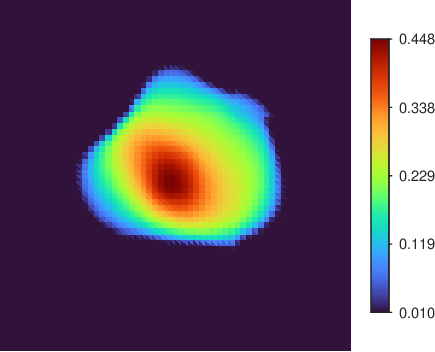}
    }
    \caption{Permeability-parametric case: high-fidelity gas saturation \(S_g\) at \(T=1\) for the permeability fields in Figure~\ref{fig:case2-perm-logk-fields}.}
    \label{fig:case2-perm-fom-saturation-fields}
\end{figure}

Figure~\ref{fig:case2-perm-logk-fields} shows representative log-permeability fields.  Figure~\ref{fig:case2-perm-fom-saturation-fields} gives the corresponding high-fidelity gas saturation at \(T=1\). 
The changes in permeability lead to visible differences in the plume shape.

\begin{figure}[H]
    \centering
    \subfloat[Pressure\label{fig:case2-perm-rom1-pressure-time}]{
        \includegraphics[width=0.48\textwidth]{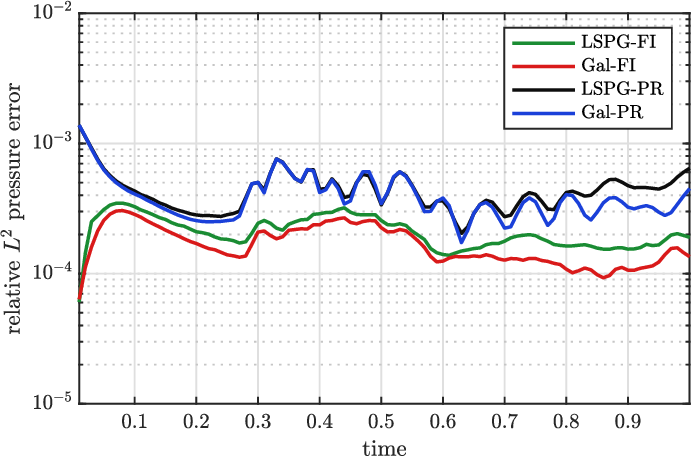}
    }
    \hfill
    \subfloat[Saturation\label{fig:case2-perm-rom1-saturation-time}]{
        \includegraphics[width=0.48\textwidth]{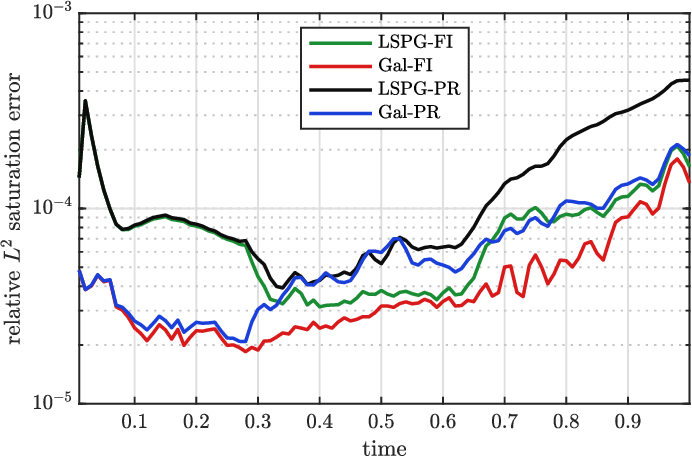}
    }
    \caption{Permeability-parametric case: pressure and saturation relative-error histories at the representative test point \(\theta_K=-0.5603\).}
    \label{fig:case2-perm-rom1-error-histories}
\end{figure}

We next compare the four ROM variants at the representative held-out point
\(\theta_K=-0.5603\). 
All errors are relative \(L^2\) errors with respect to the corresponding DD-FOM trajectory, computed at the 100 saved time levels. Figure~\ref{fig:case2-perm-rom1-error-histories} reports the temporal error histories. 

\subsection{Test case 3: fixed and moving wells}
\label{sec:num_case3}

This test case considers a more challenging CO$_2$ storage setting with a heterogeneous permeability field, stronger nonlinear transport effects, and moving well sources.

\subsubsection{Fixed-well scenario}
\label{ex3-1}

\begin{figure}[htbp]
    \centering
    \includegraphics[width=8cm]{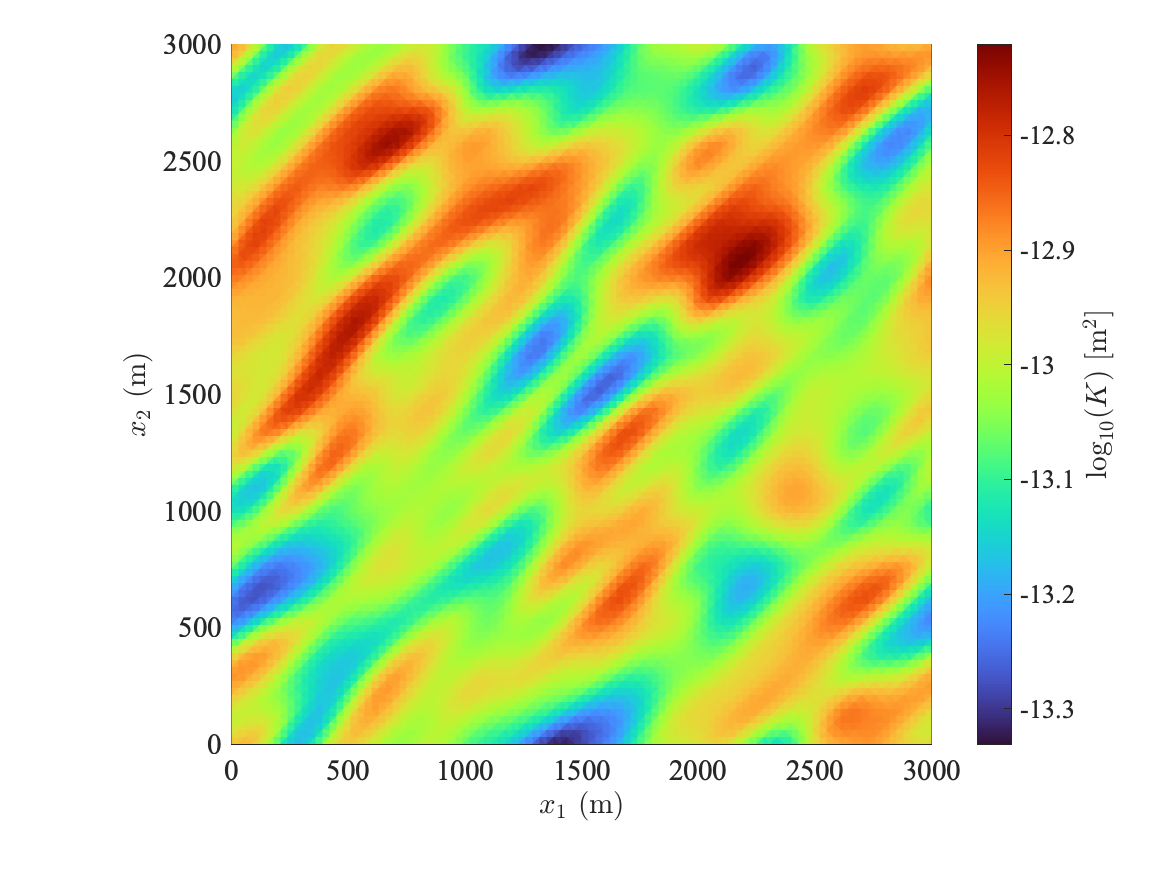}
    \caption{Test case 3: the heterogeneous anisotropic permeability field.}
    \label{fig:logK}
\end{figure}

We first consider a fixed injection well. The computational domain is
\(\Omega=[0,3]\times[0,3]~\mathrm{km}^2\), partitioned into \(5\times5\) nonoverlapping subdomains. 
Each subdomain is discretized into 800 triangular elements. CO$_2$ is injected through the rectangular well region
$
[1.32,1.68]\times[1.32,1.68]~\mathrm{km}^2,
$
located in subdomain $13$. 
The permeability is assumed to follow the log-normal distribution, and the correlation lengths are fixed to be $421.875$ m and $140.625$ m along the two anisotropy directions respectively. One realization of such permeability field is depicted in Figure \ref{fig:logK}.

\begin{figure}[H]
    \centering
    \includegraphics[width=8cm]{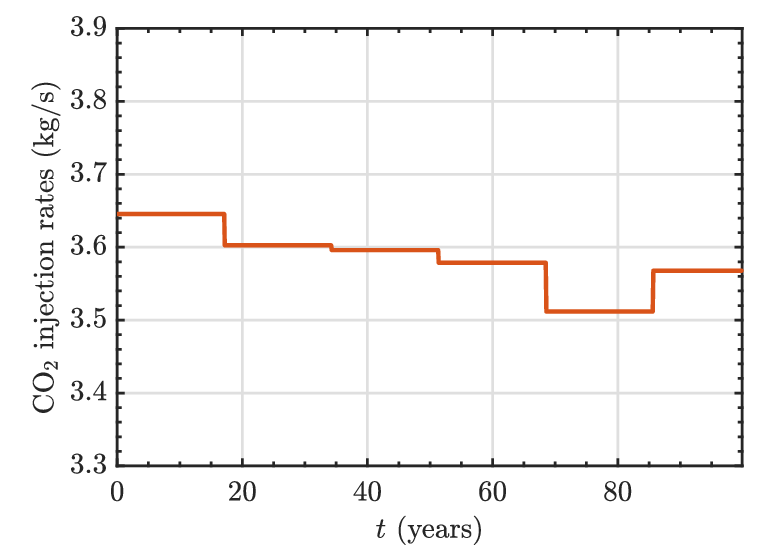}
    \caption{Test case 3: injection rate function.}
    \label{fig:inj rate}
\end{figure}

The remaining physical parameters are temperature $37^\circ\mathrm{C}$, $\phi=0.3$, $\mu_a=1\times10^{-3}\ \mathrm{Pa\,s}$, $\mu_g=5\times10^{-5}\ \mathrm{Pa\,s}$, $\rho_a=1000\ \mathrm{kg/m^3}$, $\rho_g=684\ \mathrm{kg/m^3}$, $X_a^{\mathrm{CO}_2}=0$, and $X_g^{\mathrm{CO}_2}=1$.
The density of CO$_2$ is obtained from the NIST Chemistry WebBook at 37$^\circ$C and 10 MPa, corresponding to supercritical conditions. Since an immiscible incompressible two-phase flow model is considered and dissolution effects are neglected, the CO$_2$ phase and aqueous phase are assumed to be pure phases, resulting in $X_a^{\mathrm{CO}_2}=0$ and $X_g^{\mathrm{CO}_2}=1$. Dirichlet boundary conditions are imposed on all four outer boundaries, with the prescribed pressure $p_D = 10\ \text{MPa}$ and gas saturation $S_D = 0$. 
The source is imposed as a localized CO$_2$ mass injection in the well region. The operational injection profile is the piecewise-constant curve shown in Figure~\ref{fig:inj rate}, and the total source is assigned consistently by setting $Q_t=Q_{\mathrm{CO}_2}$. The simulation is run for $100$ years with a uniform time step of $1$ year.

\begin{figure}[H]
    \centering
    \subfloat[$t=30$ years\label{fig:ex3 HF PS a}]{
        \includegraphics[width=0.31\textwidth]{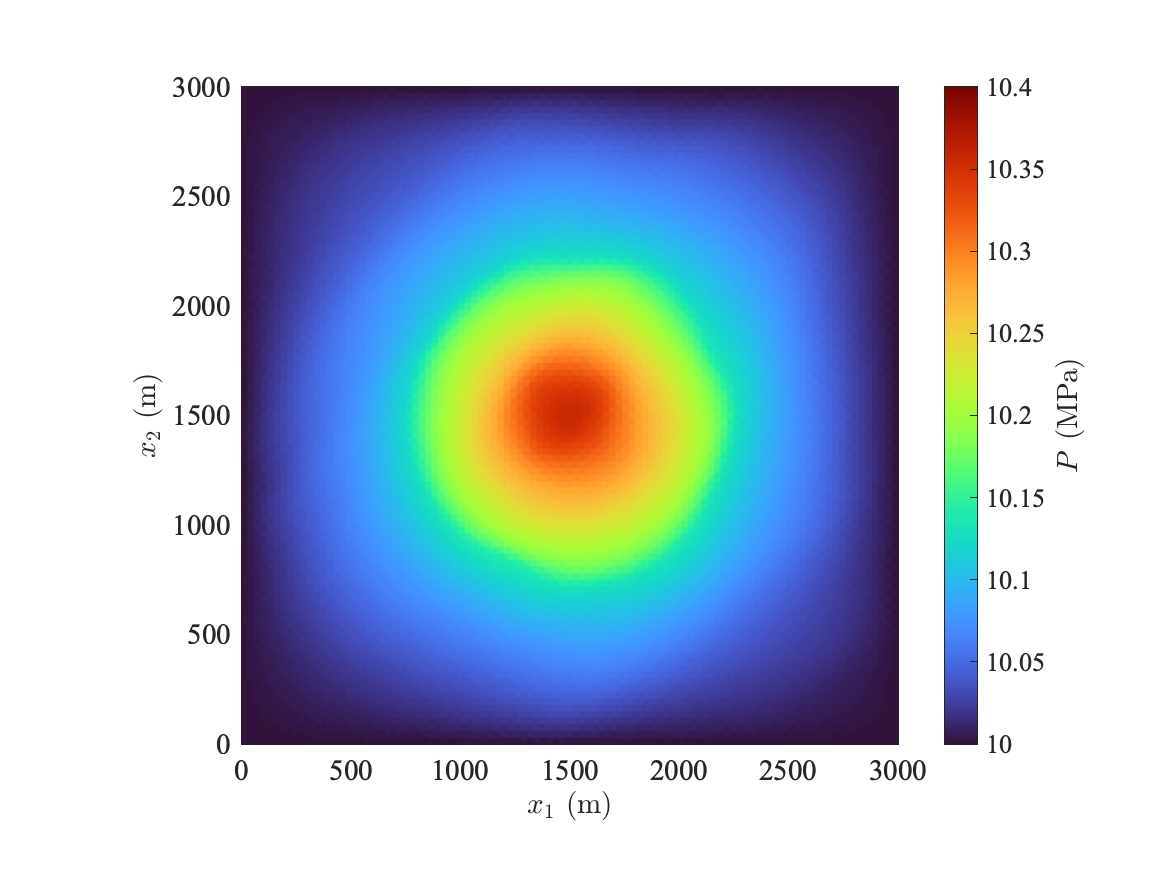}
    }
    \hfill
    \subfloat[$t=60$ years \label{fig:ex3 HF PS b}]{
        \includegraphics[width=0.31\textwidth]{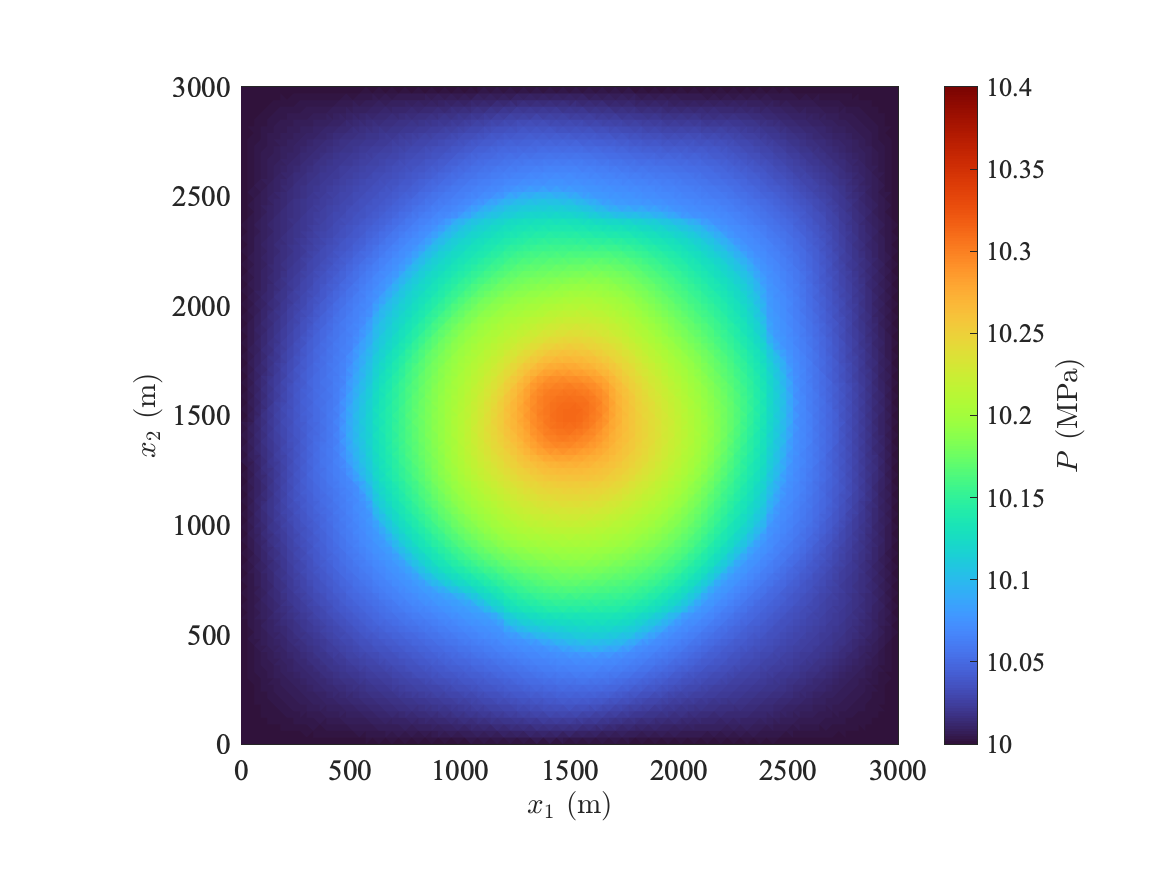}
    }
    \hfill
    \subfloat[$t=100$ years\label{fig:ex3 HF PS c}]{
        \includegraphics[width=0.31\textwidth]{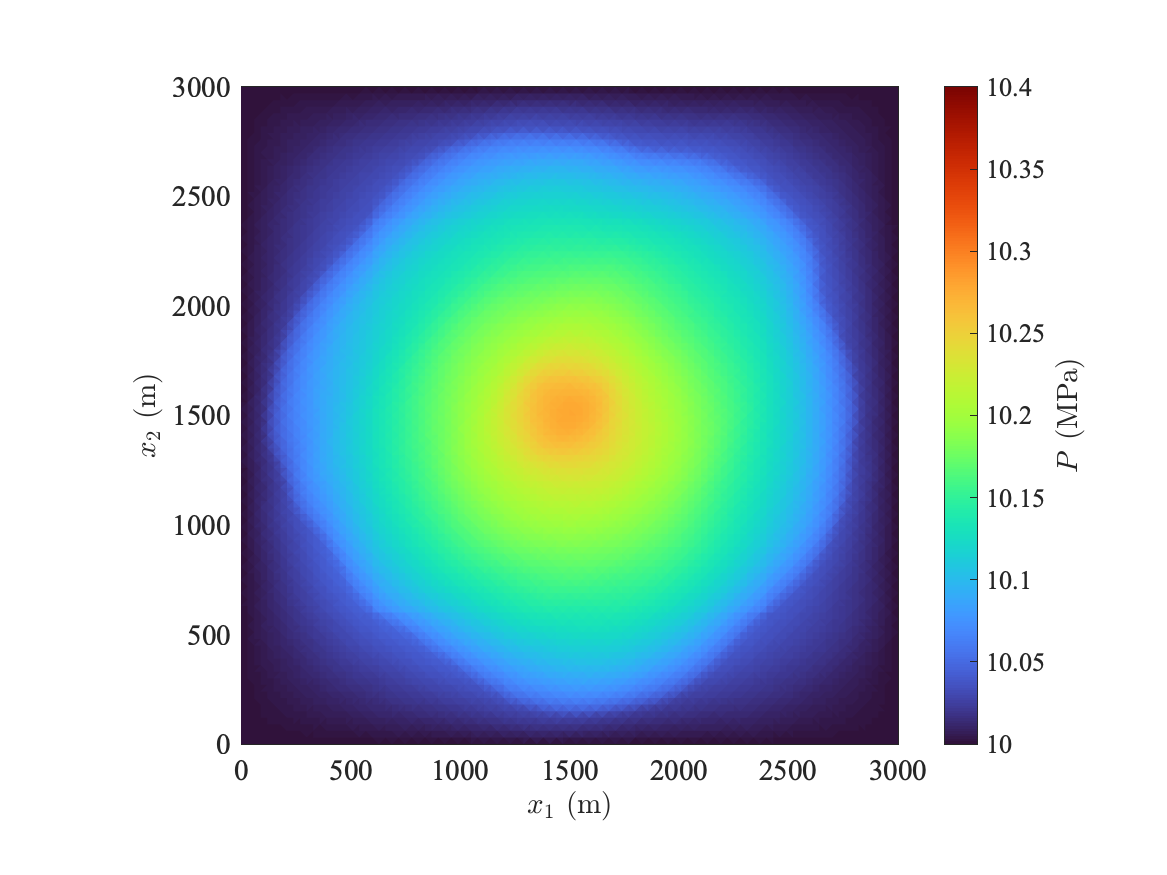}
    }
    \vspace{0.1cm} 
    \subfloat[$t=30$ years\label{fig:ex3 HF PS d}]{
        \includegraphics[width=0.31\textwidth]{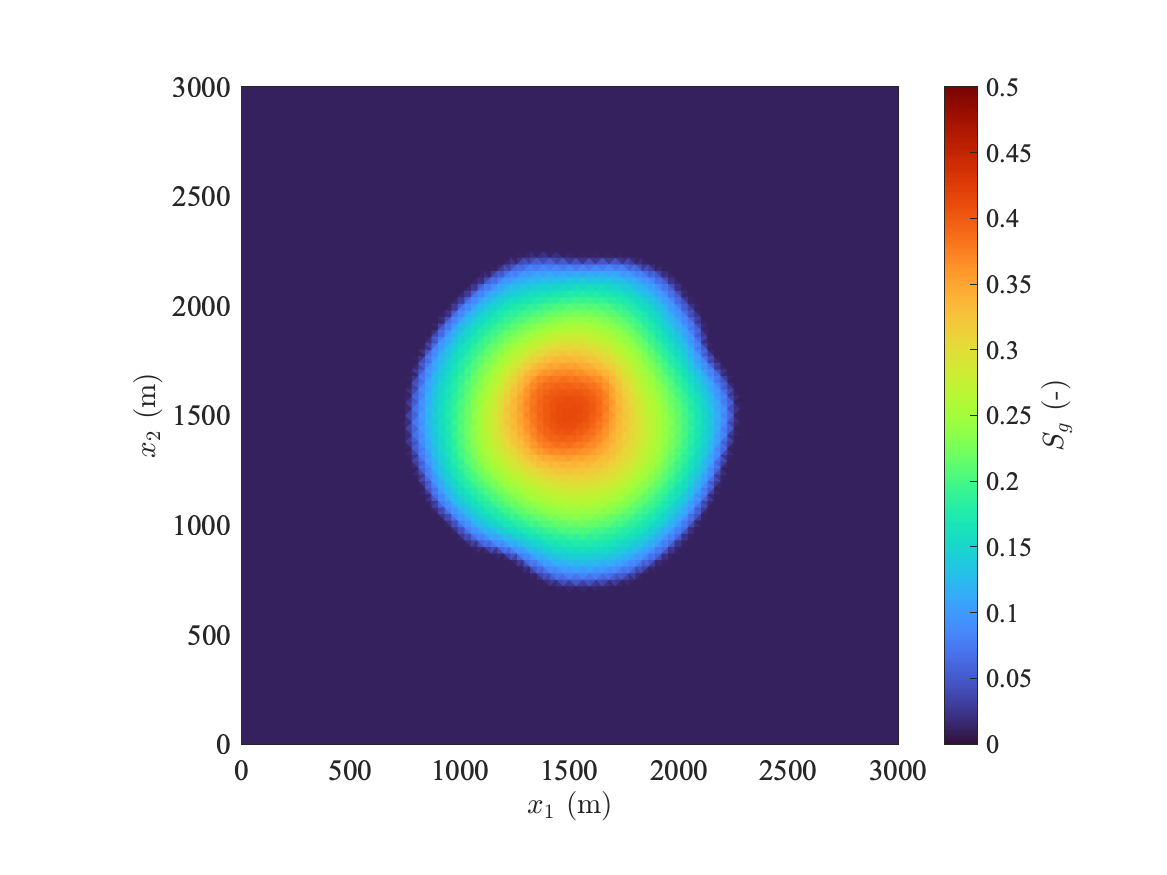}
    }
    \hfill
    \subfloat[$t=60$ years\label{fig:ex3 HF PS e}]{
        \includegraphics[width=0.31\textwidth]{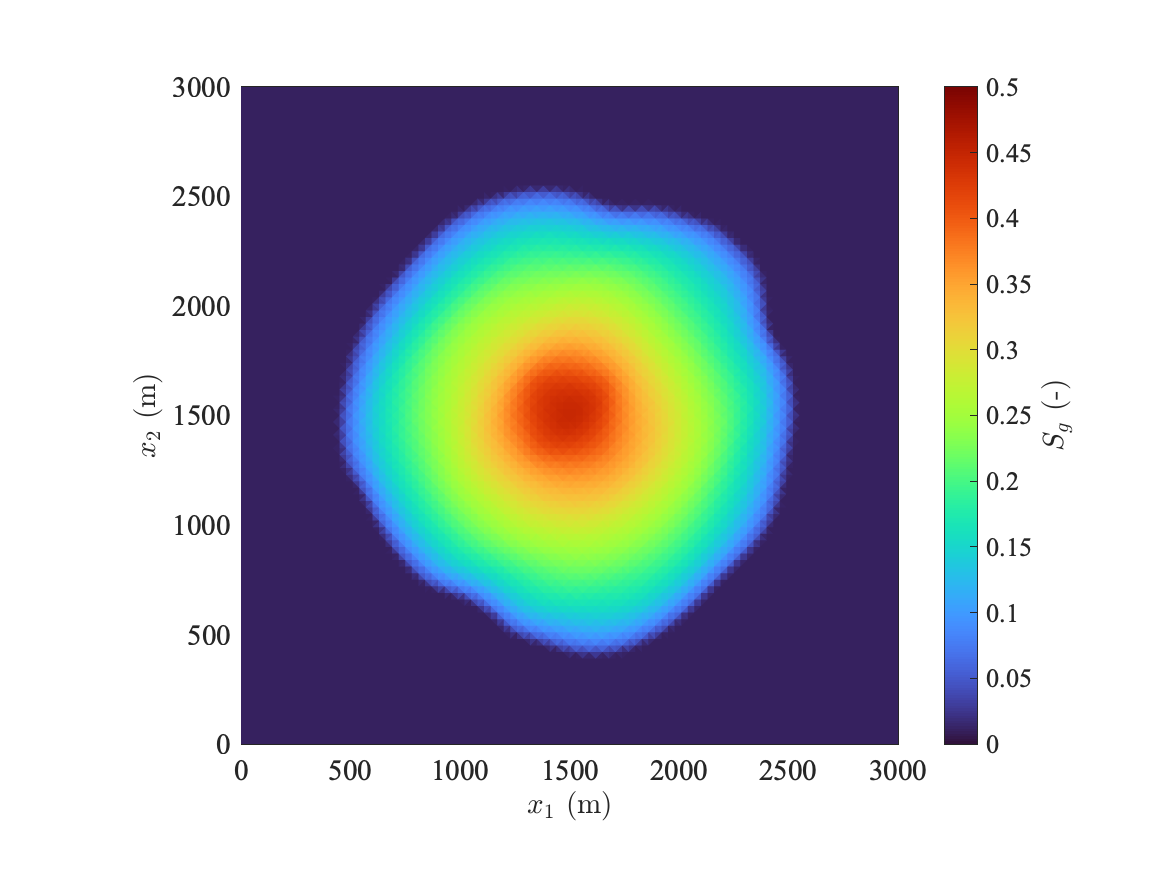}
    }
     \hfill
    \subfloat[$t=100$ years\label{fig:ex3 HF PS f}]{
        \includegraphics[width=0.31\textwidth]{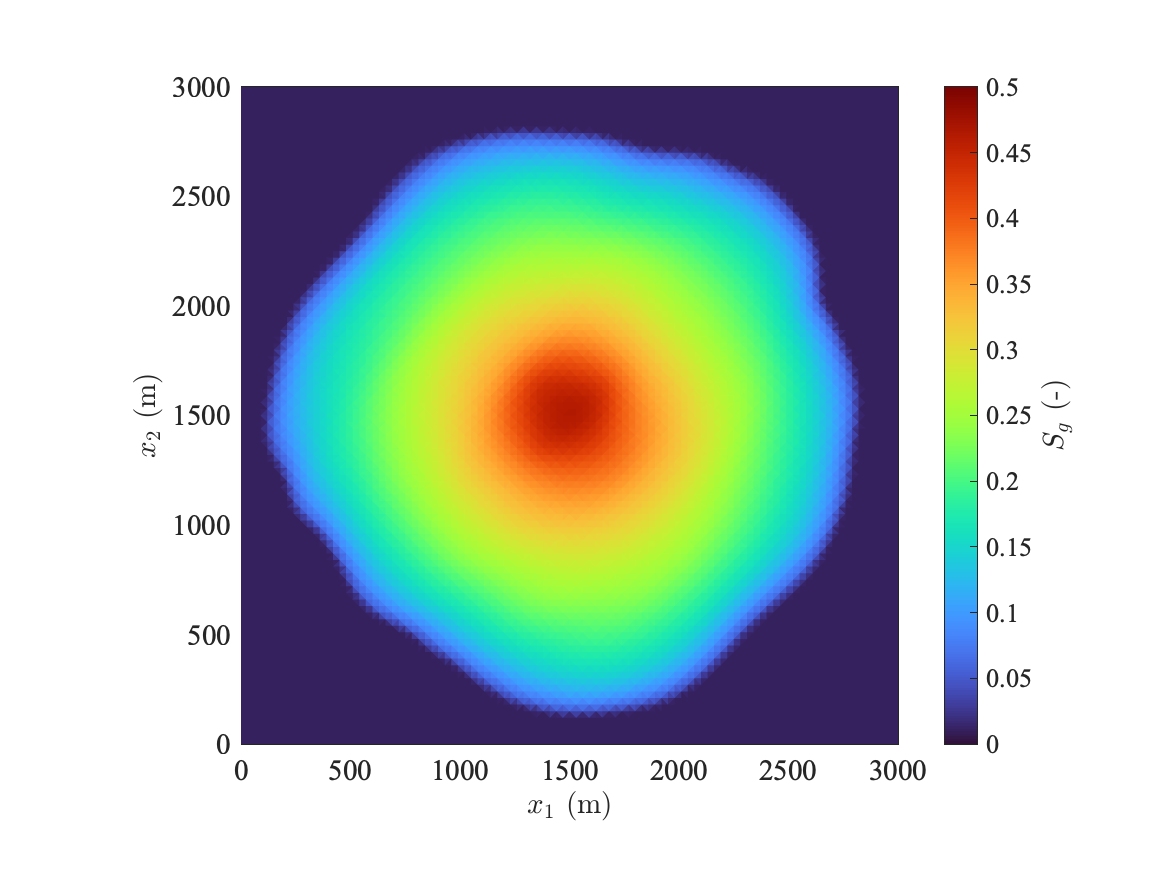}
    }
    \caption{HF results of pressure and saturation. (a)-(c) Pressure; (d)-(f) Saturation.}
    \label{fig:ex3 HF PS}
\end{figure}

Figure \ref{fig:ex3 HF PS} illustrates the evolution of the high-fidelity (HF) pressure and saturation fields at $t=30$ years, $t=60$ years and $t=100$ years. Driven by the central injection, the fluid fields exhibit strong radial displacement characteristics. 

\begin{figure}[H]
    \centering
    \subfloat[\label{fig:ex3q1 errline 2m b}]{
        \includegraphics[width=0.45\textwidth]{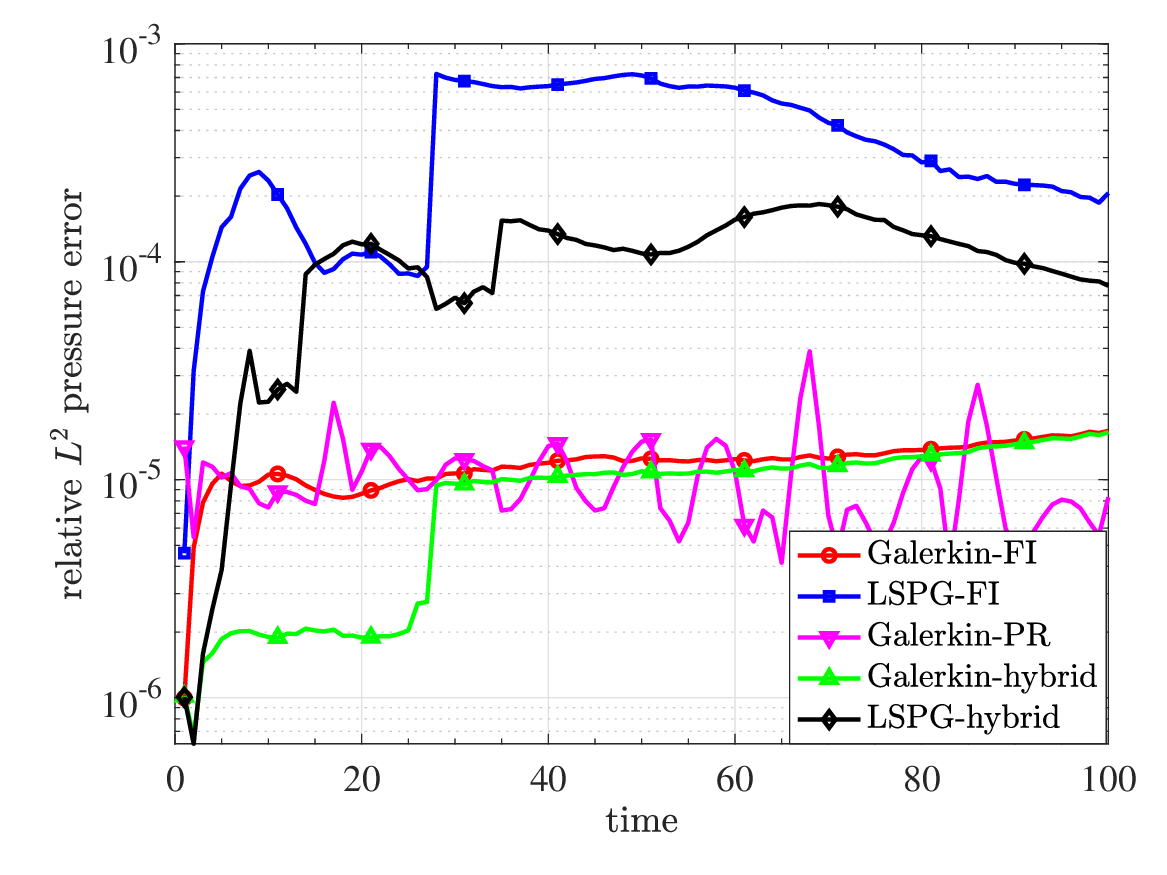}
    }
    \hfill
    \subfloat[\label{fig:ex3q1 errline 2m c}]{
        \includegraphics[width=0.45\textwidth]{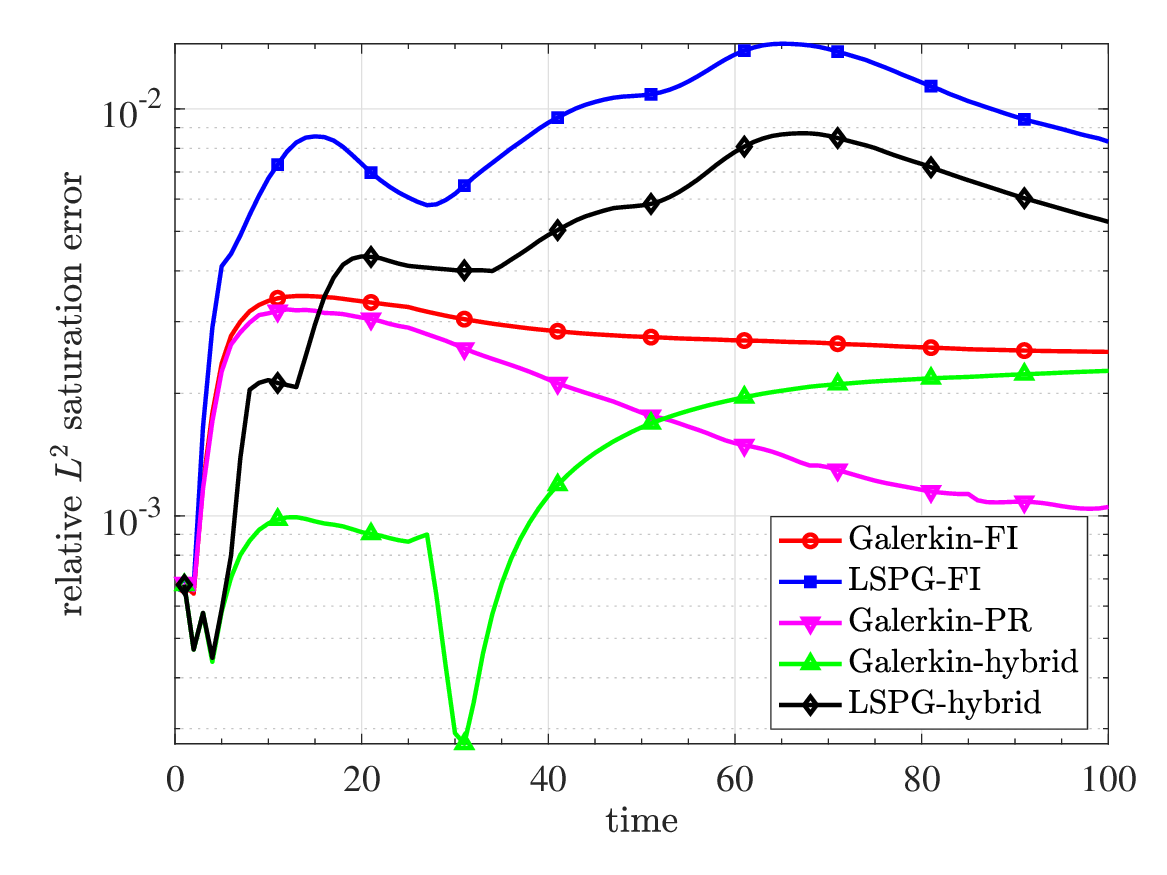}
    }
    \caption{Relative errors of LSPG-FI and LSPG-hybrid in a purely reconstructive setting.}
    \label{fig:ex3q1 errline 2m}
\end{figure}

We first test the ROMs in a reconstruction setting. The reduced spaces are trained using high-fidelity snapshots over the full 100-year interval and are then used to simulate the same time interval. Figure~\ref{fig:ex3q1 errline 2m} shows the relative error of Galerkin-FI, LSPG-FI, Galerkin-hybrid and LSPG-hybrid. All four ROM variants remain accurate over the full simulation window: the saturation errors are of order \(10^{-2}\), while the pressure errors remain below \(10^{-3}\).

\begin{figure}[H]
    \centering
    \subfloat[\label{fig:ex3q1 pre errline 2m b}]{
        \includegraphics[width=0.45\textwidth]{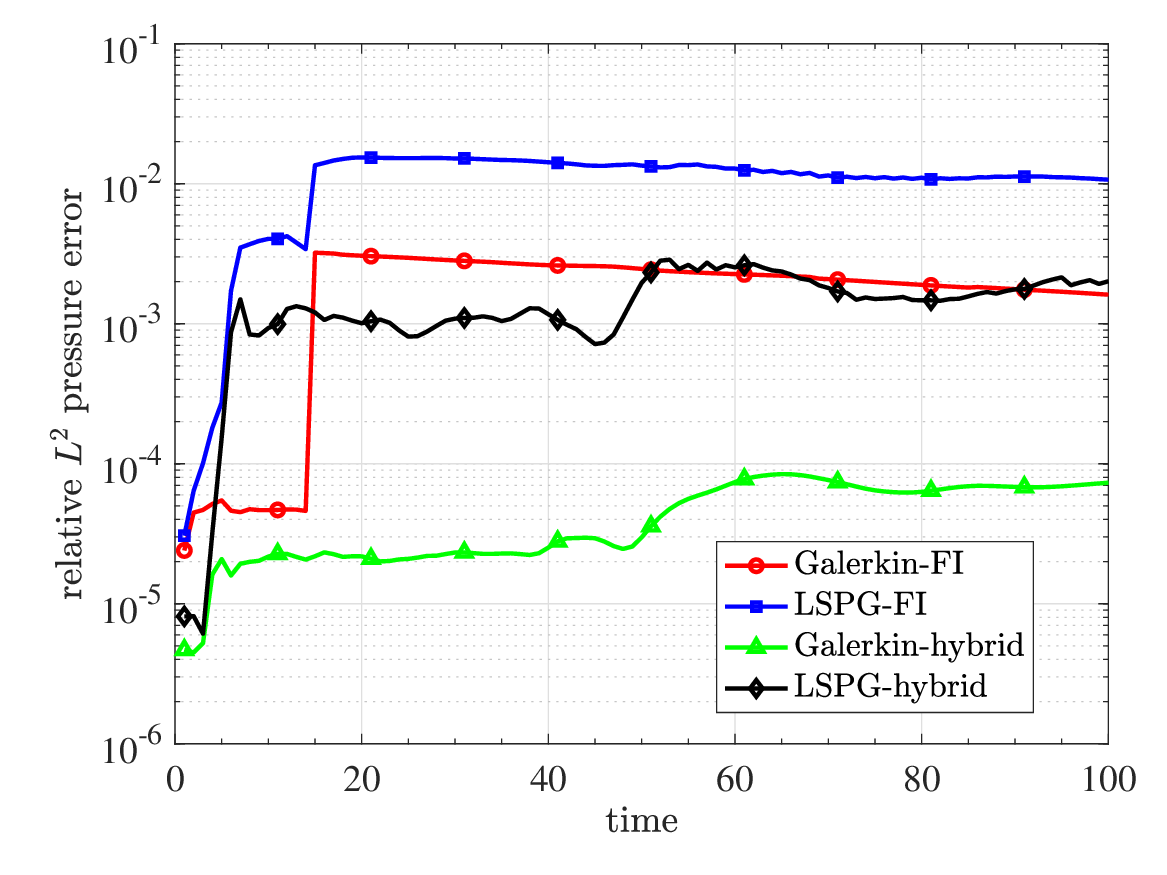}
    }
    \hfill
    \subfloat[\label{fig:ex3q1 pre errline 2m c}]{
        \includegraphics[width=0.45\textwidth]{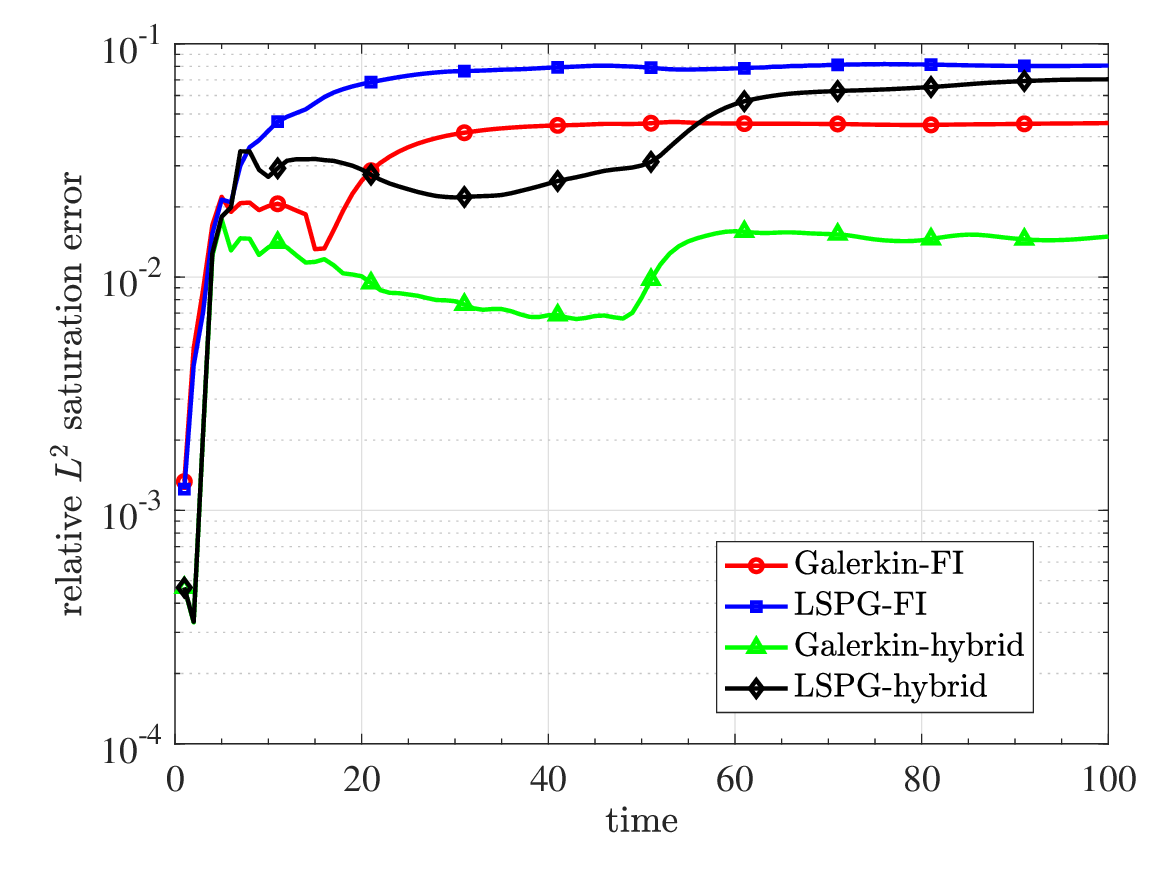}
    }
    \caption{Relative errors of LSPG-FI and LSPG-hybrid with time extrapolation.}
    \label{fig:ex3q1 pre errline 2m}
\end{figure}

We next evaluate time extrapolation. In this test, the reduced spaces are built only from the first 50 time steps and are then used to advance the ROMs to the final time. This tests whether the HAPOD bases remain robust as the plume continues to grow beyond the training window. 
Figure~\ref{fig:ex3q1 pre errline 2m} illustrates the relative errors of the Galerkin-FI, LSPG-FI, Galerkin-hybrid and LSPG-hybrid methods to evaluate their predictive capabilities. In this extrapolative regime, Galerkin-FI gives smaller errors than LSPG-FI, and the hybrid variants further improve robustness by resolving the well and near-well subdomains at high fidelity. 

\subsubsection{Moving well scenario}

We now test whether the reduced spaces trained from the original well configuration can be transferred to new well locations. All physical parameters, boundary data, permeability field, and injection profile are kept the same as in Section~\ref{ex3-1}; only the well location is changed. This setting is challenging because the dominant near-well pressure response and saturation front are displaced to subdomains that were not active as well subdomains during training.

For the moving-well tests, HAPOD is first applied to the velocity, pressure, and saturation snapshots from the original fixed-well problem. The resulting modes are organized into well-region and well-free groups and are used to build ROMs for the new well configurations. We focus on full-interface and hybrid variants. Port-reduced variants are not reported here, due to poor performance.

\begin{figure}[H]
    \centering
    \subfloat[$t=10$ years\label{fig:ex3-2 HF PS a}]{
        \includegraphics[width=0.31\textwidth]{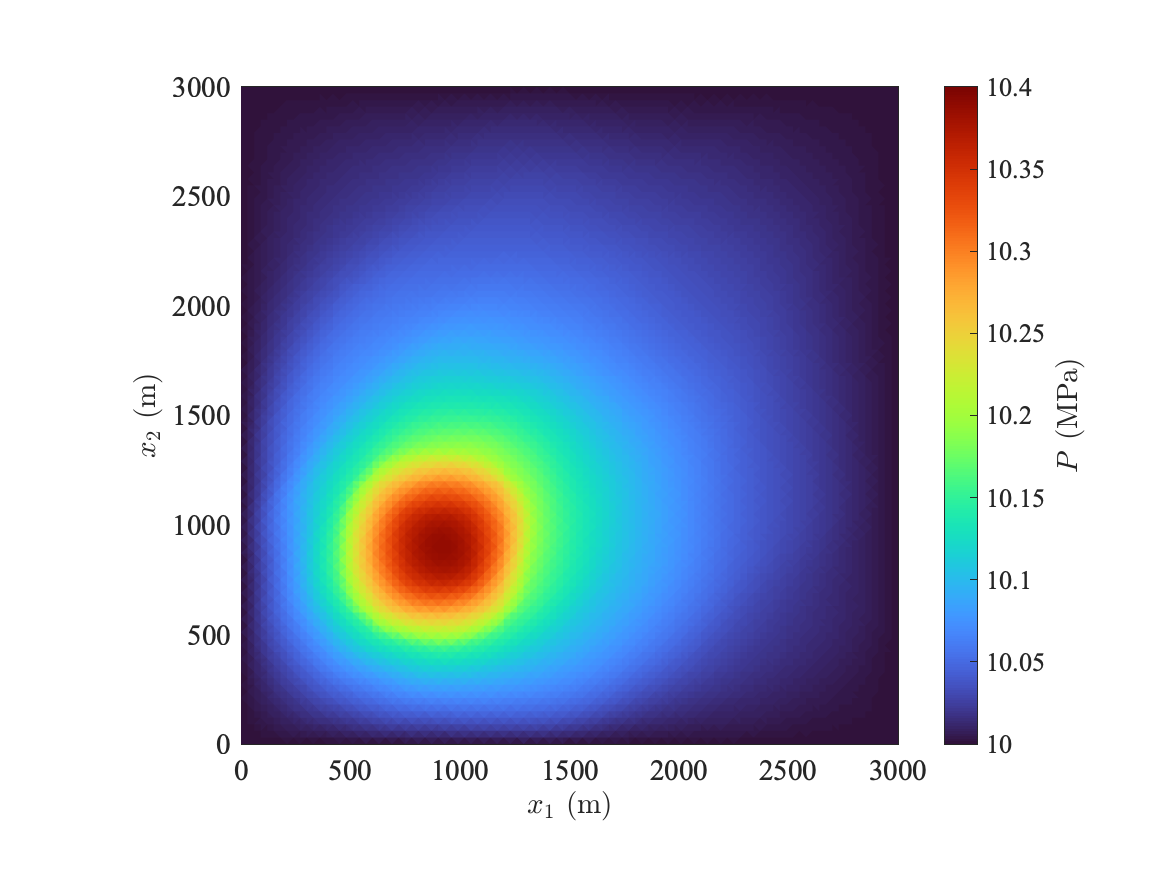}
    }
    \hfill
    \subfloat[$t=20$ years\label{fig:ex3-2 HF PS b}]{
        \includegraphics[width=0.31\textwidth]{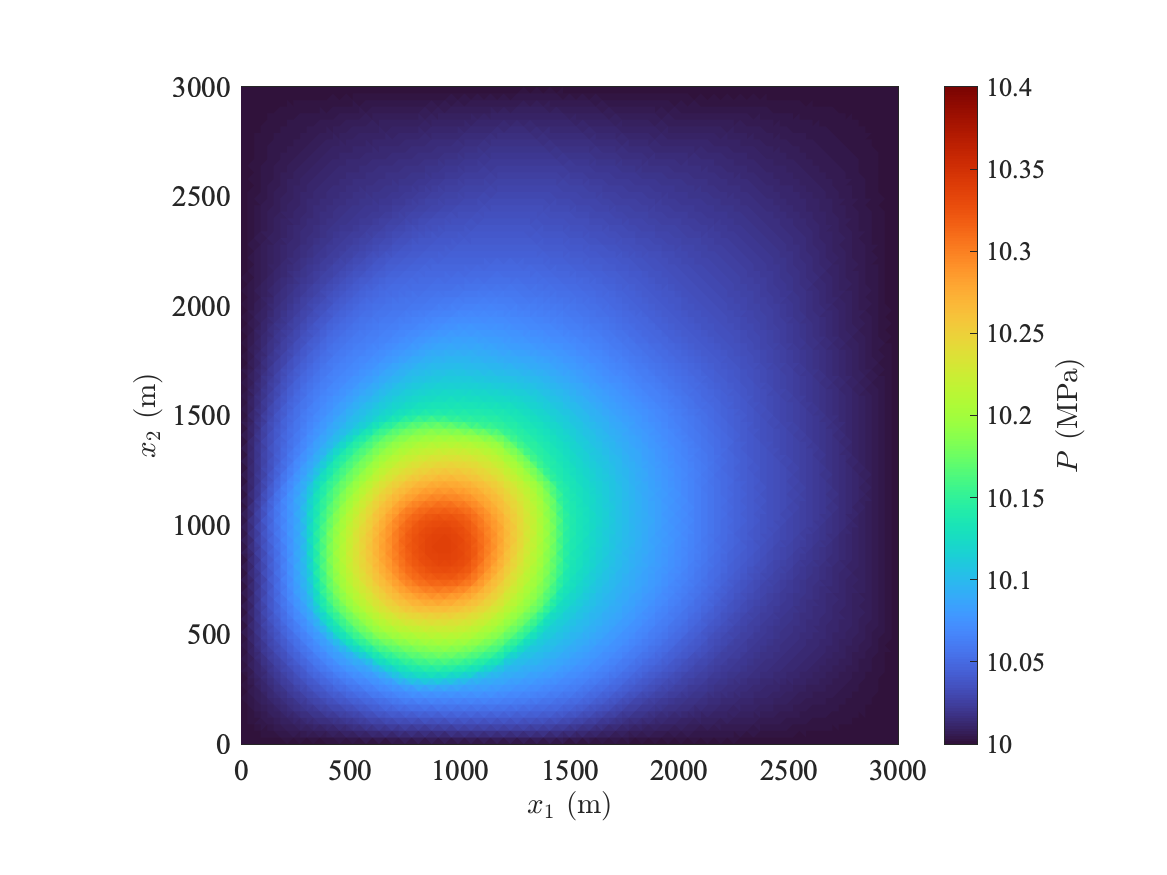}
    }
    \hfill
    \subfloat[$t=30$ years\label{fig:ex3-2 HF PS c}]{
        \includegraphics[width=0.31\textwidth]{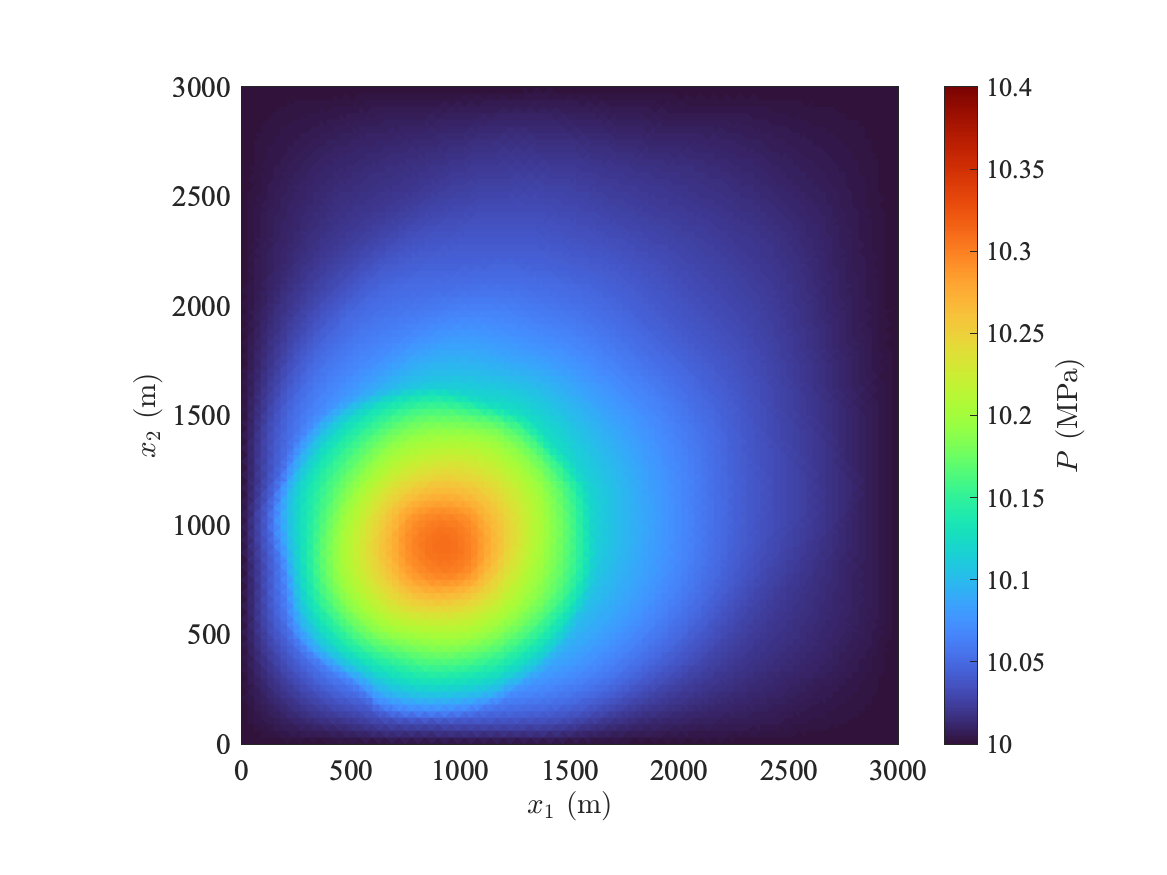}
    }
    \vspace{0.1cm} 
    \subfloat[$t=10$ years \label{fig:ex3-2 HF PS d}]{
        \includegraphics[width=0.31\textwidth]{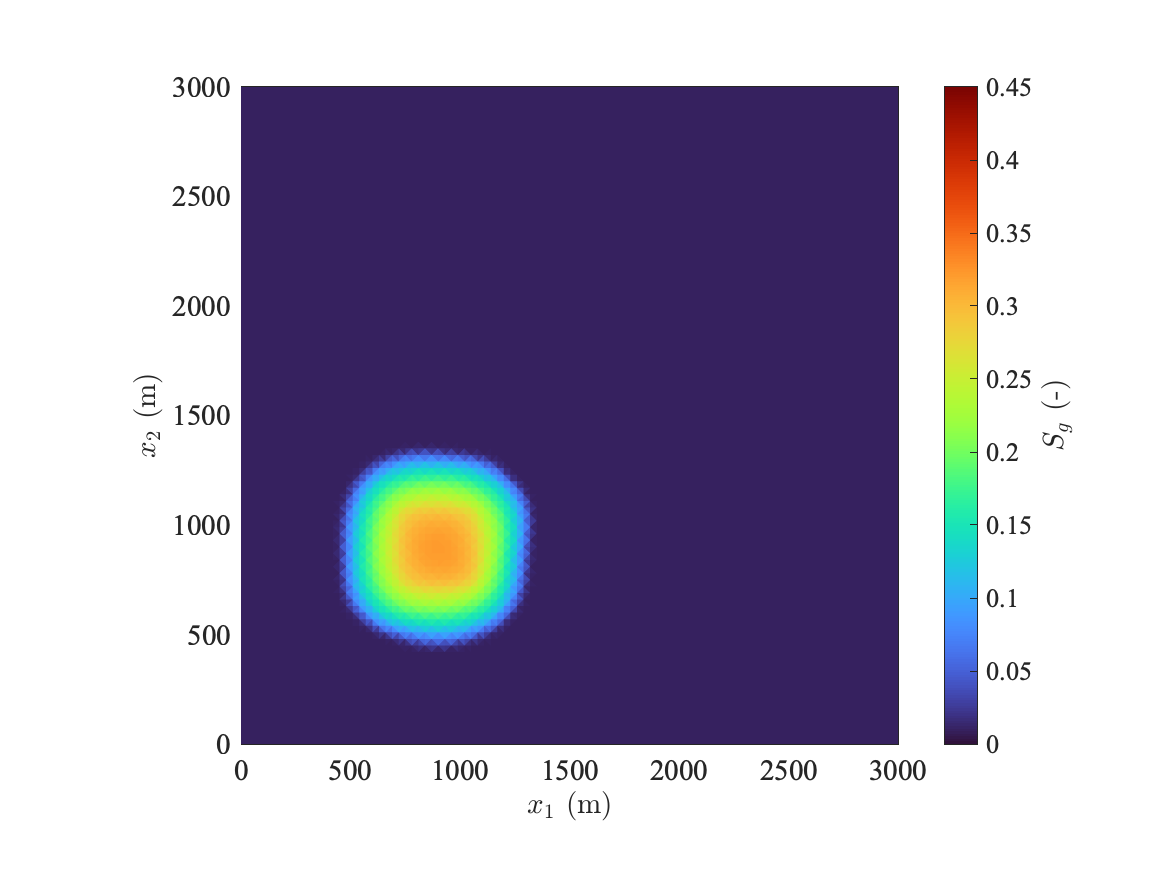}
    }
    \hfill
    \subfloat[$t=20$ years \label{fig:ex3-2 HF PS e}]{
        \includegraphics[width=0.31\textwidth]{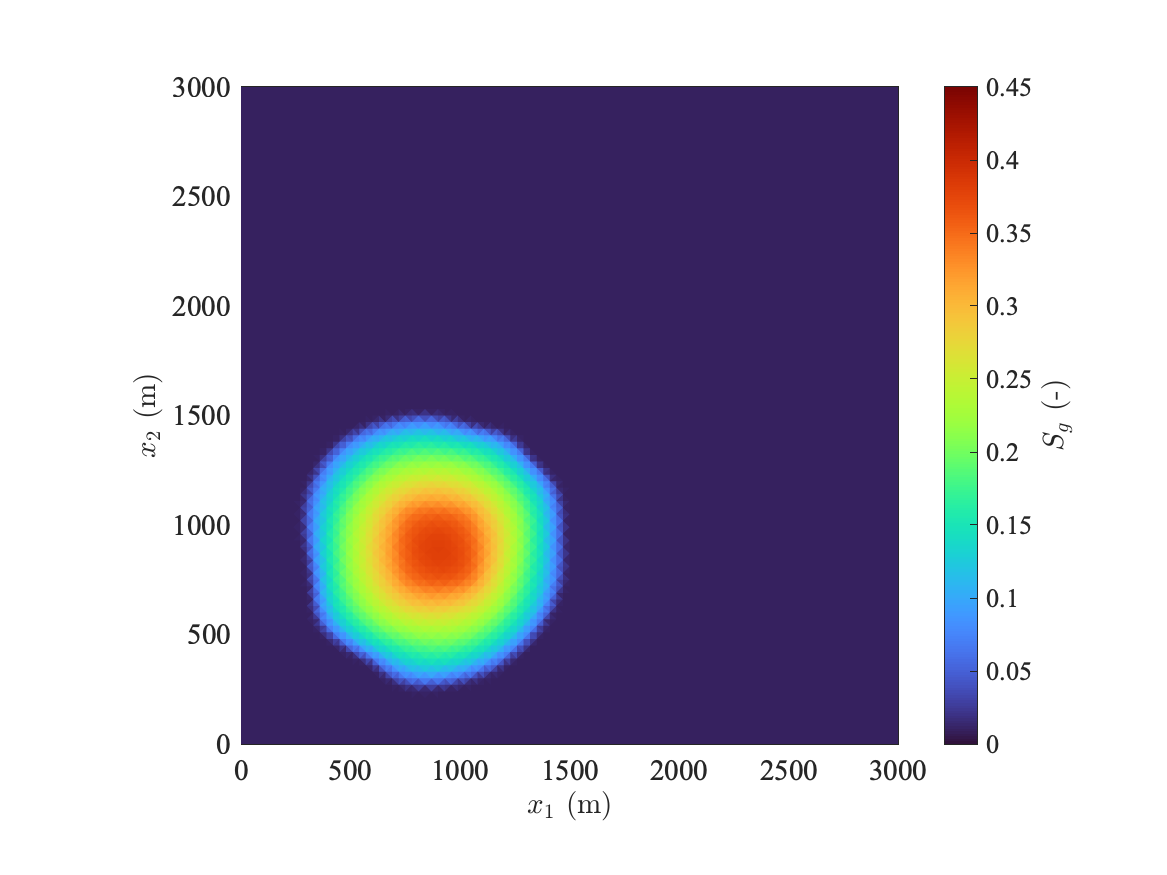}
    }
     \hfill
    \subfloat[$t=30$ years\label{fig:ex3-2 HF PS f}]{
        \includegraphics[width=0.31\textwidth]{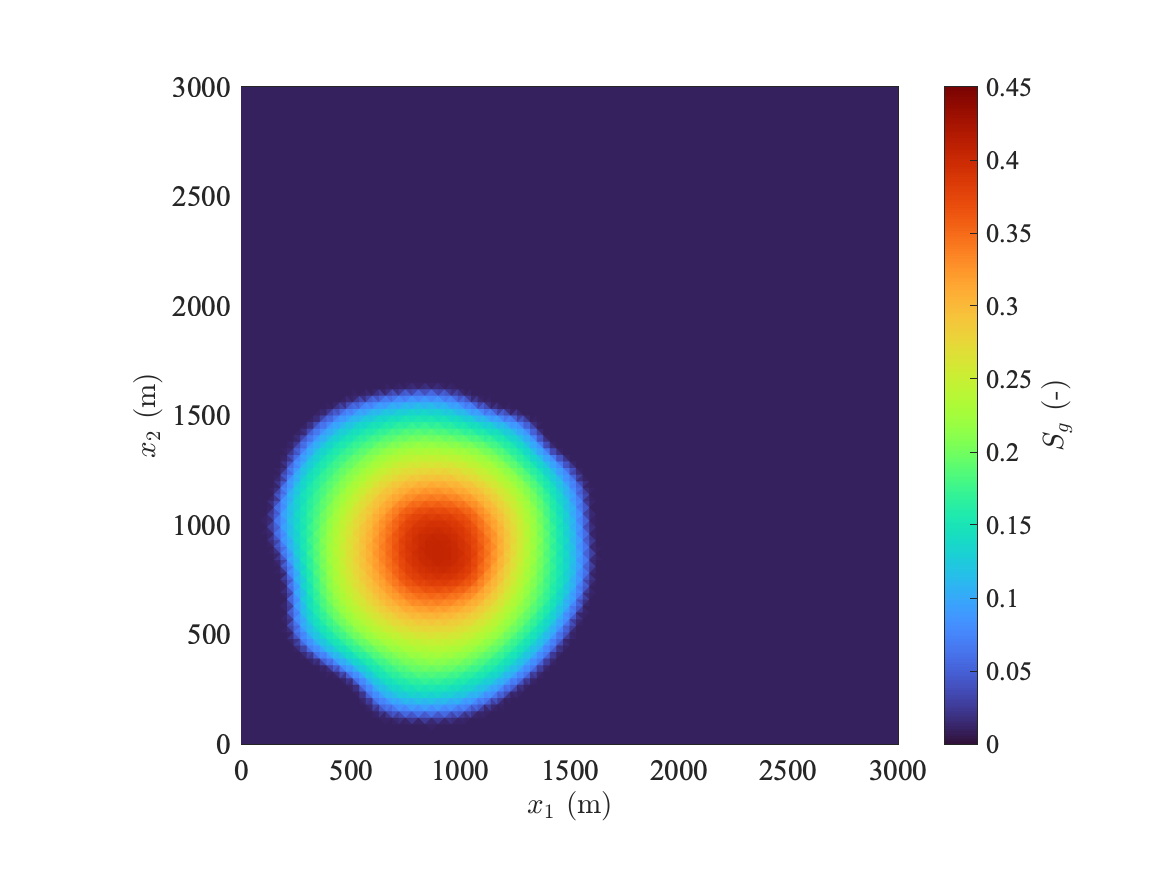}
    }
    \caption{HF results of pressure and saturation. (a)-(c) Pressure; (d)-(f) Saturation.}
    \label{fig:ex3-2 HF PS}
\end{figure}

The first relocated well is placed in
\[
[0.72,1.08]\times[0.72,1.08]~\mathrm{km}^2,
\]
which lies in subdomain 7. Figure~\ref{fig:ex3-2 HF PS} shows the corresponding high-fidelity pressure and saturation fields at \(t=10\), \(20\), and \(30\) years.

\begin{figure}[htbp]
    \centering
    \subfloat[\label{fig:ex3q2 errline 2m b}]{
        \includegraphics[width=0.45\textwidth]{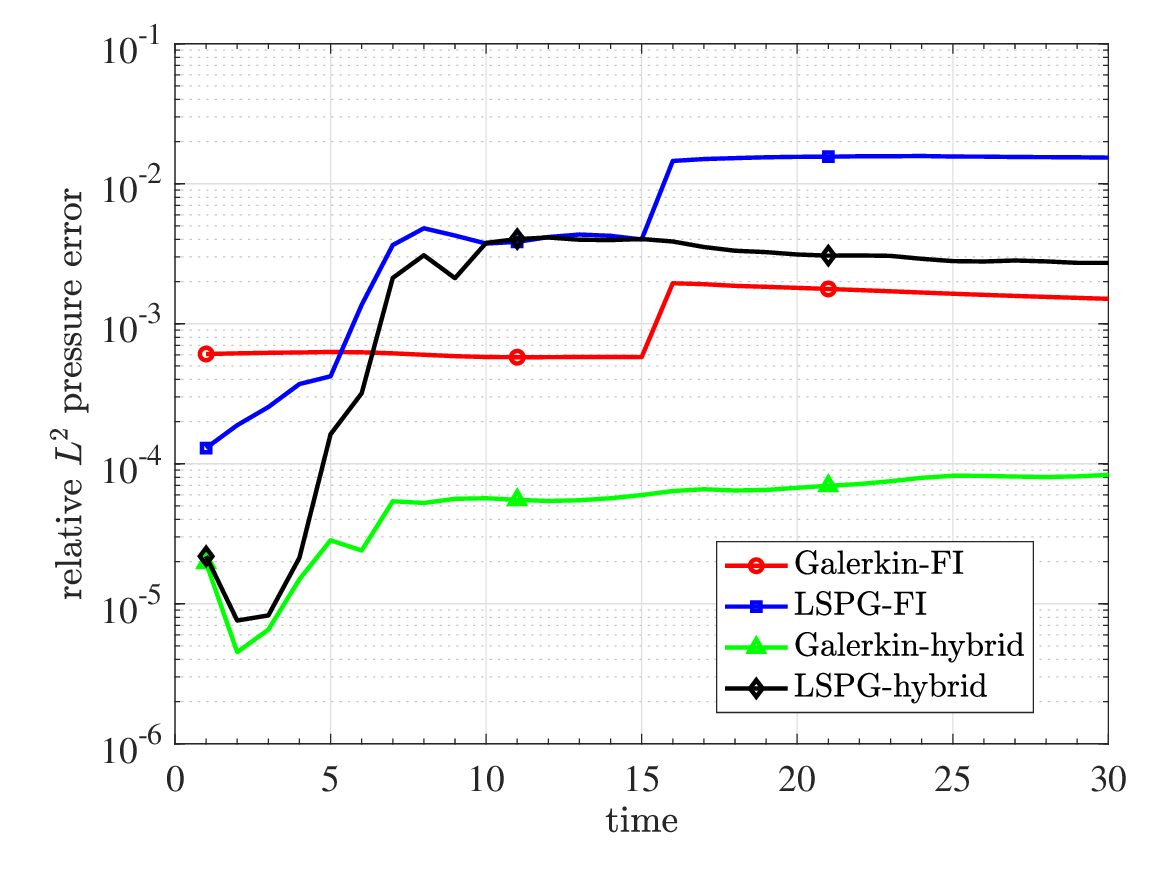}
    }
    \hfill
    \subfloat[\label{fig:ex3q2 errline 2m c}]{
        \includegraphics[width=0.45\textwidth]{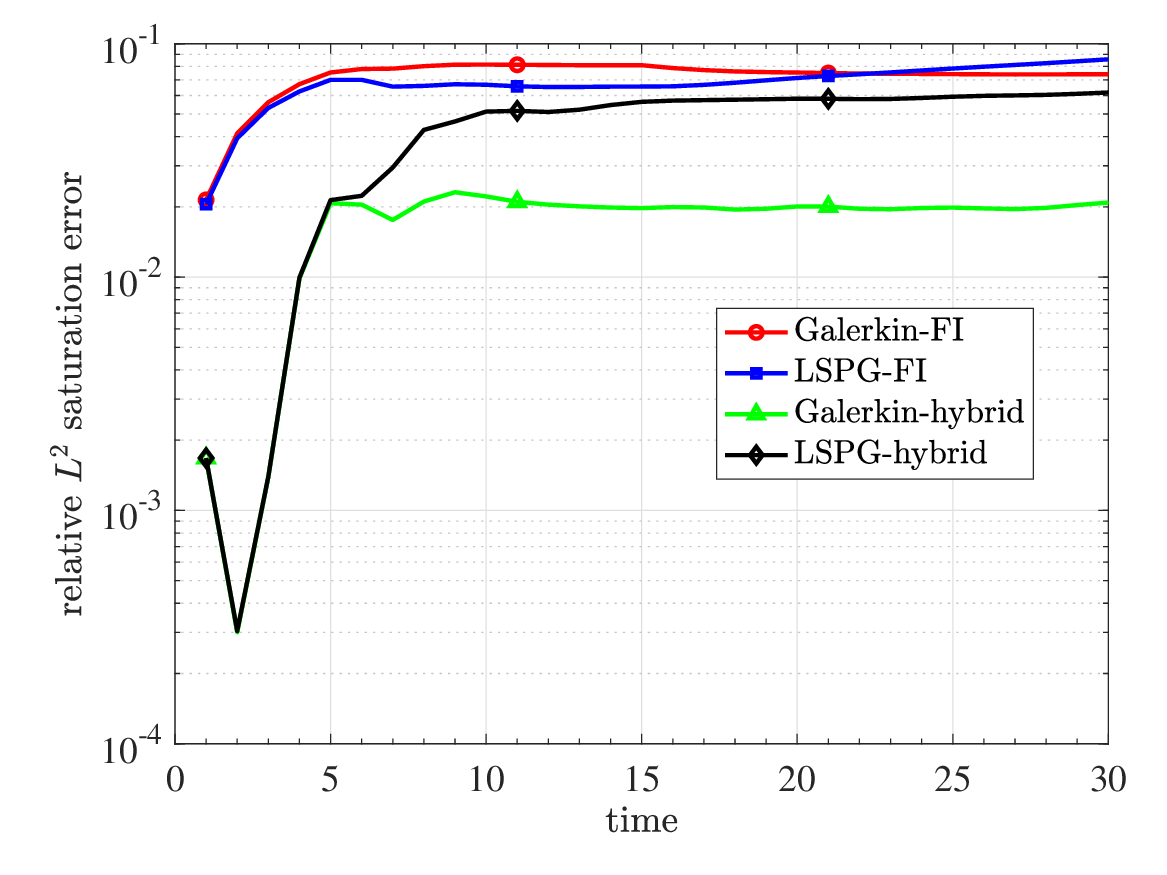}
    }
    \caption{Relative error evolution for HF results and ROM results.}
    \label{fig:ex3q2 errline 2m}
\end{figure}

Figure~\ref{fig:ex3q2 errline 2m} shows the relative errors of Galerkin-FI, LSPG-FI, Galerkin-hybrid and LSPG-hybrid. 

\begin{figure}[H]
    \centering
    \subfloat[$t=10$ years\label{fig:ex3-3 HF PS a}]{
        \includegraphics[width=0.31\textwidth]{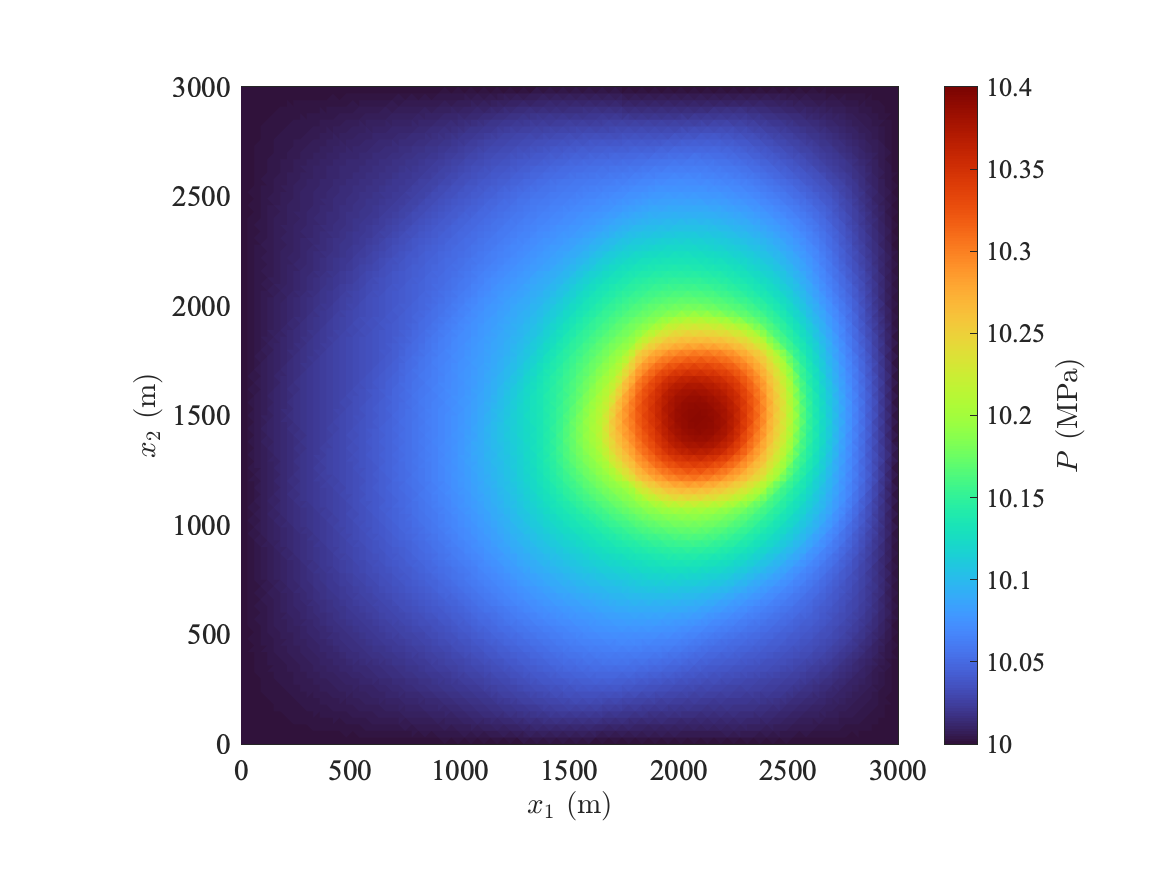}
    }
    \hfill
    \subfloat[$t=20$ years\label{fig:ex3-3 HF PS b}]{
        \includegraphics[width=0.31\textwidth]{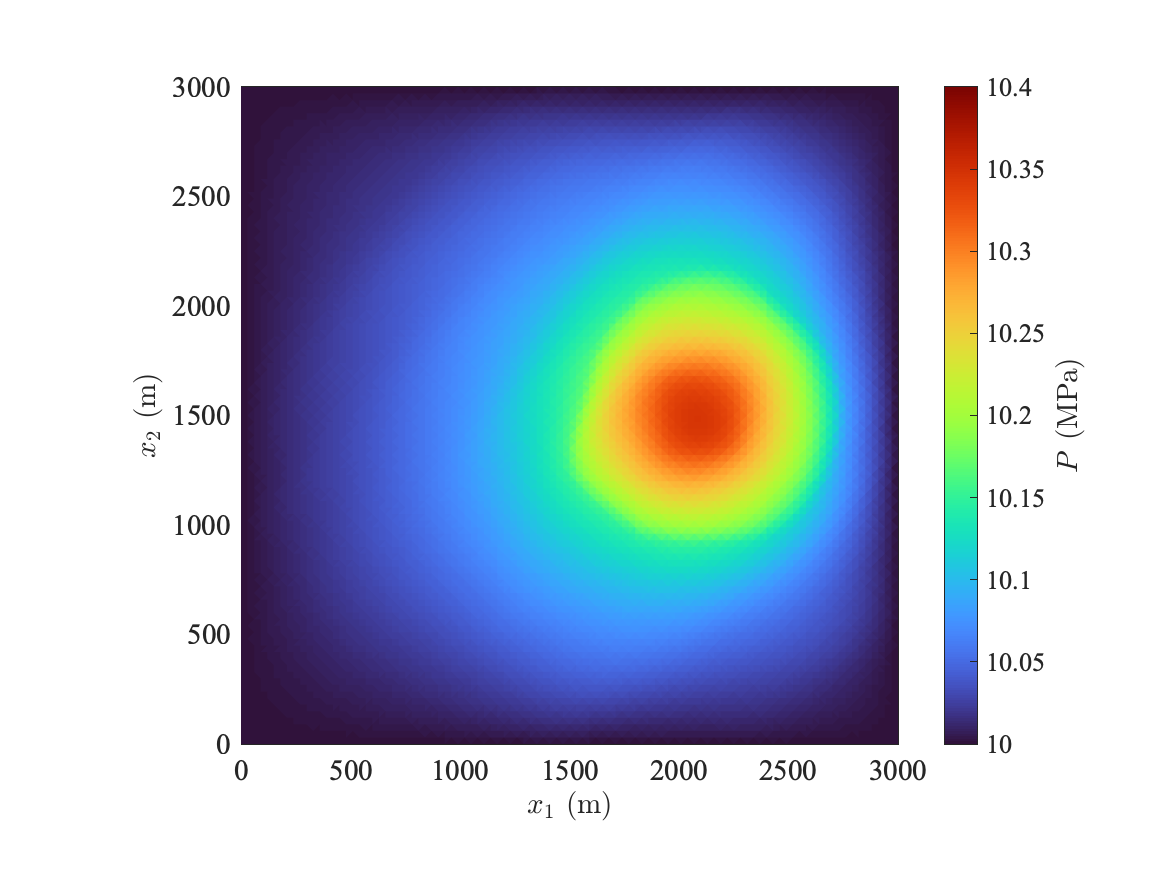}
    }
    \hfill
    \subfloat[$t=30$ years\label{fig:ex3-3 HF PS c}]{
        \includegraphics[width=0.31\textwidth]{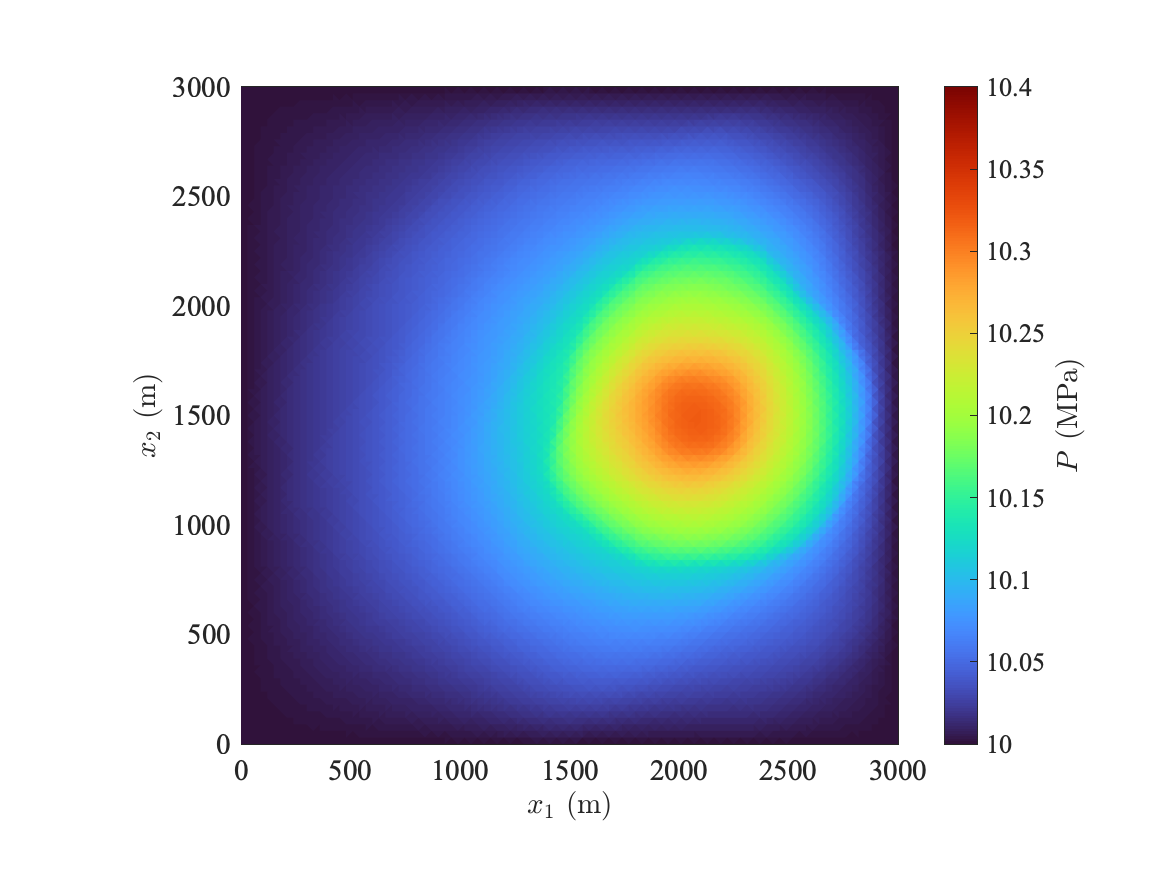}
    }
    \vspace{0.1cm} 
    \subfloat[$t=10$ years \label{fig:ex3-3 HF PS d}]{
        \includegraphics[width=0.31\textwidth]{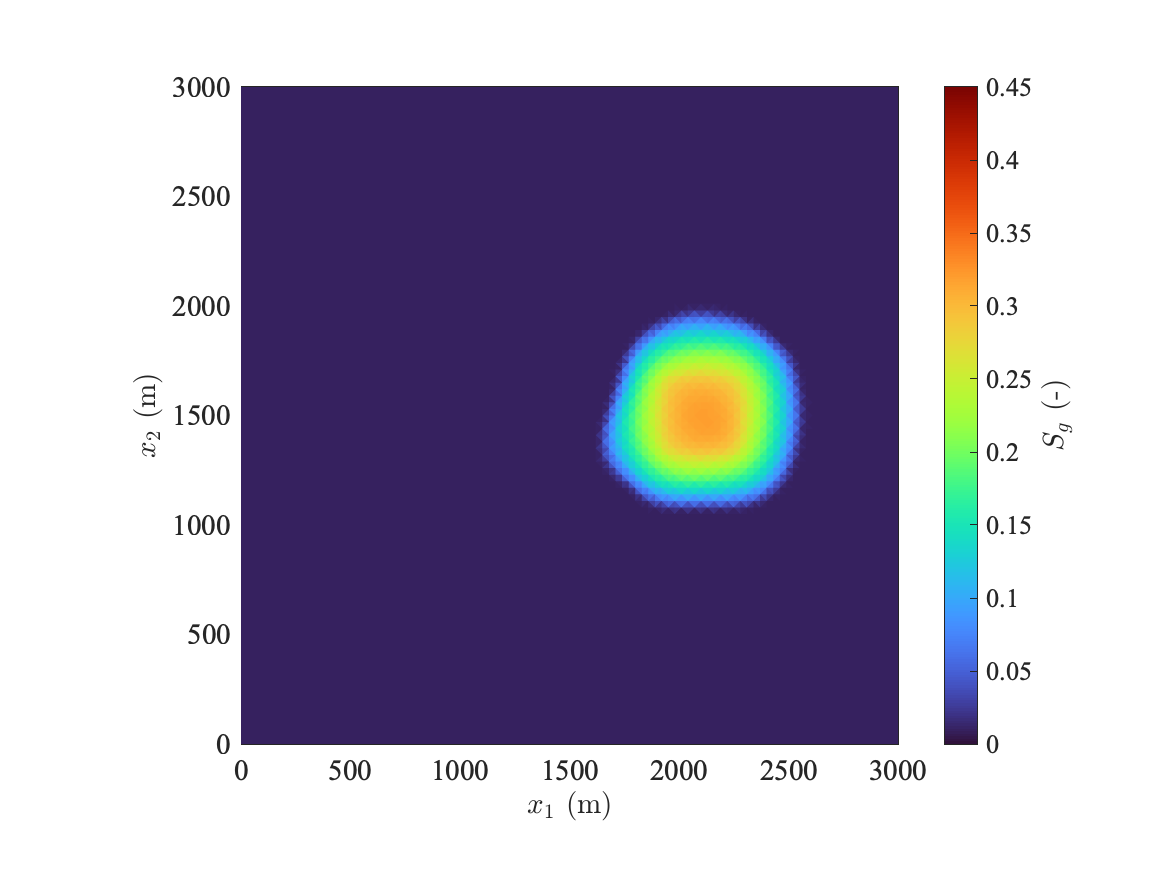}
    }
    \hfill
    \subfloat[$t=20$ years \label{fig:ex3-3 HF PS e}]{
        \includegraphics[width=0.31\textwidth]{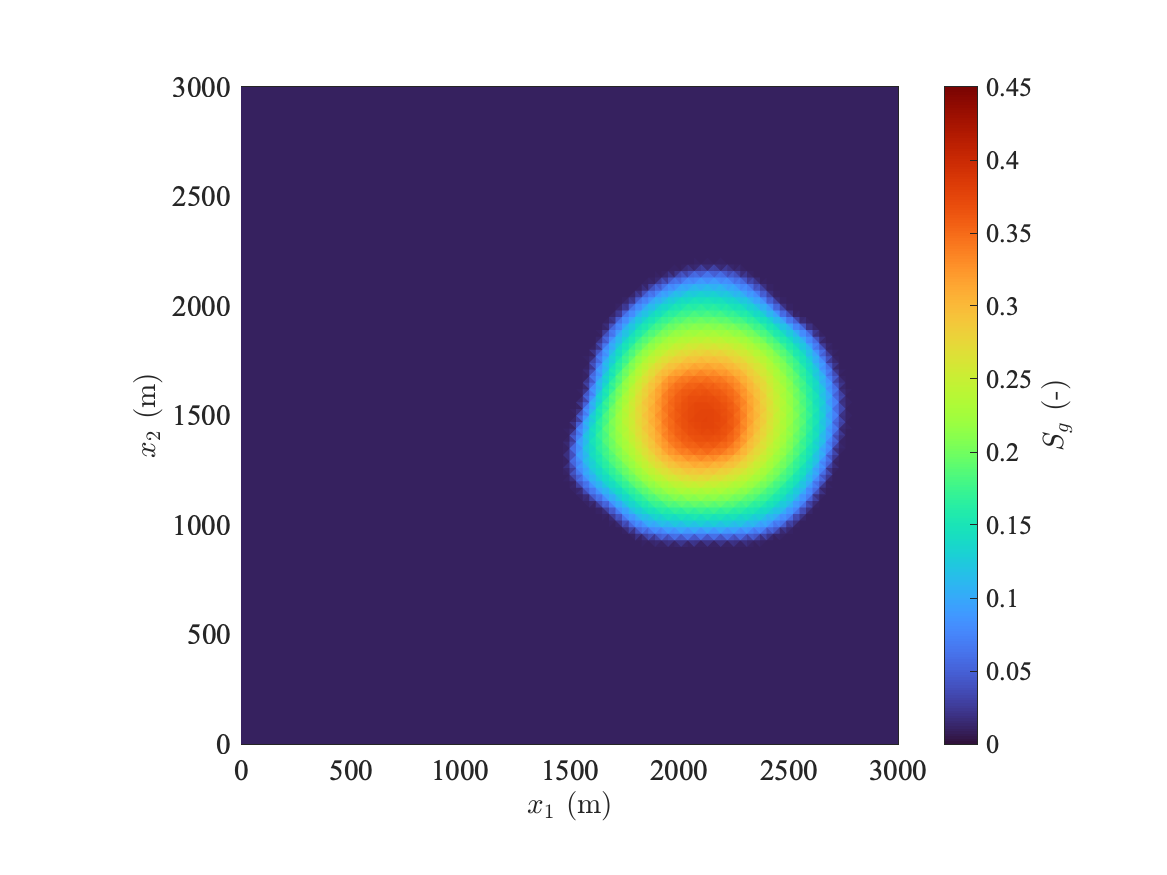}
    }
     \hfill
    \subfloat[$t=30$ years\label{fig:ex3-3 HF PS f}]{
        \includegraphics[width=0.31\textwidth]{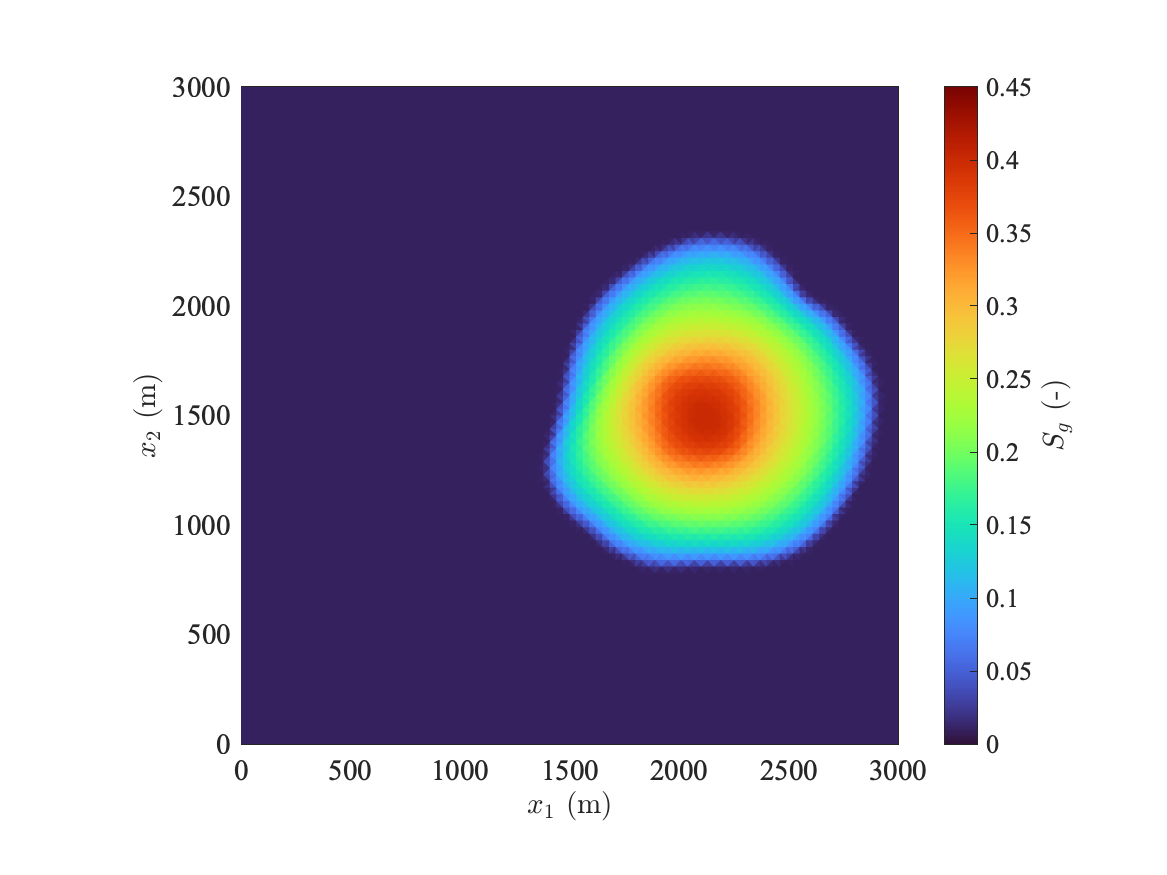}
    }
    \caption{HF results of pressure and saturation. (a)-(c) Pressure; (d)-(f) Saturation.}
    \label{fig:ex3-3 HF PS}
\end{figure}

\begin{figure}[H]
    \centering
    \subfloat[\label{fig:ex3q3 errline 2m b}]{
        \includegraphics[width=0.45\textwidth]{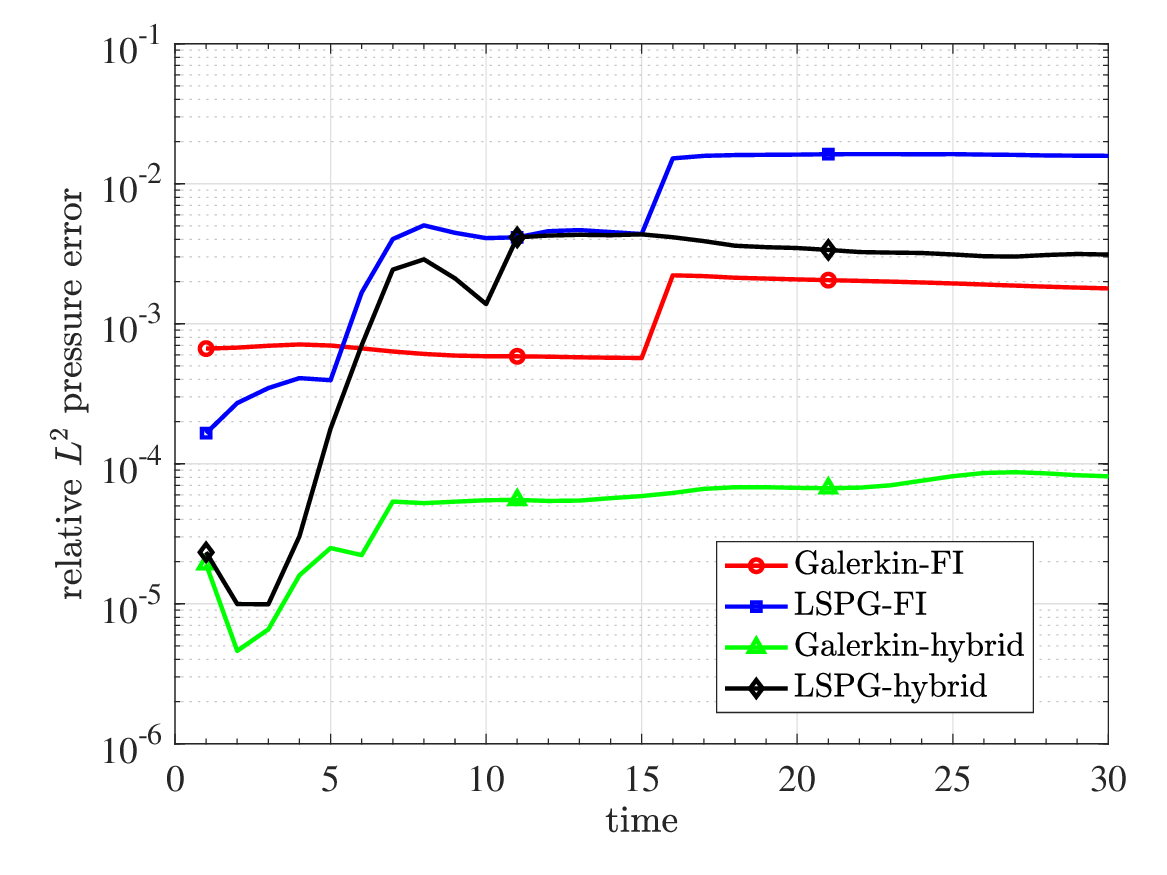}
    }
    \hfill
    \subfloat[\label{fig:ex3q3 errline 2m c}]{
        \includegraphics[width=0.45\textwidth]{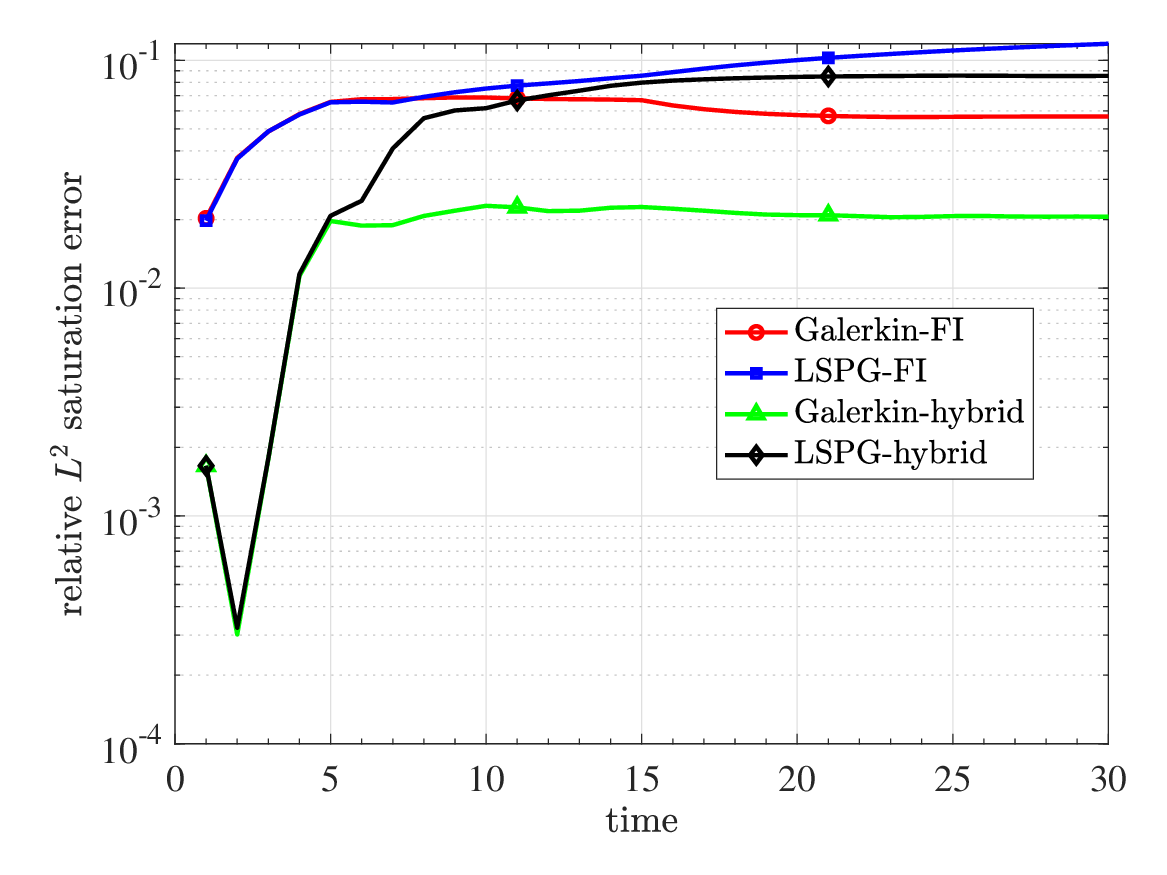}
    }
    \caption{Relative error evolution for HF results and ROM results.}
    \label{fig:ex3q3 errline 2m}
\end{figure}

The second relocated well is placed in
\[
[1.92,2.28]\times[1.32,1.68]~\mathrm{km}^2,
\]
which lies in subdomain 14. 
Figure~\ref{fig:ex3-3 HF PS} shows the HF results. Figure~\ref{fig:ex3q3 errline 2m} shows the relative error. 
In both settings, Galerkin-FI outperforms LSPG-FI in terms of accuracy, and the Galerkin-hybrid variant maintains robust and accurate predictions.

\section{Conclusions}

\label{sec:conclusion}

We developed a domain-decomposition reduced-order modeling  (DD-ROM) framework for two-phase CO$_2$ storage simulations with heterogeneous media and changing well configurations. The central idea is to use the Schur complement structure of the pressure and transport subproblems as the reduced architecture: local subdomain states are approximated by Galerkin or LSPG reduced models, while the global coupling is imposed through pressure and saturation interface residuals. 
Since the interface equations are shared by the high-fidelity and reduced formulations, the model fidelity can be selected locally, for example by using ROMs in the far field and the FOM in well or near-well subdomains, without introducing new transmission conditions. To handle moving wells, the proposed DD-ROM structure is combined with well-aware ingredients, including HAPOD-based far-field basis compression, source-response enrichment, and hybrid FOM--ROM coupling. 

The numerical results validate the main components of the methodology. The manufactured solution test verifies the high-fidelity DD solver and confirms the accuracy of the DD-ROM variants in a smooth setting. The parameterized permeability test shows that the local HAPOD bases provide accurate out-of-sample predictions in heterogeneous media. For moving wells, the results show the benefit of using full-interface coupling and hybrid FOM--ROM treatment near the relocated source, since the displaced near-well response is not always captured reliably by the original reduced spaces.

Future work will focus on adaptive and localized extensions of the present framework. We plan to develop online criteria for basis enrichment and hybridization, so that high-fidelity treatment is activated only in subdomains where the reduced approximation is insufficient, for example near relocated wells or under-resolved saturation fronts. Another direction is to reduce the offline cost by localized training and component-based enrichment \cite{smetana2023,iollo2023}, rather than relying on global high-fidelity simulations. The hyper-reduction strategy also requires further development: in the present tests, the sampled edge and cell sets remain moderately large, which limits the achievable speedup. More adaptive conservative sampling rules are therefore needed to reduce nonlinear assembly cost while protecting interface fluxes, boundary terms, and well/source contributions \cite{ebrahimi2024hyperreduced,ebrahimiYano2026adaptive}. Finally, 
extensions to more detailed compositional CO$_2$ storage models and more realistic well controls will also be pursued.

\bibliography{all_references}
\end{document}